%% file: ComMorse.tex
\documentclass[10pt]{amsart}
\usepackage{amsfonts}
\usepackage{amscd}
\usepackage{amssymb}
\usepackage{amsthm}
\usepackage{amsmath}
\usepackage{pb-diagram}
\usepackage{graphicx}
\usepackage{color}

\def\noi{\noindent}

\def\bu{_\bullet}

\def\xrarrow{\xrightarrow} 

\def\smalldisj{\ {\textstyle{\coprod}}\ }

\def\noteq{\neq}

\def\grad{\nabla}
\def\<{\left<}
\def\>{\right>}

\DeclareMathOperator{\ind}{ind}
\DeclareMathOperator{\simp}{simp}
 \DeclareMathOperator{\Tr}{Tr}
 \DeclareMathOperator{\Hom}{Hom}

\newcommand{\field}[1]{\mathbb{#1}}
\newcommand{\ZZ}{\ensuremath{{\field{Z}}}}
\newcommand{\CC}{\ensuremath{{\field{C}}}}
\newcommand{\RR}{\ensuremath{{\field{R}}}}
\newcommand{\QQ}{\ensuremath{{\field{Q}}}}

\newcommand{\length}[1]{\ensuremath {\| #1 \|}}

\def\ll{\lambda}

\newcommand{\C}{\ensuremath{{\mathcal{C}}}}

\newcommand{\F}{\ensuremath{{\mathcal{F}}}}

\newcommand{\nL}{\ensuremath{{\mathcal{L}}}}
\newcommand{\M}{\ensuremath{{\mathcal{M}}}}

\newcommand{\nP}{\ensuremath{{\mathcal{P}}}}

\newcommand{\R}{\ensuremath{{\mathcal{R}}}}

\newcommand{\Sdot}{\ensuremath{{\mathcal{S}_\bullet}}}

\newcommand{\W}{\ensuremath{{\mathcal{W}}}}
\newcommand{\Wh}{\ensuremath{{\mathcal{W}\textit{h}}}}

\def\a{\alpha}
\def\b{\beta}
\def\g{\gamma}
\def\d{\partial}

\def\e{\epsilon}
\def\f{\phi}
\def\k{\kappa}
\def\r{\rho}
\def\s{\sigma}
\def\Sig{\Sigma}
\def\t{\tau}
\def\th{\theta}

\def\z{\zeta}
\def\w{\omega}

\def\End{{\rm End}}

\def\Ind{{\rm Ind}}
\def\interior{{\rm int}}

\def\colorr{\color{red}}
\def\colorb{\color{blue}}

\def\st{\,|\,}

\def\op{^{op}}

\def\c{\psi}

\newtheorem{thm}{Theorem}[section]
\newtheorem{lem}[thm]{Lemma}
\newtheorem{cor}[thm]{Corollary}
\newtheorem{prop}[thm]{Proposition}
\theoremstyle{definition}
\newtheorem{defn}[thm]{Definition}
\newtheorem{rem}[thm]{Remark}
\newtheorem{eg}[thm]{Example}

\begin{document}

\title[Higher complex torsion and framing principle]{Higher Complex Torsion and \\the Framing Principle}

\author{Kiyoshi Igusa}

\input{abstract}



\subjclass[2000]{Primary 57R45; Secondary 57R50, 19J10}


\keywords{higher Franz-Reidemeister torsion, almost complex
manifolds, diffeomorphisms, twisted cochains, superconnections,
polylogarithms, transfer}

\thanks{Partially supported by the National Science Foundation.}

\maketitle


\input{intro}


\input{CT} 

\input{sec1} 

\input{sec2} 

\input{FramPr} 

\input{ProofFP} 

\input{FPapp} 





\bibliographystyle{amsalpha}

\bibliography{C:/bookbib3}

\end{document}

%% file: abstract.tex
\begin{abstract}
We prove the Framing Principle in full generality and use it to
compute the higher FR-torsion for all smooth bundles with oriented
closed even dimensional manifold fibers. We also show that the
higher complex torsion invariants of bundles with closed almost
complex fibers are multiples of generalized Miller-Morita-Mumford
classes.
\end{abstract}

%% file: intro.tex
\section*{Introduction}

This manuscript is a combination of two papers. The first paper is
about ``higher complex torsion,'' which is a generalization of
higher Franz-Reidemeister (FR) torsion to smooth bundles with
almost complex fibers. The rest of the manuscript is a summary of
basic facts about higher Franz-Reidemeister torsion. Among these
``basic facts'' there are several new result, such as the general
Framing Principle and its applications and the transfer theorem,
but the rest was already explained with full details in
\cite{[I:BookOne]}.

There are several reasons why these two papers appear together as
one manuscript. One is heuristic. The general Framing Principle
allows us to compute the higher FR torsion for smooth bundles with
even dimensional fibers. This is needed to show that the complex
torsion is a generalization of the ``real torsion.'' The second
reason is that the proof of the main theorem about higher complex
torsion (the ``Complex Framing Principle'') is the same as the
proof of the Framing Principle. Finally, the facts about higher FR
torsion are not widely known, so it would be difficult for anyone
to read the article about higher complex torsion without a summary
of the definition and main properties of the higher FR torsion.

\subsection{Higher FR-torsion}
Higher Franz-Reidemeister torsion is an invariant of smooth
bundles
\begin{equation}\label{intro:eq:manifold bundle}
  M\to E\to B
\end{equation}
where $M,E,B$ are compact smooth manifolds. In this introduction
we will concentrate on the case of untwisted coefficients and
assume that $\pi_1B$ acts trivially on the rational homology of
$M$. In other words, our bundle (\ref{intro:eq:manifold bundle})
is classified by a mapping
\[
    \xi_E:B\to BDif\!f_0(M)
\]
of $B$ into the classifying space of the group $Dif\!f_0(M)$ of
diffeomorphisms of $M$ which induce the identity on the rational
homology of $M$.

Here are two examples.
\begin{enumerate}
  \item $Dif\!f_0(S^n)$ is the group of orientation preserving
  diffeomorphisms of $S^n$.
  \item $Dif\!f_0(\Sig_g)\simeq\pi_0(Dif\!f_0(\Sig_g))=T_g$ is the
  \emph{Torelli group} if $\Sig_g$ is a closed oriented surface of genus
  $g$.
\end{enumerate}

Using Morse theory and Waldhausen K-theory, John Klein
\cite{[K:thesis]} and I \cite{[IK1:Borel2]}, \cite{[I:BookOne]}
constructed a cohomology invariant
\[
    \t_{2k}(E)\in H^{4k}(B;\RR)
\]
which is the pull-backs of a universal higher FR-torsion
invariant, i.e., $\t_{2k}(E)=\xi_E^\ast(\t_{2k})$ where
\[
    \t_{2k}\in H^{4k}(BDif\!f_0(M);\RR).
\]

In 1995 Bismut and Lott \cite{[Bismut-Lott95]} constructed similar
classes called \emph{higher analytic torsion} classes. Sebastian
Goette \cite{[Goette01]} explained why these cohomology classes
agree in many cases. Basically, it is because both definitions use
a canonical $A_\infty$ functor associated to a smooth fiber bundle
as we explain later in this introduction.

The homotopy type of the diffeomorphism groups of spheres was
determined by Farrell and Hsiang in the 1970's. Using the higher
FR-torsion their result can be stated as follows.


\begin{thm}[Farrell-Hsiang\cite{[Farrell-Hsiang:Diff(Dn)]}]\label{intro:thm:Farrell-Hsiang}
(a) For even dimensional spheres we have
\[
   H^\ast(BDif\!f_0(S^{2N});\RR)\cong
   \RR[\t_2,\t_4,\cdots]\cong\RR[p_1,p_2,\cdots]
\]
in the stable range (up to degree $\frac{2N}{3}-3$ according to
\cite{[I:BookOne]}).

(b) For odd dimensional spheres we have twice as much cohomology:
\[
   H^\ast(BDif\!f_0(S^{2N+1});\RR)\cong
   \RR[p_1,\t_2,p_2,\t_4,\cdots]
\]
in the stable range.
\end{thm}

\begin{rem}[Newton polynomials]
For even dimensional spheres, the relationship between the higher
torsion classes $\t_{2k}$ and the Pontrjagin classes $p_j$ is
given by
\[
   \t_{2k}=(-1)^k\z(2k+1)\left(\frac1{(2k)!}N_k(p_1,p_2,\cdots)\right)
\]
where $N_k$ is the Newton polynomial which expresses $\sum_i
x_i^k$ in terms of the elementary symmetric functions of $x_i$ and
$\z(s)=\sum_{n\geq1} n^{-s}$ is the Riemann zeta function. For odd
spheres a similar relation holds for the restriction of the
classes $\t_{2k}$ and $p_j$ to $BSO(2N+2)$. (See
Remark~\ref{FP:rem:Newton polynomial}.)
\end{rem}

The higher torsion of the Torelli group was first studied by John
Klein in \cite{[K:Torelli]}. Klein conjectured that the higher
FR-forsion classes of $T_g$ were proportional to the
Miller-Morita-Mumford classes. This was verified by R. Hain, R.
Penner and the author in 1992 but the original argument was never
published. In \cite{[I:BookOne]}, using the Framing Principle, a
precise formula was obtained:
\begin{equation}\label{intro:eq:torsion of Torelli}
 \t_{2k}(T^s_g)=(-1)^k\z(2k+1)\frac{\k_{2k}}{2(2k)!}\in
 H^{4k}(T^s_g;\RR).
\end{equation}
Here $T_g^s$ is the Torelli group of a surface $\Sig_g^s$ of genus
$g$ with $s\geq1$ marked points.

One of the purposes of this paper is to extend this formula to odd
Miller-Morita-Mumford classes.

\subsection{Construction of $\t_k$}


The basic construction behind the definition of higher
Franz-Reidemeister torsion goes back to early work of Ed Brown
\cite{[Brown59:TwistedTensor]} where he defined the concept of a
\emph{twisted cochain} $\f$. These can be written as sums of
cochains $\f=\sum_{p\geq0}\f_p$ defined on the base space $B$ of a
fiber bundle $E\to B$ so that
\[
    \delta\f=\f'\cup\f
\]
where $\f'=\sum(-1)^p\f_p$. Here $\f_p$ is a $p$-cochain with
coefficients in the degree $p-1$ part of a graded endomorphism
ring $\End_R(PR)$ where $P$ is a partially ordered graded set. If
we add the identity endomorphism of the free $R$-module $RP$ to
$\f_1$ we get an $A_\infty$-functor given on morphisms by $
\Phi_p=\f_p$ for $p\noteq1$, $\Phi_1=I+\f_1$ and on objects by
$\Phi(X)=(RP,\f_0(X))$. This is a contravariant $A_\infty$ functor
from the category of generic small simplices in $B$ to the
differential graded category of finitely generated free
$R$-complexes. (Remark~\ref{sec1:rem:A infty functors}.)

Bismut and Lott use flat $\ZZ$-graded \emph{superconnections} of
total degree $1$. These are linear operators which in local
coordinates can be written as $D=d+A$ where $A=\sum_{p\geq0}A_p$
and $A_p$ is a $p$-form on $B$ with coefficients in the degree
$1-p$ part of a graded endomorphism ring $\End(E)$ where $E$ is a
graded vector bundle over $B$. The superconnection $D$ is
\emph{flat} if it satisfies the condition
\[
    -dA=A\wedge A.
\]
In other words, $A$ is an infinitesmal twisted cochain up to sign.
An $A_\infty$-functor $\Phi$ is given on objects by
$\Phi(\s)=(E_b,A_0(b))$ where $b$ is the barycenter of
$\s\subseteq B$. To obtain $\Phi_p$ (up to sign) on the morphisms
we need to integrate $A_p$ for $p\noteq1$ and $I-A_1$ for $p=1$.
This process is independent of the choice of coordinates. For
example, $\Phi_1$ is parallel transport by the connection $d+A_1$
along the edges of the first barycentric subdivision of $B$. (See
\cite{[Goette01]}, \cite{[I:Twisted]} for more details.)

Brown showed that the homology of the total space of a fiber
bundle $E\to B$ could be given by the \emph{twisted tensor
product}
\begin{equation*}
  C_\ast(B)\otimes_\f C_\ast(F)
\end{equation*}
which is the ordinary tensor product complex with boundary map
$\d_\f$ which is twisted by a twisted cochain $\f$. If we take
coefficients in $\QQ$ (or any other field) we can replace
$C_\ast(F)$ by the graded homology $H_\ast(F;\QQ)$ which has an
$A_\infty$ functor ${\Psi}$ induced by $\Phi=I+\f$. If $\pi_1B$
acts trivially on $H_\ast(F;\QQ)$, there is an induced twisted
cochain ${\c}$ so that ${\Psi}=I+{\c}$ where $I$ is the identity
map on $H_\ast(F;\QQ)$. Then we get another twisted tensor product
\begin{equation}\label{intro:eq:algebraic twisted tensor product}
    C_\ast(B;\QQ)\otimes_{{\c}}C_\ast(F;\QQ)\simeq
    C_\ast(E;\QQ).
\end{equation}

In the case when $E\to B$ is a smooth bundle with compact manifold
fiber $M$, we can find a fiberwise ``generalized Morse function''
$f:E\to\RR$ which gives a finitely generated cellular chain
complex for the fiber over each generic small transverse simplex
in $B$. This Morse theory construction defines another twisted
cochain and another twisted tensor product
\begin{equation}\label{intro:eq:Morse theory twisted tensor product}
    C_\ast(B)\otimes_\f C_\ast(f)
\end{equation}
which is rationally fiberwise quasi-isomorphic to
(\ref{intro:eq:algebraic twisted tensor product}). In other words,
they are equivalent as $A_\infty$-functors. However, they are not
simplicially equivalent as twisted cochains on $B$. There is a
well-defined algebraic K-theory obstruction to deforming one
twisted cochain into the other. The higher Franz-Reidemeister
torsion is a numerical invariant which depends on this K-theory
class.

This numerical invariant is given by first interpolating the
twisted cochain so that it becomes an infinitesmal twisted
cochain, i.e., a flat $\ZZ$ graded superconnection and integrating
an even degree form derived from the odd degree Kamber-Tondeur
form. This process, explained in detail in \cite{[I:BookOne]}, is
summarized in subsection \ref{sec1:subsec:superconnections}.

\subsection{Framing Principle}

The Framing Principle is the main tool used to compute the higher
Franz-Reidemeister torsion. The original version of the Framing
Principle in \cite{[I:BookOne]} has a restriction on the
birth-death singularities of a family of oriented generalized
Morse functions. In this paper we remove this condition and prove
the Framing Principle for all families of oriented generalized
Morse functions.

The general Framing Principle gives us new insight into the nature
of higher Franz-Reidemeister torsion. The higher FR torsion is a
sum of two invariants: one is a \emph{tangential invariant} which
depends only on the fiber homotopy type of the vertical tangent
bundle and the other is an \emph{``exotic'' diffeomorphism
invariant} which detects different smooth structures on the same
fiberwise tangential homeomorphism type.

The tangential invariant $T_{2k}(E)$ is similar in description to
the \emph{Miller-Morita-Mumford classes} for oriented surface
bundles. This invariant is equal to the higher FR torsion for all
oriented closed even dimensional fibers and equal to zero for
oriented closed odd dimensional fibers. Equivalently, the exotic
component is zero in the first case and equal to the entire higher
torsion invariant in the second case.

For oriented even dimensional fibers with boundary, the higher FR
torsion is given by the tangential invariant plus half the torsion
of the boundary:
\begin{equation}\label{intro:eq:exotic and tangential}
    \t_{2k}(E)=\underbrace{\frac12\t_{2k}(\d E)}_{exotic\ component}+
    \underbrace{\frac12
    (-1)^k\z(2k+1)\frac1{(2k)!}T_{2k}(E)}_{tangential\ component}.
\end{equation}
For oriented odd dimensional fibers with boundary the tangential
invariant comes entirely from the boundary.

These facts can be summarized by saying that the exotic torsion
comes entirely from oriented odd dimensional fibers in all cases.
(For oriented even dimensional fibers the exotic torsion is merely
detecting the exotic smooth structure of the boundary. For
unoriented fibers we go to either the oriented cover or the
oriented disk bundle, whichever is odd dimensional.)

\subsection{Complex torsion} Finally we come the the main point of
this paper which is to give a generalization of higher FR torsion
for bundles with almost complex (AC) fibers. One of the main
motivations is to extend the Klein-Hain-I-Penner relation
(\ref{intro:eq:torsion of Torelli}) to all of the
Miller-Morita-Mumford classes (not just the even ones).

The complex torsion of a bundle $E\to B$ with AC fiber $M$ is a
sequence of cohomology classes $\t_k^\CC(E,\z)_m$ defined using
the equivariant higher FR torsion of the vertical tangent lens
space bundle $S^{2n-1}(E)/Z_m$. It is related to the
\emph{generalized Miller-Morita-Mumford classes}
\[
    T_k(E):=tr_B^E(k!ch_k(T^vE))\in H^{2k}(B;\ZZ)
\]
by the following formula assuming that $M$ is a closed AC manifold
\begin{equation}\label{intro:eq:complex torsion and generalized
MMM}
   \t_k^\CC(E,\z)_m=\frac12m^k L_{k+1}(\z)
    \frac1{k!}T_k(E)
\end{equation}
where $L_{k+1}(\z)$ is the polylogarithm function
\[
    L_{k+1}(\z)
    =\R\left(\frac1{i^k}\sum_{n=1}^\infty\frac{\z^n}{n^{k+1}}\right).
\]
As a special case we obtain the complex torsion of $M_g$, the
mapping class group of a surface of genus $g$:
\begin{equation}\label{intro:eq:complex torsion of Mg}
    \t_k^\CC(M_g,\z)_m=\frac12m^kL_{k+1}(\z)\frac{\k_k}{k!}.
\end{equation}
This is the generalization of equation (\ref{intro:eq:torsion of
Torelli}) that we were looking for.

Both of the main theorems of this paper, the general Framing
Principle and the computation of higher complex torsion, are
examples of the \emph{transfer theorem} for higher FR torsion.
This theorem says that the higher FR torsion commutes with
transfer in the sense that, if $D\to E$ is a bundle with fiber $Y$
and $\F$ is a Hermitian coefficient system on $D$ for which $Y$ is
acyclic, then the FR torsion of $D$ as a bundle over $B$ is the
transfer of the FR torsion of $D$ as a bundle over $E$:
\begin{equation}\label{intro:eq:transfer theorem}
        \t_k(D;\F)_B=tr^E_B(\t_k(D;\F)_E).
\end{equation}
The transfer theorem and the explanation of the Framing Principle
in terms of the transfer were inspired by the work of Dwyer, Weiss
and Williams \cite{[DWW]}. Analogous formulas for higher analytic
torsion were given by Xiaonan Ma \cite{[Ma97]}.

We note that the transfer formula also explains the presence of
polylogarithms in the formula for higher FR torsion. (See
subsection \ref{sec2:subsec:transfer}.)

Outline:
\begin{enumerate}
    \item Complex torsion
\begin{enumerate}
  \item Definition for closed AC fibers
  \item Generalized Miller-Morita-Mumford classes
  \item Complex Framing Principle
  \item Nonempty boundary case
\end{enumerate}
    \item Definition of higher FR torsion
\begin{enumerate}
  \item Generalized Morse functions
  \item Families of chain complexes
  \item Monomial functors
  \item Filtered chain complexes
  \item Subfunctors
  \item The Whitehead category
  \item Definition of higher FR torsion
  \item {Families of matrices as flat superconnections}
  \item Independence of birth-death points
  \item Positive suspension lemma
\end{enumerate}
  \item Basic properties of higher FR-torsion
\begin{enumerate}
   \item Basic properties
   \item Suspension theorem
  \item Additivity, Splitting Lemma
  \item Applications of the Splitting Lemma
  \item Local equivalence lemma
  \item Product formula
  \item Transfer for coverings
  \item More transfer formulas
\end{enumerate}
  \item The Framing Principle
\begin{enumerate}
  \item Statement for Morse bundles
  \item General statement
  \item Push-down/transfer
  \item The Framing Principle
\end{enumerate}
    \item Proof of the Framing Principle
\begin{enumerate}
  \item Transfer theorem
  \item Stratified deformation lemma
  \item Proof of transfer theorem
  \item Proof of Framing Principle
\end{enumerate}
    \item Applications of the Framing Principle
\begin{enumerate}
  \item Torelli group
  \item Even dimensional fibers
  \item Unoriented fibers
  \item Vertical normal disk bundle
\end{enumerate}
\end{enumerate}

I would like to thank Bernhard Keller who directed me to the work
of Kadeishvilli \cite{[Kad80]} and thereby also alerted me to the
earlier work of my colleague Ed Brown
\cite{[Brown59:TwistedTensor]}. I should also thank Sebastian
Goette for explaining his work on higher analytic torsion
\cite{[BG2]}, \cite{[Goette01]}. It was in correspondences with
Goette that I realized the importance of proving the Framing
Principle in full generality. Finally, I would like to thank Ed
Brown and Gordana Todorov who helped me to sort out my ideas about
twisted cochains and $A_\infty$-functors.

%% file: CT.tex
\vfill\eject\section{Complex torsion}\label{CT:section}

\begin{enumerate}
  \item Definition for closed AC fibers
  \item Generalized Miller-Morita-Mumford classes
  \item Complex Framing Principle
  \item Nonempty boundary case
\end{enumerate}

In this section we consider smooth bundles
\begin{equation}\label{CT:eq:bundle with AC fibers}
  M^{2n}\to E\to B
\end{equation}
where the fibers are {almost complex (AC) manifolds}, i.e., the
vertical tangent bundle $T^vE$ has an almost complex structure
$J$. We call this a \emph{bundle with almost complex fibers}. We
assume that $B$ is a connected closed manifold but $M$ is only
assumed to be compact. We consider separately the cases when $M$
is closed and when it has a boundary.

We will define ``higher complex torsion'' invariants
\[
    \t_k^\CC(E,\z)_m\in H^{2k}(B;\RR)
\]
having properties analogous to the higher FR torsion invariants
for smooth bundles with even dimensional manifold fibers. The
definition is a complexified version of Theorems
\ref{FP:thm:normal disk bundle-closed case} and \ref{FP:thm:normal
disk bundle with boundary}.

\subsection{Definition for closed AC fibers}\label{CT:subsec:definition}

Although the almost complex fibers of $p:E \to B$ are
diffeomorphic they do not necessarily have corresponding almost
complex structures (i.e., the diffeomorphism does not preserve
$J$). Thus bundles (\ref{CT:eq:bundle with AC fibers}) are not
classified by the group of diffeomorphisms of $M$ which preserve
its almost complex structure. We use instead a standard simplicial
construction for the classifying space.

\begin{defn}[the classifying space $\M_\CC(M)$]\label{CT:def:classifying space for AC bundles}
If $M^{2n}$ is an even dimensional real manifold let $\M_\CC(M)$
be the geometric realization of the simplicial set whose
$k$-simplices are manifolds $E$ diffeomorphic to $M\times\Delta^k$
together with an almost complex structure on the vertical tangent
bundle of $E$ over $\Delta^k$.
\end{defn}

There is a canonical bundle
\begin{equation}\label{CT:eq:universal bundle over MC(M)}
  M\to E_\CC(M)\to \M_\CC(M)
\end{equation}
which has an almost complex structure on its vertical tangent
bundle.

\begin{prop}[$\M_\CC(M)$ as classifying space]\label{CT:prop:MC(M) classifies AC bundles}
Any smooth bundle $M\to E\to B$ with almost complex fibers is
classified by a map $B\to\M_\CC(M)$ which is unique up to
homotopy.
\end{prop}

\begin{proof}
Any smooth triangulation of $B$ gives a map $B\to\M_\CC(M)$. Any
smooth triangulation of $B\times\d I$ extends to a triangulation
of $B\times I$ giving a homotopy between any two classifying maps.
\end{proof}

Since all complex bundle admit Hermitian metrics which are unique
up to isotopy we can construct the vertical tangent sphere bundle:
\[
    S^{2n-1}\to S^{2n-1}(E)\to E.
\]
The circle group $S^1=U(1)$ acts fiberwise on the sphere bundle
$S^{2n-1}(E)$ so we get a free action of every cyclic group
\[
    Z_m=\{\z\in\CC\st \z^m=1\}\subseteq S^1=U(1).
\]
The quotient is a bundle over $B$:
\begin{equation}\label{CT:eq:vertical tangent lens space bundle}
  S^{2n-1}(M)/Z_m\to S^{2n-1}(E)/Z_m\to B.
\end{equation}
For any $\z\in Z_m$ let $\r_\z$ be the one dimensional
representation of $\pi_1(S^{2n-1}(M)/Z_m)$ given by composing the
canonical epimorphism
\[
    \pi_1(S^{2n-1}(M)/Z_m)\twoheadrightarrow Z_m
\]
with the homomorphism $Z_m\to Z_m$ which sends the generator
$e^{2\pi i/m}$ to $\z$. Then
\[
    H_\ast(S^{2n-1}(M)/Z_m;\r_\z)=0
\]
for any $\z\noteq1$ in $Z_m$. Consequently, the higher torsion of
(\ref{CT:eq:vertical tangent lens space bundle}) with coefficients
in $\r_z$ are defined for $z\noteq1$.

\begin{defn}[higher complex torsion for closed fibers]\label{CT:def:complex torsion}
We define the \emph{higher complex torsion} of a bundle with
closed almost complex fibers (\ref{CT:eq:universal bundle over
MC(M)}) to be the family of cohomology classes
$\t_k^\CC(E,\z)_m\in H^{2k}(B;\RR)$ for all $k,m\geq1$ and
nontrivial $m$-th roots of unity $\z$ given by
\begin{equation}\label{CT:eq:def of complex torsion}
    \t_k^\CC(E,\z)_m:=-\frac12\t_k(S^{2n-1}(E)/Z_m;\r_\z).
\end{equation}
For the universal bundle (\ref{CT:eq:universal bundle over MC(M)})
for a closed $2n$-manifold $M$ this gives the \emph{universal
higher complex torsion classes}
\[
  \t_k^\CC(\z)_m\in H^{2k}(\M_\CC(M);\RR).
\]
\end{defn}

If the action of $\pi_1B$ on the rational homology of $M$ is upper
triangular then so is its action on $S^{2n-1}(M)/Z_m$ by the Gysin
sequence:
\[
    \cdots\to H_{p+1}(M;\QQ)\to H_{p+2n}(S^{2n-1}(M)/Z_m;\QQ)\to
    H_{p+2n}(M;\QQ)\to\cdots.
\]
So, $\t_k(E,1)_m$ can be defined by (\ref{CT:eq:def of complex
torsion}) in that case.

These invariants are related by the following transfer formula
which follows from the transfer formula for cyclic groups
(\ref{sec2:cor:transfer formula for cyclic groups}).
\begin{prop}[transfer formula for complex
torsion]\label{CT:prop:transfer formula} If $z$ is a nontrivial
$m$-th root of unity then
\begin{equation}\label{CT:eq:transfer formula for complex torsion}
    \t_k^\CC(E,z)_m=\sum_{\z^r=z}\t_k^\CC(E,\z)_{mr}.
\end{equation}
This also holds for $z=1$ if $H_\ast(M;\QQ)$ is $\pi_1B$-upper
triangular.
\end{prop}

The formula for higher complex torsion (\ref{CT:eq:def of complex
torsion}) is a complexified version of Theorem~\ref{FP:thm:normal
disk bundle-closed case} where, instead of taking the normal
bundle, we take the tangent bundle and changed the sign of the
torsion. The following theorem shows that higher complex torsion
is a generalization of higher FR torsion.

\begin{thm}[complex torsion generalizes real torsion]\label{CT:thm:complex torsion generalizes real torsion}
Suppose that the action of $\pi_1B$ on $H_\ast(M;\QQ)$ is upper
triangular. Then
\[
    \t_{2k}^\CC(E,1)_1=\t_{2k}(E)
\]
for all $k\geq1$.
\end{thm}

This follows from Theorem~\ref{FP:thm:normal disk bundle-closed
case} and the following lemma.

\begin{lem}[torsion of complementary sphere bundles]\label{CT:lem:tangent and normal sphere torsion}
Let $S(E)$ be the vertical normal sphere bundle of $E$ (as in
Theorem~\ref{FP:thm:normal disk bundle-closed case}) and suppose
that the action of $\pi_1B$ on $H_\ast(M;\QQ)$ is upper
triangular. Then
\[
    \t_{2k}(S(E))+\t_{2k}(S^{2n-1}(E))=0.
\]
\end{lem}

\begin{rem}[true definition of higher complex torsion]\label{CT:rem:true
def of CT using normal sphere bundle} It would not be necessary to
invoke this lemma if we used the sphere bundle of the complex
vertical normal bundle of $E$ instead of that of the vertical
tangent bundle as in Definition~\ref{CT:def:complex torsion}. In
other words, it makes more sense to use the following equivalent
definition.
\[
        \t_k^\CC(E,\z)_m=\frac12\t_k(S(E)/Z_m;\r_\z).
\]
\end{rem}

\begin{proof}
The Whitney sum of the vertical tangent bundle and the vertical
normal bundle is trivial of real dimension $2N$. This implies that
the trivial sphere bundle $E\times S^{2N-1}$ over $E$ is the
fiberwise join of the tangent sphere bundle $S^{2n-1}(E)$ and the
normal sphere bundle $S(E)$. Consequently, we can take any smooth
functions
\[
    f:S^{2n-1}(E)\to[0,1],\quad g:S(E)\to[2,3]
\]
and extend them to a smooth function
\[
    h:E\times S^{2N-1}\to[0,3]
\]
which looks locally like $f(x)+t^2$ near $S^{2n-1}(E)$ and
$g(y)-t^2$ near $S(E)$ ($t$ being the join parameter) and is
nonsingular everywhere else. This makes the cellular chain complex
$C_\ast(h)$ of $h$ into an extension of $C_\ast(f)$ by
$\Sig^{2n}C_\ast(g)$. Then by the Splitting Lemma, we have
\[
    \t_{2k}(S(E))+\t_{2k}(S^{2n-1}(E))=\t_{2k}(E\times S^{2N-1}),
\]
where
\[
    \t_{2k}(E\times S^{2N-1})=\chi(S^{2N-1})\t_{2k}(E)=0.
\]by the product formula (Lemma~\ref{sec2:lem:product formula}).
\end{proof}

\subsection{Generalized Miller-Morita-Mumford classes}\label{CT:subsec:generalized
MMM}

To state the theorem which computes the higher complex torsion for
closed fibers we need the complex version of
Definition~\ref{FP:def:general real MMM class}. This definition
arose out of a conversation with S. Morita.

\begin{defn}[generalized MMM classes]\label{CT:def:generalized complex MMM classes}
Given a smooth bundle $p:E\to B$ with almost complex fibers let
\[
    T_k(E)=tr^E_B\left(k!\,ch_k(T^vE)\right)\in
H^{2k}(B;\ZZ).
\]
\end{defn}

\begin{prop}[$T_k$ generalized MMM classes]\label{CT:prop:properties of Tk(E)}
{\rm a)} When $k$ is even, Definition~\ref{CT:def:generalized
complex MMM classes} agrees with Definition~\ref{FP:def:general
real MMM class}.

{\rm b)} When $M$ is a closed oriented surface, $T_k(E)$ is equal
to the Miller-Morita-Mumford class
  \[
    T_k(E)=\k_k(E).
  \]

{\rm c)} When $M$ is an oriented surface with boundary, $T_k(E)$
is equal to the punctured Miller-Morita-Mumford class
  \[
    T_k(E)=\k_k^\d(E).
  \]
\end{prop}

\begin{proof} By equation (\ref{FP:eq:T2k(E)=k2k(E)}).
\end{proof}

Because of Proposition~\ref{CT:prop:properties of Tk(E)} (b) we
call $T_k(E)$ the \emph{generalized Miller-Morita-Mumford classes}
of $E$. Note that these invariants are high dimensional stable
invariants.

\begin{prop}[stability of $T_k$]\label{CT:prop:stability of Tk}
\[
    T_k(E)=T_k(E\times D^{2N}).
\]
\end{prop}

Using these generalized MMM classes $T_k$ we have the following
formula for the higher complex torsion in the case of closed AC
fibers.

\begin{thm}[complex torsion for closed AC fibers]\label{CT:thm:complex torsion for closed AC fibers}
Suppose that $M^{2n}\to E\to B$ is a smooth bundle with closed
almost complex fibers. Then
\[
   \t_k^\CC(E,\z)_m=\frac12m^k L_{k+1}(\z)
    \frac1{k!}T_k(E)
\]
where $L_{k+1}(\z)$ is the polylogarithm function
\[
L_{k+1}(\z)=\R\left(\frac1{i^k}\sum_{n=1}^\infty\frac{\z^n}{n^{k+1}}\right).
\]
\end{thm}

\begin{rem}[polylogarithm term is forced]\label{CT:rem:polylogs come from transfer formula}
As we explain in section \ref{sec2:subsec:transfer}, the
polylogarithm term is forced by the transfer formula
(\ref{CT:eq:transfer formula for complex torsion}).
\end{rem}

An important special case is the higher complex torsion of the
\emph{mapping class group} $M_g$ (the group of isotopy classes of
orientation preserving self-diffeomorphisms of a surface $\Sig_g$
of genus $g$).

\begin{cor}[MMM class as complex torsion]\label{CT:cor:MMM classes as higher complex torsion}
The higher complex torsion of the canonical surface bundle over
the classifying space $BM_g$ of the mapping class group $M_g$ is
\[
    \t_k^\CC(M_g,\z)_m=\frac12m^kL_{k+1}(\z)\frac{\k_k}{k!}.
\]
\end{cor}

\begin{proof}[Proof of Theorem~\ref{CT:thm:complex torsion for closed AC
fibers}] Let
\[
    D=S^{2n-1}(E)/Z_m.
\]
This is a lens space bundle over $E$. So, by the formula for the
higher FR-torsion of a lens space bundle (\ref{sec2:cor:torsion of
lens space bundles}) we have
\[
    \t_k(D,\r_\z)_E=-m^kL_{k+1}(\z)ch_k(T^vE).
\]

By the transfer theorem (\ref{Proof:thm:transfer theorem}), this
is related to the torsion of $D$ over $B$ by
\[
    \t_k(D;\r_\z)_B=tr^E_B(\t_k(D;\r_\z)_E)=-m^kL_{k+1}(\z)\frac1{k!}T_k(E).
\]
Multiply this by $-\frac12$ to get $\t_k^\CC(E;\z)_m$.
\end{proof}

\subsection{Complex Framing Principle}\label{CT:subsec:complex FP}

In order to be able to compute the higher complex torsion we need
to rephrase Theorem~\ref{CT:thm:complex torsion for closed AC
fibers} in terms of a smooth function $f:E\to \RR$ as in the
Framing Principle.

\begin{cor}[Complex Framing Principle]\label{CT:cor:complex
framing principle} Let $E\to B$ be a smooth bundle with closed
almost complex fibers, let $f:E\to\RR$ be any generic smooth
function and let $\Sig(f)$ be the singular set of $f$. Then
\[
   \t_k^\CC(E,\z)_m=\frac12m^kL_{k+1}(\z)
    p_\ast^\Sig(ch_k(T^vE)|\Sig(f))
\]
where the push-down operator
\[
    p_\ast:H^{2k}(\Sig(f))\to H^{2k}(B),
\]
with real or integral coefficients, is given in
subsection~\ref{FP:subsec:push-down transfer}.
\end{cor}

\begin{proof}
This is equivalent to Theorem~\ref{CT:thm:complex torsion for
closed AC fibers} by the commuting diagram in
Proposition~\ref{FP:prop:push-down transfer diagram}.
\end{proof}

We give an example of how the Complex Framing Principle is used
and how it compares with the (real) Framing Principle.

Take a complex $m$-plane bundle $\xi$ over $B$. Let
\[
    \CC P^{m-1}\to P(\xi)\to B
\]
be the associated projective bundle. Since $U(m)$ is connected,
$\pi_1B$ acts trivially on the homology of the fiber.
Consequently, all of the higher real and complex torsion
invariants are defined.

We first compute the real torsion. By the splitting principle, we
may assume that $\xi$ is a sum of line bundles $\xi=\oplus\ll_i$.
Then we choose a Morse function $g$ on $\CC P^{m-1}$ with exactly
$m$ critical points at the coordinate axes of $\CC^m$ so that $g$
is $U(1)^{m-1}$ invariant. This gives a fiberwise Morse function
$f$ on $P(\xi)$ with $m$ critical points $x_1,\cdots,x_m$ of
indices $0,2,4,\cdots,2m-2$ with corresponding negative eigenspace
bundles
\[
    \g_i\cong \Hom_\CC(\ll_i,\ll_1\oplus\cdots\oplus\ll_{i-1})
    \cong\bigoplus_{j<i}\Hom(\ll_i,\ll_j).
\]
Since there can be no algebraic incidences between these critical
points, the associated twisted cochain is zero so $C_\ast(f)$ has
zero torsion. Thus the higher FR torsion of $P(\xi)$ is given by
the correction term in the Framing Principle (\ref{FP:thm:Framing
Principle}):
\begin{align*}
    \t_{2k}(P(\xi))&=(-1)^k\z(2k+1)\sum_{i>j}ch_{2k}(\Hom_\CC(\ll_i,\ll_j))\\
    &=(-1)^k\z(2k+1)\frac12ch_{2k}(\Hom_\CC(\xi,\xi))
\end{align*}

Now we compute the higher complex torsion using the Complex
Framing Principle \ref{CT:cor:complex framing principle} which
tells us to use the entire complex vertical tangent bundle along
$x_i$:
\[
    \Hom_\CC(\ll_{i+1},\ll_1\oplus\cdots\oplus\widehat{\ll_{i+1}}\oplus\cdots\oplus\ll_n).
\]
Since $\Hom_\CC(\ll_{i+1},\ll_{i+1})$ is trivial, this is stably
equivalent to $\Hom_\CC(\ll_{i+1},\xi)$. Taking the sum over all
$i$ we get
\[
    \t_k^\CC(P(\xi),\z)_m=\frac12m^k\R\left(\frac1{i^k}\nL_{k+1}(\z)\right)
    ch_k(\End_\CC(\xi)).
\]

We see in this example that the factor of $\frac12$ compensates
for the fact that we are using the entire vertical tangent bundle
instead of just the negative eigenspace bundle.

\subsection{Nonempty boundary case}\label{CT:subsec:nonempty boundary case}

We now assume that $E\to B$ is a smooth bundle with compact almost
complex fibers $M^{2n}$ having real dimension $2n$ and nonempty
boundary. We assume $B$ is closed so that $\d E$ is a bundle over
$B$ with fiber $\d M$.

When the fibers $M$ have nonempty boundary, a correction term
should be added to the definition of higher complex torsion.

\begin{defn}[complex torsion with boundary]\label{CT:def:complex torsion in boundary case}
We define the \emph{complex torsion} of $E\to B$ to be the family
of cohomology classes $\t_k^\CC(E;\z)_m\in H^{2k}(B;\RR)$ for all
$k,m\geq1$ and nontrivial $m$-th roots of unity $\z$ given by
\begin{equation}\label{CT:eq:def of CT in boundary case}
    \t_k^\CC(E;\z)_m:=
    \frac1{2m}\t_k(\d E)
    -\frac12\t_k(S^{2n-1}(E)/Z_m;\r_\z).
\end{equation}
This differs from the closed fiber case by the \emph{boundary
correction term}:
\begin{equation}\label{CT:eq:boundary correction term}
  \frac1{2m}\t_k(\d E).
\end{equation}
This correction term is zero for $k$ odd and, for $k$ even, is
defined only under the assumption that $H_\ast(\d M;\QQ)$ is
$\pi_1B$-upper triangular.
\end{defn}

As in the closed fiber case these invariants are related by the
following transfer equation whenever the terms are defined.
\[
        \t_k^\CC(E,z)_m=\sum_{\z^r=z}\t_k^\CC(E,\z)_{mr}.
\]

The boundary correction term comes from Theorem~\ref{FP:thm:normal
disk bundle with boundary} and its proof. Without this correction
term, the following theorem would not be true.

\begin{thm}[complex vs real torsion in boundary
case]\label{CT:thm:complex vs real torsion in boundary case}
Suppose that the action of $\pi_1B$ on $H_\ast(M;\QQ)$ and
$H_\ast(\d M;\QQ)$ is upper triangular. Then
\[
    \t_{2k}^\CC(E,1)_1=\t_{2k}(E)
\]
for all $k\geq 1$.
\end{thm}

\begin{proof}
In Theorem~\ref{FP:thm:normal disk bundle with boundary} and its
proof we show that
\[
    \t_{2k}(E)=\frac12\t_{2k}(\d
    E)+\frac12\t_{2k}(S(E))-\frac12\t_{2k}(S(E)|_{\d E})
\]
and we also note that the last term is zero by
Theorem~\ref{FP:thm:even dim closed fibers} since $S(E)|_{\d E}$
has closed even dimensional stably parallelizable fibers. This
agrees with the definition of $\t_{2k}^\CC(E,1)_1$ by
Corollary~\ref{CT:lem:tangent and normal sphere torsion}.
\end{proof}

As we remarked in \ref{CT:rem:true def of CT using normal sphere
bundle}, a more natural way to write the definition of complex
torsion is the following
\begin{equation}\label{CT:eq:correct def of CT}
    \t_k^\CC(E;\z)_m=
    \frac1{2m}\t_k(\d E)
    +\frac12\t_k(S(E)/Z_m;\r_\z)
\end{equation}
where $S(E)$ is the sphere bundle of the complex vertical normal
bundle of $E$.

As in the case of the generalized MMM classes, the higher complex
torsion is stable when defined.

\begin{thm}[stability of complex torsion in boundary
case]\label{CT:thm:stability in boundary case} If $H_\ast(M;\QQ)$
and $H_\ast(\d M;\QQ)$ are $\pi_1B$ upper triangular then
\[
    \t_k^\CC(E,\z)_m=\t_k^\CC(E\times D^{2N},\z)_m.
\]
\end{thm}

\begin{proof}
By stability of $T_k$, we only need to show that the boundary
correction terms for $\t_k^\CC(E,\z)_m$ and $\t_k^\CC(E\times
D^{2N} ,\z)_m$ agree.

When $M$ has no boundary the boundary correction term of
$\t_k^\CC(E\times D^{2N},\z)_m$ is
\[
    \frac1{2m}\t_k^\CC(E\times
    S^{2N-1},\z)_m
    =\frac1{2m}\chi(S^{2N-1})\t_k(E)=0.
\]
When $M$ has a boundary this correction term is given by the
gluing formula (\ref{sec2:cor:gluing formula}) as follows.
\begin{align*}
    \frac1{2m}\t_k(\d(E\times D^{2N}))
    &=\frac1{2m}\left(\t_k(\d E\times D^{2N})
    +\t_k( E\times S^{2N-1})
    -\t_k(\d E\times S^{2N-1})\right)\\
    &=\frac1{2m}\t_k(\d E).
\end{align*}
The other two terms are zero by the product formula (since
$\chi(S^{2N-1})=0$).
\end{proof}

Finally, we note that the Complex Framing Principle still holds
with the same proof but it does not compute the complex torsion in
the boundary case because of the boundary correction term.

\begin{thm}[Complex Framing Principle in boundary
case]\label{CT:thm:CFP in boundary case} Let $f:E\to\RR$ be any
generic smooth function whose fiberwise gradient points outward
along $\d E$. Let $\Sig(f)$ be the singular set of $f$. Then
\begin{align*}
    \t_k^\CC(E,\z)_m&=\frac12 m^kL_{k+1}(\z)\frac1{k!}T_k(E)+\frac1{2m}\t_k(\d
    E)\\
    &=\frac12m^kL_{k+1}(\z)p_\ast^\Sig(ch_k(T^vE)|\Sig(f))+\frac1{2m}\t_k(\d
    E).
\end{align*}
\end{thm}

%% file: sec1.tex
\vfill\eject\section{Definition of higher FR-torsion}

\begin{enumerate}
  \item Generalized Morse functions
  \item Families of chain complexes
  \item Monomial functors
  \item Filtered chain complex
  \item Subfunctors
  \item The Whitehead category
  \item Definition of higher FR torsion
  \item {Families of matrices as flat superconnections}
  \item Independence of birth-death points
  \item Positive suspension lemma
\end{enumerate}

In this section we review the definition of higher
Franz-Reidemeister torsion and, more generally, the cellular chain
complex functor associated to any fiberwise oriented generalized
Morse function.

Suppose we have a smooth bundle
\begin{equation*}
  M\to E\to B
\end{equation*}
where $B$ is a closed manifold and $M$ is a compact $n$-manifold
with $\d M=\d_0M\coprod\d_1M$. Let $\d E=\d_0E\coprod\d_1E$ where
$\d_iE$ is a bundle over $B$ with fiber $\d_iM$.

Given a fiberwise generalized Morse function
  \[f:(E,d_0E,\d_1E)\to (I,0,1)\]
  and a representation $\r:\pi_1E\to G$ where $G$ is a subgroup of
  the group of units of a ring $R$ we construct a family of
cellular chain complexes $C_\ast(f;R)$, with twisted
$R$-coefficients, giving an $(R,G)$-expansion functor
  \[
    C_\ast(f;R):B\to \Wh\bu(R,G)
  \]
  where $\Wh\bu(R,G)$ is the Whitehead category.
  The expansion functor consists of a poset functor and a
  \emph{twisted cochain} in the sense of
  \cite{[Brown59:TwistedTensor]}.

In the special case when $(R,G)=(M_r(\CC),U(r))$, we have
cohomology classes
  \[
    \t_k\in H^{2k}(\Wh\bu^h(M_r(\CC),U(r));\RR)
  \]
  which we call the \emph{universal Franz-Reidemeister torsion invariants}.
  Here, $\Wh\bu^h(-,-)$ is the simplicial full subcategory of
  $\Wh\bu(-,-)$ given by acyclic chain complexes. Thus, in the
  basic construction, we obtain higher FR torsion invariants for
  $(E,\d_0E)$ in the case when we have a Hermitian coefficient
  system $\F$ on $E$ with respect to which $(M,\d_0M)$ is acyclic.
  They are given by
  \[
    \t_k(E,\d_0E;\F)=C_\ast(f;\F)^\ast(\t_k)\in H^{2k}(B;\RR).
  \]

\subsection{Generalized Morse functions}\label{sec1:subsec:GMFs}

A \emph{generalized Morse function} (GMF) on a smooth $n$-manifold
$M$ is a smooth function
\[
    f:M\to\RR
\]
having Morse (nondegenerate) and \emph{birth-death} (b-d)
singularities. We recall that, at a Morse singularity, the
function $f$ can be written in the form
\begin{equation}\label{sec1:eq:Morse singuarlity}
  f(x,y)=-\length{x}^2+\length{y}^2+C
\end{equation}
with respect to local coordinates $(x,y)\in \RR^i\times \RR^{n-i}$
for $M$ and $C$ is a constant (the critical value). At a
birth-death point, $f$ can be written as
\begin{equation}\label{sec1:eq:birth-death of index i-1}
  f(x,y)=-\length{(x_1,\cdots,x_{i-1})}^2+x_i^3+\length{y}^2+C.
\end{equation}

We choose a Riemannian metric on $M$ which we will assume to be
\emph{standard} in a coordinate neighborhood of each critical
point. By this we simply mean that the $x$-plane is perpendicular
to the $y$-plane at each point in a small neighborhood. In the
case of a birth-death point we also require the $x_i$ direction be
perpendicular to the tangent plane spanned by $\frac{\d}{\d x_j}$
for $j=1,\cdots,i-1$ in a small neighborhood of the b-d point.

To begin, we choose the metric on $M$ so that the coordinate chart
given by $(x,y)$ is an isometry. We say in that case that the
coordinates are \emph{isometric}. Then the gradient of $f$ will be
the transpose of the derivative, i.e,
\[
    \grad f(x,y)=(-2x,2y)^t=-\sum 2x_j\frac{\d}{\d x_j}+\sum 2
    y_k\frac{\d}{\d y_k}
\]
near a Morse point and
\[
    \grad
    f(x,y)=(-2x_1,\cdots,-2x_{i-1},3x_i^3,2y_1,\cdots,2y_{n-i})^t
\]
near a b-d point.

If the coordinates are not isometric but the metric is still
standard then the subsets of the coordinate neighborhood given by
$x=0$ and $y=0$ will still be totally geodesic in both cases as
will the subset given by $x_i=0$ and the $x_i$-axis in the b-d
case. Consequently, the union of the set of trajectories of $\grad
f$ converging to the critical point is unchanged.

The b-d singularity given in (\ref{sec1:eq:birth-death of index
i-1}) has \emph{index} $i-1$ since this is the index of the second
derivative $D^2 f$ of $f$. The vectors $\frac{\d}{\d x_j}$ span
the \emph{negative eigenspace} of $D^2 f$. The negative eigenspace
is a subspace of the tangent space of $M$ at any critical point
which depends on the choice of metric. However, it is unique up to
convex choice. An \emph{orientation} for $f$ is defined to be an
orientation for the negative eigenspace of $D^2f$ at every
critical point of $f$.

In a generic $p$-parameter family of GMF's birth-death points will
occur on a codimension one subset of the parameter space. Each b-d
point will be \emph{generically unfolded} in the sense that the
family of functions $f_t$ has the form:
\begin{equation}\label{sec1:eq:unfolding of bd point}
  f_t(x,y)=-\length{(x_1,\cdots,x_{i-1})}^2
  +x_i^3+t_0x_i+\length{y}^2+C
\end{equation}
with respect to parameter coordinates $t_0,t_1,\cdots,t_{p-1}$ and
$t$-dependent local coordinates $(x,y)$ for $M$. The singular set
$\Sig(f)$ is a smooth manifold given in this coordinate patch by
\[
    t_0=-3x_i^2
\]
(and all other coordinates equal to zero). Choose a $t$-dependent
metric for $M$ so that the coordinates $x,y$ are isometric for
each $t$. Then, for each $t$, the metric will be standard at the
Morse points of $f_t$. [If we reparametrize the $x_i$ coordinate
so that the nondegenerate quadratic singularities that occur when
$t_0<0$ have the form $\pm x_i^2$ as required by
(\ref{sec1:eq:Morse singuarlity}) then the coordinate $x_i$ will
no longer be isometric, but it will still be perpendicular to the
other coordinates.]

Given an oriented Morse function
\[
    f:M\to[0,1]
\]
with the property that the boundary of $M$ is the disjoint union
of $\d_0 M = f^{-1}(0)$ and $\d_1M=f^{-1}(1)$ we get a
\emph{cellular chain complex} $C_\ast(f)$ which is a free
$\ZZ[G]$-complex ($G=\pi_1M$) with a basis which is well-defined
up to permutation and multiplication by elements of $G$. The
homology of the cellular chain complex is the relative homology
\[
    H_k(C_\ast(f))=H_k(M,\d_0M;\ZZ[G]).
\]
In a generic family of GMF's, the orientations of negative
eigenspaces should be chosen in such a way that the incidence
between the cancelling critical points at a birth-death point is
positive ($+g$ for some $g\in G$).

\begin{defn}[fiberwise oriented GMF]\label{sec1:def:fiberwise oriented GMF}
A \emph{fiberwise oriented GMF} on a smooth bundle $E\to B$ is
defined to be a smooth function
\[
    f:E\to\RR
\]
which is a GMF on each fiber $M$ together with a continuous family
of orientations for the negative eigenspace of $D^2f_t$ so that
the incidence coefficients of cancelling critical points is
positive.
\end{defn}

There are two problems with fiberwise oriented GMF's. The first is
that they may not exist and the second is that they may not be
unique. For the existence we know that a fiberwise oriented GMF
exists for a smooth bundle $E\to B$ whose fiber dimension is
greater than or equal to its base dimension:
\[
    \dim M\geq\dim B.
\]
There are also special arguments in the case $\dim M\leq2$. To
make the fiber dimension arbitrarily large, we can take the
product with a large even dimensional sphere $S^{2N}$. This
multiplies all higher torsion invariants by $2$ (the Euler
characteristic of $S^{2N}$) so we can divide by $2$ to get the
invariant for the original bundle. (This procedure is justified by
Lemma~\ref{sec2:lem:product formula} below.) We could also take
the product with any large dimensional disk and round off the
corners.

For the uniqueness problem we take a special type of fiberwise
oriented GMF called a fiberwise ``framed function.''

\begin{defn}[framed function]\label{sec1:def:framed function}
A \emph{framed function} on a smooth manifold $M$ is a GMF $f:M\to
\RR$ together with a framing of the negative eigenspace of $D^2f$
at each critical point. At a birth-death point of index $i-1$ we
also choose the $i$-th framing vector to be \[v_i=\frac{\d}{\d
x_i}\]in standard coordinates (\ref{sec1:eq:birth-death of index
i-1}). A \emph{fiberwise framed function} on a smooth bundle $E\to
B$ is a smooth function $f:E\to\RR$ together with continuous
families of tangent vectors $v_1,v_2,\cdots$ giving the
restriction $f_t:M_t\to\RR$ of $f$ to each fiber the structure of
a framed function.
\end{defn}

\begin{rem}[orientation of a framed function]
\label{sec1:rem:orientation of a framed function}
A fiberwise framed function gives a fiberwise oriented GMF if we
take the orientation given by reversing the framing vectors:
$(v_i,\cdots,v_1)$.
\end{rem}

\begin{thm}[Framed Function Theorem\cite{[I:FF]}]\label{sec1:thm:framed function theorem} The space of framed functions
on a compact smooth $n$-manifold $M$ is $n-1$ connected.
Consequently, a smooth bundle $E\to B$ admits a fiberwise framed
function $f:E\to\RR$ if $\dim M\geq\dim B$ and $f$ is unique up to
framed homotopy if $\dim M>\dim B$ (i.e., given two fiberwise
framed functions $f,g$ there is a fiberwise framed function on the
bundle $E\times I\to B\times I$ which restricts to $f,g$ on
$E\times 0$ and $E\times1$).
\end{thm}

\begin{rem}[$C^1$-local framed function theorem]\label{sec1:rem:C1 local FFT}
In \cite{[I:C1-local]} it was show (under the same dimension
restrictions) that the fiberwise framed function $f$ can be chosen
to be arbitrarily $C^1$-close to any given smooth function $g$.
Furthermore, $f$ can be chosen to be equal to $g$ outside of an
arbitrarily small neighborhood of the set on which $g$ is not
fiberwise framed. In particular $\Sig(f)$ will lie in an
arbitrarily small neighborhood of $\Sig(g)$.
\end{rem}

\subsection{Families of chain
complexes}\label{sec1:subsec:families of chain complexes}

Given a fiberwise oriented GMF $f$ on the total space of a smooth
bundle $E\to B$ and a \emph{vertical metric} on $E$ (i.e., a
metric on the \emph{vertical tangent space} $T^vE$) which we
always assume is standard near each critical point, we get a
family of chain complexes parametrized by $B$. We call it
$C_\ast(f)$.

If $f$ is chosen to be fiberwise framed then $C_\ast(f)$ will be
unique up to homotopy by the Framed Function Theorem (and the
uniqueness up to homotopy of the vertical metric) and we use it to
define the higher FR torsion invariants. If $f$ is not framed we
can still recover the higher FR torsion from the Framing
Principle. Thus, it is useful to consider $C_\ast(f)$ even when it
is not canonical.

In this subsection we will describe what we mean by a ``family of
chain complexes.'' For the technical details, we refer to
\cite{[I:BookOne]}. To (hopefully) simplify the presentation we
reverse the direction of the arrows. Since every category is
homotopy equivalent to its opposite category, this does not change
anything.

The ideas is simple. We take a small $k$ simplex $\s$ in the base
$B$. Assume first that the function $f$ is Morse over every point
in $\s$. Then the singular set forms a disjoint union of sheets
\[
    \Sig(f|\s)\cong \Sig(f_b)\times\Delta^k
\]
where $b\in\s$ is, say, the barycenter. If $f$ has b-d points over
$\s$ we delete all components of $\Sig(f|\s)$ which contain b-d
points. The remaining set of sheets forms a \emph{poset}
(partially ordered set) $P$ with a \emph{standard ordering} given
by $x<y$ if the critical value of $x$ is less than the critical
value of $y$ over every point in $\s$:
\[
    x<y\Leftrightarrow (\forall t\in\s) f_t(x)<f_t(y).
\]

In this paper we will not use the standard partial ordering on
$P$. We use instead the (adjusted) \emph{minimal partial ordering}
which is given as follows. Let $P^+(\s)$ be the set of all
components of $\Sig(f|\s)$. We assume that $\s$ is transverse to
the \emph{bifurcation set} of $f$ (the set of all $t\in B$ for
which $f_t$ is not Morse). So, this will be a finite set. We take
the transitive relation on $P^+(\s)$ generated by the relation
$x<y$ if for some $t\in\s$ a trajectory of $\grad f_t$ goes from a
point in $x$ up to a point in $y$.

\begin{prop}[minimal partial ordering]\label{sec1:prop:minimal
ordering} If $\s$ is sufficiently small this transitive relation
will also be anti-reflexive. A fortiori, it will satisfy the
following. If $x<y$ in $P^+(\s)$ then
\[
    f_t(x_i)<f_t(y_j)
\]
for every $t\in\s$, $x_i\in x\cap \Sig(f_t)$ and $y_j\in
y\cap\Sig(f_t)$.
\end{prop}

Let $P=P(\s)$ be the subset of $P^+(\s)$ of Morse components with
the induced partial ordering. Then the proposition above tells us
that the standard ordering on $P$ is a refinement of the minimal
ordering. The minimal ordering depends on the vertical metric
which we always assume is standard near the critical points.
Later, in subsection~\ref{sec1:subsec:independence of bd points},
we will adjust the partial orderings of both $P(\s)$ and
$P^+(\s)$.

Note that $P$ is \emph{graded} since every element $x\in P$ has a
\emph{degree} given by the index of $f$ at $x$. However, it is
useful not to assume any relationship between the partial ordering
and the grading.

Over each vertex $v_i$ of $\s$ we obtain, by standard Morse
theory, a cellular chain complex $C_\ast(f_{v_i})$. Recall that
the definition of this complex requires that a path be chosen from
the base point of $E$ to each critical point of $f_{v_i}$. Since
$\Delta^k$ is contractible, this can be done consistently over
every point in $\s$. We also need an orientation of the
topological cells associated to each critical point. This choice
is given by the orientation of $f$.

With these choices, $C_\ast(f_{v_i})$ is a free based chain
complex over some ring $R$. ($R=\ZZ[G]$ with $G=\pi_1E$ is the
universal case.) The basis set will be the graded poset $P$ for
each vertex $v_i$. However the boundary map will depend on $i$. We
call the boundary map $\f_0(i)$. Thus
\[
    C_\ast(f_{v_i})=(PR,\f_0(i))
\]
where $PR=R^P$ denotes the free right $R$-module generated by $P$.
Note that $\f_0(i)$ is an endomorphism of $PR$ satisfying the
following.
\begin{enumerate}
  \item $\f_0(i)$ is homogeneous of degree $-1$.
  \item $\f_0(i)(x)$ is an $R$-linear combination of elements $y<x$
  of $P$.
  \item $\f_0(i)^2=0$
\end{enumerate}
Any endomorphism of $PR$ satisfying condition (2) will be called
\emph{strictly upper triangular}. (If the condition is weakened to
$y\leq x$ we call it \emph{upper triangular}.)

Over each edge $(v_i,v_j)$ in $\s$ with $i<j$ we have an upper
triangular chain isomorphism
\[
    \Phi_1(i,j):C_\ast(f_{v_j})=(PR,\f_0(j))\xrarrow{\approx}(PR,\f_0(i))=C_\ast(f_{v_i}).
\]
This chain isomorphism is ``close'' to the identity map on $PR$ so
we write it as
\[
    \Phi_1(i,j)=1+\f_1(i,j).
\]
Then $\f_1(i,j)$ will be a strictly upper triangular endomorphism
of $PR$ of degree $0$ satisfying the equation
\[
  \f_0(i)(1+\f_1(i,j))=(1+\f_1(i,j))\f_0(j),
\]
which can also be written as
\begin{equation}\label{sec1:eq:e2(i,j) equation}
    \f_0(i)\f_1(i,j)-\f_1(i,j)\f_0(j)=\f_0(j)-\f_0(i).
\end{equation}

Given three vertices $v_i,v_j,v_k$ with $i<j<k$, the two-parameter
family of functions $f$ over the corresponding triangular face of
$\s$ gives us a strictly upper triangular chain homotopy
\[
    \f_2(i,j,k):\Phi_1(i,j)\Phi_1(j,k)\simeq \Phi_1(i,k).
\]
This can be written in the form
\begin{multline}\label{sec1:eq:e2(i,j,k) equation}
  \f_0(i)\f_2(i,j,k)-\f_1(i,j)\f_1(j,k)+\f_2(i,j,k)\f_0(k)\\=\f_1(j,k)-\f_1(i,k)+\f_1(i,j).
\end{multline}

\begin{defn}[$\Delta^k$-family of chain
complexes]\label{sec1:def:Delta k family of chain complexes} Let
$R$ be any ring. Then by a \emph{$\Delta^k$-family of chain
complexes} over $R$ we mean a pair $(P,\f)$ where $P$ is a finite
nonnegatively graded poset and $\f$ is a system of endomorphism
\[
    \f_p(a_0,\cdots,a_p):PR\to PR
\]
for all $0\leq a_0<a_1<\cdots<a_p\leq k$ which are homogeneous of
degree $p-1$, which are {strictly upper triangular} ($\f_p(a)(x)$
is a linear comb. of $y<x$ in $P$) and which satisfy the following
(where $\f_{-1}()=0$).
\begin{equation}\label{sec1:eq:defining equation for families of chain complexes}
  \sum_{i=0}^p(-1)^i\f_i(a_0,\cdots,a_i)\f_{p-i}(a_i,\cdots,a_p)
  =
  \sum_{i=0}^p(-1)^i\f_{p-1}(a_0,\cdots,\widehat{a_i},\cdots,a_p).
\end{equation}
Let $\Delta^k(R)$ denote the set of all $\Delta^k$-families of
$R$-complexes where $P$ lies in some fixed universe $\Omega$.
These form a simplicial set $\Delta^\bullet(R)$ in an obvious way.
\end{defn}

Condition (\ref{sec1:eq:defining equation for families of chain
complexes}) can be expressed in the following succinct form where
$\f=\sum \f_p$ is a sum of cocycles on $\Delta^k$ with
coefficients in the graded ring $\End_R(PR)$ and
$\f'=\sum(-1)^p\f_p$.
\begin{equation}\label{sec1:eq:de = e u e}
    \f'\cup \f=\delta \f.
\end{equation}

\begin{rem}[twisted cochains]\label{sec1:rem:Brown discovered twisted cochains}
Ed Brown \cite{[Brown59:TwistedTensor]} was the first to consider
sums of cochains $\f$ satisfying (\ref{sec1:eq:defining equation
for families of chain complexes}), (\ref{sec1:eq:de = e u e}). He
called them \emph{twisted cochains}. To show that these are the
same, we can rewrite the conditions as
\begin{equation}\label{sec1:eq:def of twisted cochain}
  \d\f_p=\f_{p-1}\d-\sum_{i=1}^{p-1}(-1)^i\f_i\cup\f_{p-i}
\end{equation}
where $\d\f_p=\f_0\f_p+(-1)^p\f_p\f_0$. Since the differential
$\f_0$ varies from point to point, $\f$ is, strictly speaking, a
twisted cochain with coefficients in a differential graded (DG)
category whose DG Hom sets have the form
\begin{equation}\label{sec1:eq:DG Hom sets}
  \Hom(A,B)=\Hom((PR,\f_0(A)),(PR,\f_0(B))).
\end{equation}
To express the relationship between $\f,R,P$ we say that $\f$ is a
twisted cochain \emph{over} $R$ which is \emph{subordinate} to
$P$.
\end{rem}

\begin{rem}[$A_\infty$-functors]\label{sec1:rem:A infty functors}
Let $\Phi_p$ be given by $\Phi_1=1+\f_1$ and $\Phi_p=\f_p$ for
$p\noteq1$. If $\f$ is a twisted cochain on the nerve of a
category $\C$ then $\Phi$ is an \emph{$A_\infty$-functor} from
$\C\op$ to the DG category described above. To see this we convert
notation as follows.
\begin{align*}
    \delta\f-\f'\cup\f&=\sum_{i=0}^p(-1)^i\f_{p-1}\d_i-\f'\cup\f\\
        &=\sum_{i=1}^{p-1}(-1)^i\Phi_{p-1}\d_i-\Phi'\cup\Phi
        \end{align*}
Then note that:
        \begin{align*}
    \Phi_{p-1}\d_i &=\Phi_{p-1}(1^{i-1},m_2,1^{p-i-1})\\
    \Phi'\cup\Phi&=\sum_{i=0}^p(-1)^im_2(\Phi_i,\Phi_{p-i})\\
        &=\sum_{i=1}^{p-1}(-1)^im_2(\Phi_i,\Phi_{p-i})+m_1(\Phi_p),
\end{align*}
where $m_2(f,g)=g\circ f$ ($f\circ g$ in $\C\op$) and $m_1=\d$.
Thus (\ref{sec1:eq:de = e u e}) becomes:
\[
    \sum_{i=1}^{p-1}(-1)^i\Phi_{p-1}(1^{i-1},m_2,1^{p-i-1})=m_1(\Phi_p)
    +\sum_{i=1}^{p-1}(-1)^im_2(\Phi_i,\Phi_{p-i}).
\]
\end{rem}

For simplicity we follow the twisted cochain formalism. We will
vary the poset to obtain what we call ``monomial functor'' and
later we will generalize the construction to an ``expansion
functor.'' Just as in the case of $A_\infty$ functors, these are
not true functors. However, we will see that they can be
reinterpreted as true functors in
Remark~\ref{sec1:rem:reinterpreting monomial functors} below.

\subsection{Monomial functors}\label{sec1:subsec:monomial
functors}

Suppose that $f:E\to \RR$ is a fiberwise oriented Morse function
and $R$ is an untwisted coefficient ring. Then, we can apply the
construction above to every small smooth simplex $\s:\Delta^p\to
B$ to obtain $p$-cochains $\f_p$ on $B$ with coefficients in the
strictly upper triangular degree $p-1$ part of the graded ring
$\End_R(PR)$ whose sum $\f=\sum \f_p$ satisfies (\ref{sec1:eq:de =
e u e}) above where $P$ is a functor from the category $\simp B$
of small smooth simplices in $B$ to the category $\nP$ of
partially ordered graded sets. (Morphisms in $\nP$ are degree $0$
order preserving bijections $\a:P\to Q$, i.e., $x<
y\Rightarrow\a(x)<\a(y)$ and $\deg(x)=\deg\a(x)$.)

We recall that the \emph{category of simplices} $\simp X$ of any
simplicial set $X$ has as objects all pairs $(\s,[k])$ where
$k\geq0$ and $\s\in X_k$ and morphisms $(\t,[j])\to (\s,[k])$ are
defined to be (nondecreasing) maps $a:[j]\to[k]$ so that
$\t=a^\ast(\s)$.

This construction needs to be generalized to twisted coefficients
$R$ (with $G\subseteq R^\times$ containing the image of $\pi_1E$)
as follows.

\begin{defn}[the category $\nP(G)$]\label{sec1:def:category P(G)}
Let $\nP(G)$ be the category whose objects are finite partially
ordered graded sets. A morphism in $\nP(G)$ is a right
$G$-equivariant mapping
\[
    \widetilde{\a}:P\times G\to Q\times G
\]
so that the induced map on orbits ${\a}:P\to Q$ is a degree $0$
order preserving bijection. We say that $\a:P\to Q$ is a
\emph{$G$-morphism} to indicate that it is covered by a
$G$-equivariant mapping $\widetilde{\a}$.
\end{defn}

If $G\subseteq R^\times$, the free right $R$-module $PR$ is a
functor on $\nP(G)$. A $G$-morphism $\a:P\to Q$ induces an
isomorphism ${\a}_\ast:PR\to QR$ by
\[
    {\a}_\ast\left(\sum x_ir_i\right)=\sum \pi\widetilde{\a}(x_i,1)r_i
\]
where $\pi(x,g)r=x(gr)$. The inverse of $\a_\ast$ will be denoted
$r$ (for \emph{restriction}). Conjugation by $\a_\ast$ gives a
ring isomorphism
\[
    (\,)_\#:\End_R(PR)\to\End_R(QR)
\]
making $\End_R(PR)$ into a functor on $\nP(G)$. If $f$ is strictly
upper triangular then so is $f_\#$.

If $f:E\to\RR$ is a fiberwise oriented Morse function and $R$ is a
coefficient ring twisted by a representation $\pi_1E\to G\subseteq
R^\times$ then the singular set of $f$, with the minimal partial
ordering, gives a functor
\[
    P=\Sig(f):\simp B\to \nP(G)
\]
and for each $p\geq0$ we get a $p$-cochains $\f_p$ on $B$ with
coefficients in the strictly upper triangular degree $p-1$ part of
$\End_R(PR)$ whose sum $\f=\sum \f_p$ satisfies (\ref{sec1:eq:de =
e u e}). The functor $P$ together with the twisted cochain $\f$
forms what we call an \emph{$(R,G)$-monomial functor} $(P,\f)$ on
$B$.

\begin{defn}[$(R,G)$-monomial functor]\label{sec1:def:RG monomial functor}
An \emph{$(R,G)$-monomial functor} on a simplicial set $X$ is
defined to be a pair $(P,\f)$ where $P$ is a functor
\[
    P:\simp X\to \nP(G)
\]
and $\f=\sum \f_p$ where $\f_p$ is a $p$-cochain on $X$ with
coefficients in the strictly upper triangular degree $p-1$ part of
$\End_R(PR)$ so that $\delta \f=\f'\cup \f$. We say that $\f$ is a
\emph{twisted cochain} \emph{subordinate} to the functor $P$.
\end{defn}

To clarify this definition, when we say that $\f_p$ is a
$p$-cochain on $X$ we mean that $\f_p$ is a function that takes
every element $x\in X_p$ to a strictly upper triangular
endomorphism $\f_p(x)\in\End_R(P(x)R)$. When we write $\delta
\f=\f'\cup \f$ we mean:
\[
    \sum_{i=0}^n(-1)^i\f_{n-1}(\d_ix)_\#
    =\sum_{p+q=n}(-1)^p\f_p(f_p(x))_\#\f_q(b_q(x))_\#
\]
in $\End_R(P(x)R)$. In other words, we are inducing the
endomorphisms on the faces of $x$ back up to $x$ and the condition
is that these form an element
\begin{equation}\label{sec1:eq:Delta k family xi(x)}
  (P(x),\f_\#)\in\Delta^k(R).
\end{equation}

\begin{defn}[twisted tensor product]\label{sec1:def:twisted tensor
product=total complex} Ed Brown's \emph{twisted tensor product}
(giving the homology of the total space
\cite{[Brown59:TwistedTensor]}) is called the \emph{total complex}
of $(P,\f)$ in \cite{[I:BookOne]}:
\begin{equation}\label{sec1:eq:total complex E(P,e)}
  C_\ast(X)\otimes_\f (PR,\f_0)=E_\ast(P,\f)
  =\bigoplus_{n\geq0}\bigoplus_{x\in X_n}x\otimes
  P(x)R
\end{equation}
with boundary $\d_\f$ given by:
\[
    \d_\f (x\otimes y)=\d x\otimes r(y)-\sum_{p+q=n}
    (-1)^pf_p(x)\otimes r(\f_q(b_q(x))_\#(y))
\]
where each $r$ represents the appropriate restriction map. For
example,
\[
    \d x\otimes r(y)=\sum (-1)^i\d_ix\otimes\d_i^\ast(y).
\]
\end{defn}

\subsection{Filtered chain complexes}\label{sec1:subsec:filtered chain complex}

This subsection describes the transition from Morse theory to the
algebra of twisted cocycles. The basic idea is to take the total
singular complex of $E|\s$ for each small $k$-simplex $\s$
together with the filtration given by a fiberwise oriented GMF.
The theory of filtered chain complexes is joint work of Klein and
the author. This includes the existence of minimal filtered
subcomplexes. Minimal filtered chain complexes are shown to be
isomorphic to twisted cochains in \cite{[I:BookOne]}. However, I
was not aware of the connection to Brown's paper
\cite{[Brown59:TwistedTensor]} until after I talked to Bernhard
Keller.

Suppose that, as above, we have a fiberwise oriented GMF
$f:E\to\RR$ and a small smooth $k$-simplex $\s$ in $B$ which is
transverse to the {bifurcation set} of $f$. Then we have an open
subset of $E|\s$ associated to closed subsets of the poset
$P^+(\s)$ of all components of $\Sig(f|\s)$.

A subset $Q$ of a poset $P$ is called \emph{closed} if it contains
every element of $P$ which is less than an element of $Q$. The
complement of $Q$ in $P$ is denoted $P/Q$. (It is the
\emph{quotient}.) For $P=P^+(\s)$, this means that there are no
trajectories of $\grad f_t$ going up from a point in $P/Q$ up to a
point in $Q$ for any $t\in\s$.

For every closed subset $Q^+$ of $P^+(\s)$ let $E^{Q^+}|\s$ be the
union over all $t\in\s$ of the set all points $x$ in $M_t$ so that
the trajectory of $\grad f_t$ passing through $x$ contains an
element of the set
\begin{equation}\label{sec1:eq:d0M cup Q+}
  \d_0M_t\smalldisj Q^+
\end{equation}
in its closure. In other words, either $x$ is an element of
(\ref{sec1:eq:d0M cup Q+}) or the trajectory of $\grad f_t$
through $x$ comes from an element of (\ref{sec1:eq:d0M cup Q+}).
Another description of the same set $E^{Q^+}|\s$ is that it is the
complement in $E|\s$ of the set of points multiply incident over
$P^+(\s)/Q^+$.

Note that in the special case of the empty set $Q^+=\emptyset$,
the space $E^\emptyset|\s$ is homotopy equivalent to $\d_0E|\s$.

\begin{defn}[filtered chain complex]\label{sec1:def:filtered chain complex of f}
The \emph{filtered chain complex} of $f$ over $\s$ is defined to
be the total singular complex $S_\ast(E|\s;\r)$ of $E|\s$ with
coefficients $R$ twisted by a representation
\[
    \r:\pi_1E\to G\subseteq R^\times
\]
together with the array of subcomplexes given by the total
singular complex of $E^{Q^+}|\t$ for all closed subsets $Q^+$ of
$P^+(\s)$ and all faces $\t$ of $\s$ and the \emph{dual classes}
for \emph{Morse layers} defined below.
\end{defn}

Suppose that $Q^+=A\coprod \{x\}$ where $A$ is a closed subset of
$Q^+$ and $x$ is a maximal element. Then the \emph{layer} of the
filtered chain complex corresponding to $(A,x,\t)$ is the
subquotient complex
\begin{equation}\label{sec1:eq:layer of (A,x,tau)}
    \frac{S_\ast(E^{Q^+}|\t;\r)}
    {S_\ast(E^A|\t;\r)+\sum S_\ast(E^{Q^+}|\d_i\t;\r)}.
\end{equation}
Depending on whether $x$ is a Morse component of $\Sig(f|\s)$ or a
component containing b-d points, we call (\ref{sec1:eq:layer of
(A,x,tau)}) a \emph{Morse} or \emph{birth-death layer}.

The birth-death layers are all acyclic. Morse layers are, however,
nontrivial in exactly one degree equal to the index $i$ of $x$
plus the dimension $d$ of $\t$. In degree $i+d$ the Morse layer
(\ref{sec1:eq:layer of (A,x,tau)}) has homology isomorphic to $R$.
The generator is dual to the $d+i$ dimensional cohomology class
given by the ascending disks of the Morse critical point of $f_t$
corresponding to $x$ where $t$ ranges through all interior points
of $\t$. Note that the orientation of the function $f$ gives a
$\ZZ$-orientation of the normal bundle of these ascending disks.
If we choose paths from the base point of $E$ to each component of
$\Sig(f|\s)$ we get an $R$-orientation, i.e., a generator of the
$d+i$ cohomology of the layer (\ref{sec1:eq:layer of (A,x,tau)}).
We call this cohomology generator the \emph{dual class}. It
depends of the choice of paths from $x$ to the base point of $E$.

The following theorem is the theorem from \cite{[IK2:FR-torsion1]}
on the existence of a ``minimal subcomplex'' of a filtered chain
complex reinterpreted in \cite{[I:BookOne]} and expressed in the
language of \cite{[Brown59:TwistedTensor]}.

\begin{thm}[existence of twisted
cochains]\label{sec1:thm:existence of twisted cochains} There
exists a $\Delta^k$-family of chain complexes $(P(\s),\f(\s))$
which is uniquely determined up to simplicial homotopy by the fact
that its total complex (the twisted tensor product),
\[
    C_\ast(\Delta^k)\otimes_{\f(\s)} P(\s)R
\]
(Definition~\ref{sec1:def:twisted tensor product=total complex}),
together with the filtration given by closed subsets of $P(\s)$
and faces of $\s$, is filtered quasi-isomorphic to the total
singular complex of the pair $(E|\s,\d_0E|\s)$ with twisted
$R$-coefficients and so that the basis element of each layer in
the twisted tensor product is dual to the dual class of the
filtered chain complex.
\end{thm}

Here, \emph{filtered quasi-isomorphic} means there is a filtration
preserving chain map which induces a quasi-isomorphism on every
subcomplex (and thus on every layer). The twisted tensor product
is \emph{minimal} in the sense that every Morse layer is free of
rank $1$, every birth-death layer is zero as is the base layer
corresponding to the empty set.

\begin{rem}[definition of FR torsion is an algebraic
problem]\label{sec1:rem:def of FR torsion is algebraic} Note that
Brown's theorem (that the twisted tensor product gives the
homology of the total space) holds by construction over each
simplex and thus, by Mayer-Vietoris, for the entire bundle. The
higher Franz-Reidemeister torsion is the obstruction to deforming
the twisted tensor product given by a fiberwise framed function on
$E$ to one given by a canonical algebraically defined twisted
cochain in situations when the latter exists. Therefore, it is a
purely algebraic problem to determine under what conditions higher
FR-torsion is defined.
\end{rem}

\subsection{Subfunctors}\label{sec1:subsec:subfunctors}

A subfunctor of a monomial functor $(P,\f)$ is given by a system
of closed subsets of the poset functor $P$.

If $\f$ is a twisted cochain over $R$ subordinate to the constant
poset $P$ then any closed subset $Q$ of $P$ defines a twisted
cochain $\f_Q$ with coefficients in $\End(QR)$: For every vertex
$v$, $Q$ generates a subcomplex $(QR,\f_Q(v))$ of $(PR,\f(v))$.
For every edge $(v_0,v_1)$, we get a commuting diagram of chain
complexes:
\begin{displaymath}
\begin{CD}
(PR,\f(v_1))  @>1+\f(v_0,v_1)>\approx>  (PR,\f(v_0))\\
@A\subseteq AA               @A\subseteq AA \\
(QR,\f_Q(v_1))   @>1+\f_Q(v_0,v_1)>\approx>    (QR,\f_Q(v_2))
\end{CD}
\end{displaymath}
and so on. We also get a \emph{quotient twisted cochain} $\f/\f_Q$
subordinate to $P/Q$ giving the quotient complexes
$(PR,\f(v))/(QR,\f_Q(v))$ at vertices, etc.

Given a poset functor $P:\C\to\nP(G)$ a \emph{subfunctor} $Q$ is
defined to be a closed subset $Q(A)$ of $P(A)$ for each object $A$
of $\C$ which is \emph{invariant} in the sense that
\[
    {f_\ast}(Q(A))=Q(B)
\]
for all morphisms $f:A\to B$ in $\C$. We also have the
\emph{quotient functor} $P/Q$ given by $(P/Q)(A)=P(A)/Q(A)$.

Finally, if $\xi=(P,\f)$ is an $(R,G)$-monomial functor then a
poset subfunctor $Q$ of $P$ gives a \emph{monomial subfunctor}
\[
    \eta=(Q,\f_Q)
\]
of $\xi$ and we also have the \emph{monomial quotient functor}
\[
    \xi/\eta=(P/Q,\f/\f_Q).
\]
We say that $\xi$ is an \emph{extension} of $\eta$ by $\xi/\eta$.

We need one more definition. If a subset $S$ of a poset $P$ is
both open and closed we call it \emph{independent}. In that case
$P$ is a disjoint union of two closed subsets $Q,S$, i.e., the
elements of $Q$ are unrelated to the elements of $S$. We write
\[
    P=Q\vee S
\]
(although posets have no base points). If $(P,\f)$ is a monomial
functor and $P=Q\vee S$ (a disjoint union of poset subfunctors)
then $$\f=\f_Q\oplus\f_S$$ since upper triangular endomorphisms of
$PR$ have no cross terms ($QR$ does not map to $SR$ and $SR$ does
not map to $QR$). Therefore,
\[
    (P,\f)=(Q,\f_Q)\oplus(S,\f_S).
\]

\subsection{The Whitehead category}\label{sec1:subsec:Wh(R,G)}

We now need to construct two simplicial categories $\M\bu(R,G)$
and $\Wh\bu(R,G)$ which will serve as the target category for the
$(R,G)$-monomial functors defined above and for the more general
$(R,G)$-expansion functors, respectively. The arrows are reversed
from the definitions of the same terms in \cite{[I:BookOne]}. We
will see that this makes no problems.

\begin{defn}[Whitehead category]\label{sec1:def:Wh(R,G)}
The \emph{Whitehead category} $\Wh\bu(R,G)$ is defined to be the
simplicial category with simplicial set of objects equal to
$\Delta^\bullet(R)$. A morphism $(P,\f)\to(Q,\c)$ in $\Wh_k(R,G)$
is defined to be a decomposition $P=A\vee B$ and a morphism
$\a:A\to Q$ in $\nP(G)$ so that
\begin{enumerate}
  \item The restriction $\f_A$ of $\f$ to $A$ corresponds to $\c$ under the ring
  isomorphism
  \[
      \End_R(AR)\to \End_R(QR)
  \]
  induced by $\a$, i.e., $(\f_A)_\#=\c$.
  \item $B=S_1\vee \cdots \vee S_n$ where each $S_i=\{x_i^+,x_i^-\}$ so that
  $\f_0(j)(x_i^+)=x_i^-$ (i.e., $\f_0(j)(x_i^+,1)
  =(x_i^-,g)$ for some $g\in G$) for every $j\in[k]$.
\end{enumerate}
Let $\M\bu(R,G)$ denote the simplicial subcategory of
$\Wh\bu(R,G)$ having all of the objects but only \emph{monomial
morphisms} (i.e., those with $P=A$ and $B=\emptyset$). The pairs
$S_i$ are called \emph{collapsing pairs} in $P$. They are
independent by assumption.
\end{defn}

Morphisms in the Whitehead category factor uniquely:
\[
    (A\vee B,\f_A\oplus\f_B)\xrarrow{c}(A,\f_A)\xrarrow{\a}(Q,\c).
\]
The second morphism is monomial. The first is called a
\emph{collapse} (or \emph{inverse expansion}). (In the definition
above, the morphism is an collapse if $A=Q$ and $\a$ is the
identity morphism.) The induced mapping
\[
    \End_R(AR\oplus BR)\twoheadrightarrow \End_R(AR)\xrarrow{\approx}\End_R(QR)
\]
is not a ring homomorphism. The first map is given by
\[
    \begin{pmatrix}
        a&b\\
        c&d
    \end{pmatrix}
    \mapsto
    a.
\]
The second map is the isomorphism induced by $\a$. However, we do
get a homomorphisms of rings without unit on the strictly upper
triangular endomorphisms because, on that set, $b$ and $c$ in the
matrix above are zero.

\begin{rem}[reinterpreting monomial
functors]\label{sec1:rem:reinterpreting monomial functors} An
$(R,G)$-monomial functor $\xi=(P,\f)$ on a simplicial set $X$ can
now be reinterpreted as an actual functor
\[
    \xi:\simp X\to diag\simp \M\bu(R,G)
\]
from the category of simplices of $X$ to the diagonal of the
bicategory of simplices in the simplicial category $\M\bu(R,G)$.
It takes an object $(x,[k])$ of $\simp X$ to the $\Delta^k$ family
\[
    \xi(x,[k])=(P(x),\f_\#)
\]
given in (\ref{sec1:eq:Delta k family xi(x)}) above. For each map
$a:[j]\to[k]$ giving a morphism $a_\ast:(a^\ast(x),[j])\to(x,[k])$
in $\simp X$ we get a monomial morphism
\[
    (P(a^\ast(x)),\f_\#)\to a^\sharp(P(x),\f_\#)
\]
where $a^\sharp$ is the simplicial operator on $\Delta^\bullet(R)$
induced by $a$. Writing $y=a^\ast(x)$, we can simplify the
expression to
\[
    g_{yx}:\xi(y)\to \xi(x)|y.
\]
The functoriality of $\xi$ can be expressed by the equations:
\begin{equation}\label{sec1:eq:functoriality of xi}
  g_{zx}=g_{yx}|z \circ g_{zy},\quad g_{xx}=id_{\xi(x)}.
\end{equation}
To simplify the notation we write monomial functors simply as
\[
    \xi:\simp X\to \M\bu(R,G).
\]
\end{rem}

\begin{defn}[expansion functor]
An \emph{$(R,G)$-expansion functor} on a simplicial set $X$ is
defined to be a functor, written as
\[
    \xi:\simp X\to \Wh\bu(R,G),
\]
from $\simp X$ to the diagonal of $\simp\Wh\bu(R,G)$ which has the
form $\xi(x,[k])=(P(x),\f(x))\in\Delta^k(R)$ so that each map
$a:[j]\to[k]$ induces a morphism
\[
    g_{yx}:\xi(y,[j])=(P(y),\f(y))\to
    a^\sharp(P(x),\f(x))=\xi(x,[k])|y
\]
in $\Wh_j(R,G)$, where $y=a^\ast(x)$, satisfying
(\ref{sec1:eq:functoriality of xi}) above.
\end{defn}

Suppose that we have a fiberwise oriented GMF $f:E\to\RR$ and a
group homomorphism
\[
    \r:\pi_1E\to G\subseteq R^\times.
\]
Let $\simp B$ be the category of small smooth transverse simplices
in $B$. Then, disregarding for a moment the problem of
independence of b-d points, we get the \emph{cellular chain
complex functor}
\begin{equation}\label{sec1:eq:cellular chain complex functor}
  C_\ast(f;\r):\simp B\to \Wh\bu(R,G).
\end{equation}
This is the \emph{$(R,G)$-expansion functor} which sends each
small $k$-simplex $\s^k$ in $B$ to the cellular chain complex of
$f|\s$ and each face or degeneracy $(\t,[j])\to(\s,[k])$ to a
morphism
\[
    C_\ast(f;\r)(\t)\to C_\ast(f;\r)(\s)|\t
\]
in $\Wh_j(R,G)$. It induces a continuous mapping
\[
    |C_\ast(f;\r)|:B\simeq|\simp B|\to |\Wh\bu(R,G)|.
\]
Conceptually, $C_\ast(f;\r)$ represents a family of based free
$R$-complexes, one for each point in $B$, which locally varies by
$G$-monomial change of basis, elementary expansions, elementary
basis change and other higher chain homotopies. So expansion
functors
\[
    \xi:\simp B\to \Wh\bu(R,G)
\]
will also be called \emph{families of chain complexes} on $B$.

Let $\Wh\bu^h(R,G)$ be the simplicial full subcategory of
$\Wh\bu(R,G)$ in which all chain complexes are acyclic. Then we
have the following theorem. (Since any category is homotopy
equivalent to its opposite category, the theorem remains true even
though the arrows are reversed.)

\begin{thm}[Igusa-Klein \cite{[IK2:FR-torsion1]},
\cite{[I:BookOne]}]\label{sec1:thm:IK main theorem} There is a
natural homotopy fiber sequence
\[
    |\Wh\bu^h(R,G)|\to\Omega^\infty\Sig^\infty(BG_+)\to
    BGL(R)^+\times\overline{\ZZ}
\]
where $()^+$ is the Quillen plus construction, $()_+$ means add a
disjoint base point and $\overline{\ZZ}$ is the image of $\ZZ$ in
$K_0R$.
\end{thm}

\begin{rem}[$Q(BG_+)\simeq \Omega|\Sdot\Wh\bu(R,G)\op|$]\label{sec1:rem:QBG is the group completion of Wh(R,G)}
The middle term $\Omega^\infty\Sig^\infty(BG_+)$ is the
\emph{Waldhausen K-theory} or \emph{group completion} of (a
Waldhausen category homotopy equivalent to) $\Wh\bu(R,G)$ and the
map from $|\Wh\bu^h(R,G)|$ into it is induced by inclusion. The
first term $|\Wh\bu^h(R,G)|$ is already an infinite loop space so
is equal to its group completion.
\end{rem}

In the special case when $R=\ZZ,G=1$ we conclude that

\begin{cor}[rational homotopy of $\Wh\bu^h(\ZZ,1)$]\label{sec1:cor:homotopy type of Wh h(Z,1)}
\[
    \pi_i\Wh\bu^h(\ZZ,1)\otimes\QQ=
  \begin{cases}
    \QQ & \text{if $i=4k$ with $k>0$}, \\
    0 & \text{otherwise}.
  \end{cases}
\]
\end{cor}

This implies that the real cohomology ring of $\Wh\bu^h(\ZZ,1)$ is
a polynomial ring in generators
\[
    \t_{2k}\in H^{4k}(\Wh\bu^h(\ZZ,1);\RR).
\]
A formula for these classes is given (more generally) as follows.

\subsection{Definition of higher FR torsion}\label{sec1:subsec:def
of higher FR torsion}

First we need to reduce to the two-index case using the
\emph{two-index theorem}.

\begin{thm}[two-index theorem]\label{sec1:thm:2-index thm}
For any ring $R$, subgroup $G$ of $R^\times$ and integers $0\leq
i<j$, the simplicial full subcategory $\Wh\bu^{h[i,j]}(R,G)$ of
$\Wh\bu^h(R,G)$ consisting of acyclic $R$-complexes with basis
posets having elements only in degrees $i,i+1,\cdots,j$ is a
deformation retract.
\end{thm}

Consequently, we have:
\[
    \Wh\bu^{h[0,1]}(R,G)\simeq\Wh\bu^h(R,G).
\]
We may regard $\Wh\bu^{h[0,1]}(R,G)$ as a space of invertible
matrices with coefficients in $R$ which are only well-defined up
to left and right multiplication by $G$-monomial matrices.

Now restrict to the case $(R,G)=(M_r(\CC),U(r))$. Then
\begin{equation}\label{sec1:eq:Wh h 01}
  \Wh\bu^{h[0,1]}(R,G)=\Wh\bu^{h[0,1]}(M_r(\CC),U(r))
\end{equation}
is a space of invertible complex matrices which are well-defined
up to left and right multiplication by unitary matrices. We modify
the definition so that every $k$-simplex is actually equal to a
smooth family of invertible matrices over $\CC$ and morphisms
multiply these by unitary matrices. We will come back to this
point in the next subsection.

Next, we use the \emph{Kamber-Tondeur form}
(\cite{[Kamber-Tondeur74]}, \cite{[Dupont76]}). This gives a
simplicial $2k$-cochain $c_{2k}$ on the smooth version of the
simplicial category (\ref{sec1:eq:Wh h 01}) given on a smooth
$2k$-simplex $f_t,t\in\Delta^k$, by
\begin{equation}\label{sec1:eq:Kamber-Tondeur}
  c_{2k}(f_t)=\frac1{2i^k(2k+1)!}\int_{(t,u)\in \Delta^{2k}\times
    I}\Tr((h_t^{-u}dh_t^u)^{2k+1})
\end{equation}
where $h_t=f_tf_t^\ast$. By adding certain polynomial correction
terms to $c_{2k}$ we get a cocycle $D_{2k}$. We do not need to
know what these correction terms are because the higher FR torsion
is a rational multiple of a polylogarithm and polylogarithms are
linearly independent from polynomials. (\cite{[I:BookOne]}, Lemma
7.7.2.)

The result is that we have well-defined real cohomology classes
\begin{equation}\label{sec1:eq:universal FR torsion invariants}
  \t_k=[D_{2k}]\in H^{2k}(\Wh\bu^h(M_r(\CC),U(r));\RR)
\end{equation}
which we call the \emph{universal Franz-Reidemeister torsion
invariants}. (It is independent of $r$ is a natural sense.) To get
invariants of smooth bundles we proceed as follows.

Let $(E,\d_0E)\to B$ be a smooth bundle pair with fiber
$(M,\d_0M)$ and let $\F$ be a Hermitian coefficient system on $E$.
This is given by a unitary representation $\r:\pi_1E\to U(r)$.
Suppose further that
\begin{equation}\label{sec1:eq:acyclicity assumption}
  H_\ast(M,\d_0 M;\F)=0.
\end{equation}
Suppose there is a fiberwise framed function
\[
    f:(E,\d_0E)\to (I,0)
\]
which is unique up to framed homotopy. The Framed Function Theorem
(\ref{sec1:thm:framed function theorem}) says this holds if $\dim
M>\dim B$. Then, setting aside a problem with the b-d points, the
cellular chain complex functor (\ref{sec1:eq:cellular chain
complex functor}) of $f$ with coefficients in $\F$ is an expansion
functor
\[
    C_\ast(f;\F):\simp B\to\Wh\bu^h(M_r(\CC),U(r)).
\]
So, we can pull back the universal higher FR-torsion invariant
$\t_k$ to obtain a well-defined cohomology class
\[
    \t_k(E,\d_0E;\F)=C_\ast(f;\F)^\ast(\t_k)\in H^{2k}(B;\RR).
\]
This is called the \emph{higher Franz-Reidemeister torsion
invariant} of the pair $(E,\d_0E)$ with coefficients in $\F$.

\subsection{Families of matrices as flat
superconnections}\label{sec1:subsec:superconnections}

We will explain what we mean by a smooth family of invertible
matrices. By differentiating the components of the structure we
will obtain what is known as a flat $\ZZ$-graded superconnection
of total degree 1 \cite{[Bismut-Lott95]}. (See also Quillen's
original article on superconnections
\cite{[Quillen:superconnections]}.) This is explained in detail in
\cite{[I:BookOne]} although the superconnection formalism is not
mentioned there. The relation to superconnections was observed by
Goette \cite{[Goette01]}.

\begin{defn}[smooth category of invertible
matrices]\label{sec1:def:W diff(Cr,U(r))} The \emph{smooth
category of invertible matrices} is defined to be the simplicial
category $\W\bu^{dif\!f}(\CC^r,U(r))$ whose objects in degree $k$
are triples $(P,C_\ast,E)$ where
\begin{enumerate}
  \item $P=P_0\coprod P_1$ is a graded poset with $|P_0|=|P_1|=n$
  for some $n\geq0$.
  \item $C_\ast(t)=(P\CC\otimes\CC^r,\f_0(t)), t\in\Delta^k$ is
  a smooth $\Delta^k$-family of acyclic chain complexes of the
  form
  \[
     P_1\CC\otimes\CC^r\xrarrow{\f_0(t)}P_0\CC\otimes\CC^r
  \]
  where the boundary map $\f_0(t)$ is strictly upper triangular
  and gives a smooth mapping
  \[
    \f_0:\Delta^k\to GL_{nr}(\CC).
  \]
  \item $E(s,t)$ is a smooth family of upper triangular chain
  isormorphisms
  \[
    E(s,t):C_\ast(t)\xrarrow{\approx}C_\ast(s)
  \]
  for $(s,t)\in\Delta^k\times\Delta^k$ so that $E(t,t=I$ (the
  identity map) and $E(s,t)-I$ is strictly upper triangular for
  all $s,t$.
  \item There exists a unique family of strictly upper triangular
  chain homotopies
  \[
    S(s,t,u):E(s,u)\simeq E(s,t)E(t,u)
  \]
  for $(s,t,u)\in\Delta^k\times\Delta^k\times\Delta^k$ so that
  $S(s,t,u)=0$ if $s,t,u$ are not distinct.
\end{enumerate}
The morphisms will be given later.
\end{defn}

\begin{rem}[homotopy equivalence with the Whitehead category]
The simplicial category $\W\bu^{dif\!f}(\CC^r,U(r))$ and its
continuous counterpart $\W\bu^{cont}(\CC^r,U(r))$ are both
homotopy equivalent to the combinatorial space
$\Wh\bu^{h[0,1]}(M_r(\CC),U(r))$. This is Theorem 2.5.11(a) in
\cite{[I:BookOne]}. The homotopy equivalence is given by forgetful
functors:
\[
    \W\bu^{dif\!f}(\CC^r,U(r))\xrarrow{\simeq}
    \W\bu^{cont}(\CC^r,U(r))\xrarrow{\simeq}
    \Wh\bu^{h[0,1]}(M_r(\CC),U(r)).
\]
By the two index theorem (\ref{sec1:thm:2-index thm}) these are
all homotopy equivalent to the Whitehead category
$\Wh\bu^h(M_r(\CC),U(r))$.
\end{rem}

A smooth family of chain complexes $(P,C_\ast)$ is equivalent to a
flat $\ZZ$-graded superconnection on $P\CC\otimes\CC^r$ of total
degree $1$. This was observed by Goette \cite{[Goette01]}. The
idea is to differentiate the structure maps. In our case there are
only three structure maps $\f_0(t),E(s,t)$ and $S(s,t,u)$.

Since $E(s,t)$ gives a chain isomorphisms $C_\ast(t)\to C_\ast(s)$
we have:
\[
    E(s,t)\f_0(t)=\f_0(s)E(s,t).
\]
Differentiating both sides with respect to $s$ and put $s=t$ we
get
\begin{equation*}
  (d_1E)\f_0=d\f_0+\f_0d_1E
\end{equation*}
where $d_1$ means differential with respect to the first variable
(so $d_1E+d_2E=0$ since $E(t,t)=I$ is constant). If we denote the
matrix $1$-form $d_1E$ by $\f_1$ we get
\begin{equation}\label{sec1:eq:d phi 0}
  d\f_0=\f_1\f_0-\f_0\f_1.
\end{equation}

The chain maps $E(s,t)$ are only required to form a functor up to
homotopy by $S(s,t,u)$:
\[
    E(s,t)E(t,u)=E(s,u)+\f_0(s)S(s,t,u)+S(s,t,u)\f_0(u).
\]
When this equation is differentiated with respect to $s$ and $t$
and the variables are equated, we get:
\[
    d_1d_2E+d_1E\,d_1E=\f_0d_1d_2S+(d_1d_2S)\f_0.
\]
Denoting the matrix 2-form $d_1d_2S$ by $\f_2$ this becomes:
\begin{equation}\label{sec1:eq:d phi 1}
    -d\f_1=\f_0\f_2-\f_1^2+\f_2\f_0.
\end{equation}

\begin{lem}[derivative of top form]\label{sec1:lem:derivative of
top form} The following equation follows from (\ref{sec1:eq:d phi
0}) and (\ref{sec1:eq:d phi 1}).
\begin{equation}\label{sec1:eq:d phi 2}
  d\f_2=\f_1\f_2-\f_2\f_1.
\end{equation}
\end{lem}

\begin{proof}
Differentiate (\ref{sec1:eq:d phi 1}) and replace $d\f_0$ and
$d\f_1$ by their values as given in (\ref{sec1:eq:d phi 0}) and
(\ref{sec1:eq:d phi 1}). The result can then be written as
\[
    Z\f_0=\f_0Z
\]
where
\[
    Z=d\f_2-\f_1\f_2+\f_2\f_1.
\]
But $Z$ maps $C_0$ to $C_1$ so we must have $Z\f_0=\f_0Z=0$ which
implies that $Z=0$.
\end{proof}

A $\ZZ$-graded superconnection $D=d+A_0+A_1+A_2$ of total degree 1
is \emph{flat} if
\[
    -dA=A^2.
\]
However, Quillen's notation \cite{[Quillen:superconnections]}
places the forms on the left of the matrices. If we use the
supercommutator rules to put $A^2$ in \emph{normal form} by moving
both forms to the left of both matrices (See \cite{[I:Twisted]}
for more details about this point) we get the following.
\begin{enumerate}
  \item $dA_0=A_0A_1-A_1A_0$
  \item $-dA_1=A_0A_2+A_1^2+A_2A_0$
  \item $dA_2=-A_1A_2+A_2A_1$
\end{enumerate}
This matches (\ref{sec1:eq:d phi 0}), (\ref{sec1:eq:d phi 1}) and
(\ref{sec1:eq:d phi 2}) if we let $A_0=\f_0,A_1=-\f_1$ and
$A_2=-\f_2$. So, we get the following.

\begin{prop}[families of matrices give superconnections]
\[D=d+\f_0-\f_1-\f_2\] is a flat superconnection on
$\Delta^k$.
\end{prop}

\begin{rem}[notation change]\label{sec1:rem:notation of
superconnection} In \cite{[I:BookOne]} the notation is different.
We used the symbols $X,U,V,W,f$ to denote:
\begin{enumerate}
  \item $f=\f_0$. (This is an isomorphism
  $C_1\xrarrow{\approx}C_0$.)
  \item $(U,V)=\f_1$. (So $U=\f_1|C_0$, $V=\f_1|C_1$.)
  \item $W=\f_2\f_0$. (Or $\f_2=Wf^{-1}$.)
  \item $X=\f_0^{-1}d\f_0$. (Or $d\f_0=fX$.)
\end{enumerate}
\end{rem}

\begin{rem}[role of superconnections in definition of higher FR
torsion]\label{sec1:rem:role of X,V,W in def of D2k} The
$2k$-cocycle $D_{2k}$ which gives the higher torsion is equal to
$C_{2k}$ given in (\ref{sec1:eq:Kamber-Tondeur}) plus a linear
combination of integrals of traces of products of $X,V,W$ which
are endomorphisms of $C_1$. (The term $U$ is an endomorphism of
$C_0$ and occurs only in the form $f^{-1}Uf=-X-V$.)
\end{rem}

Finally, we give the definition of the morphisms in the category
$\W\bu^{dif\!f}(\CC^r,U(r))$. They are just as in the Whitehead
category. A morphism $(P,\f)\to(Q,\c)$ is a decomposition $P=A\vee
B$ together with a morphism $\a:A\to Q$ in $\nP(U(r))$ so that
\begin{enumerate}
  \item For all $t\in\Delta^k$, the restriction $\f_A(t)$ of $\f(t)$ to $A$ corresponds to $\c$
  under the isomorphism
  \[
    \Omega(\Delta^k,\End(A\CC\otimes\CC^r))\cong
    \Omega(\Delta^k,\End(Q\CC\otimes\CC^r))
  \]
  induced by $\a$.
  \item $B=S_1\vee\cdots\vee S_n$ where each $S_i$ is a collapsing
  pair $x^+,x^-$ with constant $\f_0(t)\in U(r)$ and trivial $\f_1$ and
  $\f_2$.
\end{enumerate}

\subsection{Independence of birth-death
points}\label{sec1:subsec:independence of bd points}

We will discuss briefly a problem with the b-d points and how it
is solved. The problem is that the definition of an expansion
functor
\[
    \simp B\to\Wh\bu(R,G)
\]
assumes that collapsing pairs are independent. The independence of
b-d points would automatically imply this algebraic condition for
the cellular chain functor constructed above. However, even if b-d
points are not independent, the condition can still be achieved in
a canonical algebraic way.

A b-d point $x$ of $M_t$ is called \emph{independent} if $x$ is
not incident to any other critical point of $f_t$. In other words,
$\{x\}$ is an independent set in the poset $P^+(t)=\Sig(f_t)$.
Suppose for a moment that all b-d points are independent. Since
independence is an open condition, this would imply that, for any
small smooth transverse simplex $\s$ in $B$, any component $x\in
P^+(\s)-P(\s)$ (i.e., $x$ is a component of $\Sig(f|\s)$
containing a b-d point), will be independent in the poset
$P^+(\s)$.

Since $x\notin P(\s)$, it does not directly appear as part of the
expansion functor $\xi(\s)=(P(\s),\f_\sharp(\s))$. However, it may
produce a collapsing pair on a face $\t$ of $\s$. This happens if
$x|\t$ has no b-d points. Then $x|\t$ is necessarily a disjoint
union of two Morse components $x^+,x^-$ which form a collapsing
pair. The independence of $x$ implies that $S=\{x^+,x^-\}$ is
independent in $P(\t)$ as required by
Definition~\ref{sec1:def:Wh(R,G)}.

The condition of independence of b-d points does not hold in
general. For example, if $M_t$ is a closed manifold, it is
impossible for b-d points to be independent. One way to get around
this problem is by a combination of positive and negative
suspension. However, this produces manifold bundles with corners
which complicates the geometry. The method that we prefer is
algebraic.

In \cite{[I:BookOne]}, a rather complicated explicit algebraic
formula is given to deform a nonindependent collapsing pair into
an independent one. The reason it is complicated is because this
deformation must be carried out simultaneously and in a continuous
way on all collapsing pairs which appear in the expansion functor
$\xi$.

Here, we will use the local splitting lemma
(Lemma~\ref{sec2:lem:local splitting lemma} below) which tells us
that such a deformation exists and is unique up to homotopy since
collapsing pairs form an acyclic subquotient complex.

The preferred solution to our problem is as follows. We do not
assume independence of b-d points. Instead we choose the twisted
cochains on the simplices $\t$ which are faces of some $\s$ as
described above so that the collapsing pairs are algebraically
independent. Then we alter the partial ordering on $P(\t)$ so that
each collapsing pair is poset independent.

We also need to alter the partial ordering of $P^+(\s)$ in such a
way that collapsing pairs and ``expanded pairs'' act in unison.
Here, a \emph{collapsing pair} in $P^+(\t)$ is a pair of Morse
components $x^+,x^-$ which collapse (to one b-d component) in
another simplex $\s$ of which $\t$ is a face. A pair of Morse
components $x^+,x^-$ in $P(\s)$ which restricts to a collapsing
pair on some face $\t$ is an \emph{expanded pair}. To make these
pairs \emph{act in unison} we increase the partial ordering so
$x^-<y$ implies $x^+\leq y$ and $z<x^+$ implies $z\leq x^-$. The
poset $P^+(\s)$ is used to construct the filtered chain complex of
$f$ over $\s$. We are just discarding those subcomplexes which
separate the members of a collapsing or expanded pair, i.e., we
combined these into one layer.

With this solution of the independence of b-d point problem, the
only topological problem we have in the definition of higher FR
torsion is the existence of fiberwise framed functions. For the
moment, we assume that the higher FR-torsion is only defined when
the fiber dimension is large. Later, using properties of the
higher torsion, we will be able to define the higher torsion of
$(E,\d_0E)$ to be half the torsion of $(E\times S^{2N},\d_0\times
S^{2N})$ for $N$ large. (See Remark~\ref{sec2:rem:the product
formula is a key lemma}.)

In the absolute case, when $\d_0E$ is empty ($\d E=\d_1E$), we can
define the higher FR-torsion of $E$ to be that of $E\times D^N$
(with corners rounded). This is well-defined by the positive
suspension lemma proved below.

\subsection{Positive suspension lemma}\label{sec1:subsec:pos
suspension lemma}

Suppose that
\[
    f:(E,\d E)\to (I,1)
\]
is a fiberwise oriented GMF. Let $C_\ast(f;\r)$ be the
$(R,G)$-expansion functor for $f$ with coefficients in a
representation $\r:\pi_1E\to G\subseteq R^\times$.

\begin{defn}[positive suspension]\label{sec1:def:positive suspension}
Let $\pi:D\to E$ be the unit disk bundle of an $m$-plane bundle on
$E$. Then the \emph{positive suspension}
\[
    h=\s_+f:D\to\RR
\]
of $f$ is the map given by
\[
    h(y)=f(x)+\length{y}^2
\]
where $x=\pi(y)$.
\end{defn}

Then $\s_+f$ is an oriented fiberwise GMF with the same cellular
chain complex as $f$. We will prove this obvious statement. The
proof will explain something about the definition of $C(f;\r)$ and
will also give the definition of higher torsion for low
dimensional fibers as promised at the end of the last subsection.

\begin{lem}[positive suspension lemma]\label{sec1:lem:positive
suspension lemma} The $(R,G)$-expansion functors of $f$ and
$\s_+f$ are equal.
\[
    C_\ast(\s_+f;\r)=C_\ast(f):\simp B\to\Wh\bu(R,G).
\]
\end{lem}

\begin{rem}[vertical metric on $D$]\label{sec1:rem:vertical metric
on D} To define $C_\ast(\s_+f;\r)$ we need a vertical metric on
$D$ as a bundle over $B$. As an $O(m)$-disk bundle, $D$ has a
vertical metric over $E$. Choose a smooth horizontal distribution
for $D$ over $E$ given by an $O(m)$-connection. Then the
horizontal directional derivatives of $\length{y}^2$ will be zero.
At each point $y\in D$ let $K_y$ be the kernel of the composition
\[
    T_yD\xrarrow{D\pi}T_xE\xrarrow{Dp}T_{p(x)}B
\]
restricted to the horizontal plane at $y$. Then $K_y$ maps
isomorphically to the vertical tangent plane of $E$ at $x$:
\begin{equation}\label{sec1:eq:K iso TvE}
  D\pi:K_y\xrarrow{\approx}T^v_xE.
\end{equation}
Pull back the vertical metric on $E$ (the one which was used to
define $C(f;\r)$) and make $K_y$ perpendicular to the fiber of $D$
over $x$. This defines a vertical metric for $D$ over $B$.
\end{rem}

\begin{proof} By definition of the vertical metric of $D$ over
$B$, the gradient of $h=\s_+f$ decomposes into vertical and
horizontal components
\begin{equation}\label{sec1:eq:grad s+(f)}
  \grad h(y)=2y+\pi^\ast\grad f(x)
\end{equation}
where $\pi^\ast$ is the inverse of $D\pi$ in (\ref{sec1:eq:K iso
TvE}). Therefore, the critical points of $h$ lie in the zero
section of $E$ in $D$ and coincide with the critical points of $f$
(if we identify $E$ with its zero section in $D$). Furthermore,
$\grad h=\grad f$ on this zero section. Thus the inclusion of the
zero section is filtration preserving and thus induces a morphism
of filtered chain complexes from the total singular complex of
$E|\s$ to that of $D|\s$ for every small transverse simplex $\s$
in $B$. By looking at the layers, we see that this inclusion map
is a filtered quasi-isomorphism preserving the dual classes. So
$C(f;\r)$ and $C_\ast(h;\r)$ can be taken to be equal. (They are
only well defined up to simplicial homotopy.)
\end{proof}

This lemma can be used to define the cellular chain complex
functor on any smooth bundle $E\to B$ with $\d_0E=\emptyset$.

\begin{thm}[cellular chain complex for $\d_0E$
empty] \label{sec1:thm:def of cellular chain complex for d0E
empty} If $E\to B$ is a smooth bundle and $\r$ is a representation
of $\pi_1E$ in $G\subseteq R^\times$ then we get an
$(R,G)$-expansion functor
\[
    \xi:\simp B\to\Wh\bu(R,G)
\]
which is well-defined up to homotopy by taking $\xi=C_\ast(f;\r)$
where $f$ is a fiberwise framed function on $E\times D^N$ for $N$
large.
\end{thm}

\begin{proof}
The statement is that $\xi$ does not depend on $N$ and $f$ (up to
homotopy). So suppose that $D$ is another disk bundle over $E$ and
$g:D\to\RR$ is a fiberwise framed function whose gradient point
outward on all faces of the boundary. Then we get two positive
suspensions
\[
    \s_+f,\s_+g:D\times D^N\to \RR.
\]
By the Framed Function Theorem, these are homotopic, so
\[C_\ast(f;\r)=C_\ast(\s_+f;\r)\simeq
    C_\ast(\s_+g;\r)=C_\ast(g;\r).
\]
\end{proof}

%% file: sec2.tex
\vfill\eject\section{Properties of higher FR
torsion}\label{sec2:section 2}

\begin{enumerate}
   \item Basic properties
   \item Suspension theorem
  \item Additivity, Splitting Lemma
  \item Applications of the Splitting Lemma
  \item Local equivalence lemma
  \item Product formula
  \item Transfer for coverings
  \item More transfer formulas
\end{enumerate}

We go over the basic properties of higher Franz-Reidemeister
torsion. Properties stated without proof are proved in
\cite{[I:BookOne]}.

\subsection{Basic properties}\label{sec2:subsec:basic properties}

First we state the obvious.

\begin{prop}[torsion of trivial bundle]\label{sec2:prop:torsion of
trivial bundle} If $E=\d_0E\times I$ then the relative torsion
with respect to any coefficient system is defined and equal to
zero:
\[
    \t_k(\d_0E\times I,\d_0E;\F)=0.
\]
\end{prop}

Another obvious fact is the naturality of torsion.

\begin{prop}[naturality of torsion]\label{sec2:prop:naturality of
torsion} The higher FR torsion invariants are natural on those
smooth bundles $E\to B$ on which they are defined, i.e., if $E\to
B$ is a smooth bundle and $f:B'\to B$ is a smooth mapping then the
torsion of the pull-back $f^\ast(E)$ of $E$ is the pull-back of
the torsion:
\[
    \t_k(f^\ast(E),f^\ast(\d_0E);f^\ast(\F)=f^\ast(\t_k(E,\d_0E;\F))
\]
assuming that $\t_k(E,\d_0E;\F)$ is defined.
\end{prop}

The naturality of higher FR torsion extends to the case when the
coefficient sheaf $f^\ast(\F)$ is replaced by a locally equivalent
sheaf. When we say that two Hermitian coefficient systems $\F$,
$\F'$ on $E$ are \emph{locally equivalent} we mean that $\F$ and
$\F'$ are unitarily equivalent on $E|\s$ for every small simplex
$\s$ in $B$. For example, if $\F$ is the pull-back of a
coefficient system on $B$ then it is locally equivalent to the
trivial coefficient sheaf with the same fiber. We only need this
in the case when the fiber is acyclic.

\begin{prop}[locally equivalent
coefficients]\label{sec2:prop:locally equivalent coefficients}
Suppose that $\F$ and $\F'$ are locally equivalent Hermitian
coefficient systems on $E$ and $H_\ast(M;\F)=H_\ast(M;\F')=0$.
Then
\[
    \t_k(E;\F)=\t_k(E;\F').
\]
\end{prop}

\begin{proof}
The reason is simple. Higher FR torsion is defined locally and is
invariant under unitary change of bases. First, we have to deform
the cellular chain complex functor into two indices. This can be
done over each simplex of $B$. Thus, the results will still be
locally equivalent. This means we have families of complex
matrices which are conjugate by unitary matrices, i.e.,
$g_t'=Ug_tV$ where $U,V$ are fixed unitary matrices (depending
only on the simplex $\s$).

But the higher FR torsion is given by integrals of traces of
products of $g_t$, $g_t^\ast$ and their inverses and derivatives
($d(g_tg_t^\ast)^u$). These are always arranged in a natural way
so that the unitary matrices $U$, $V$ cancel and do not appear
when we replace $g_t$ with $g_t'$. Thus the trace terms are
unchanged. So their integrals are equal.
\end{proof}

The following lemma allows us to generalize the condition
(\ref{sec1:eq:acyclicity assumption}) under which $\t_k$ is
defined. We need one definition. If $K$ is a field and $\pi$ is a
group, we say that a finite dimensional $K[\pi]$-module $M$ is
\emph{upper triangular} if it contains a filtration by
$K[\pi]$-submodules (i.e., a composition series)
\[
    M=M_n\supseteq M_{n-1}\supseteq\cdots\supseteq 0
\]
so that $\pi$ acts trivially on the subquotients $M_i/M_{i-1}$.

\begin{lem}[canonical cone]\label{sec2:lem:canonical cone}
Let $\xi$ be a family of chain complexes
\[
    \xi:B\to\Wh\bu(R,G)
\]
satisfying either of the following two conditions where $v\in B$
is any base point.
\begin{enumerate}
  \item The homology groups $H_n(\xi(v);R)$ are projective
  $R$-modules with trivial $\pi_1B$ action.
  \item $R$ is a field and each $H_n(\xi(v);R)$ is an upper triangular $\pi_1B$-module.
\end{enumerate}
Then there exists a family of acyclic chain complexes
\[
    C(\xi):B\to\Wh\bu^h(R,G)
\]
containing $\xi$ as a subfunctor which is canonically defined up
to simplicial homotopy.
\end{lem}

We call $C(\xi)$ the \emph{canonical cone} of $\xi$.

\begin{thm}[when $\t_k$ is defined]\label{sec2:thm:when higher FR torsion is defined}
Suppose that $(E,\d_0E)\to B$ is a smooth bundle with fiber
$(M,\d_0M)$ and let $\F$ be an $r$-dimensional Hermitian
coefficient system over $E$. Suppose also that either
\begin{enumerate}
  \item $\pi_1B$ acts trivially on $H_\ast(M,\d_0M;\F)$ or
  \item $r=1$ and $H_n(M,\d_0M;\F)$ is an upper triangular
  $\pi_1B$-module for all $n\geq0$.
\end{enumerate}
Then the higher Franz-Reidemeister torsion invariants
\[
    \t_k(E,\d_0E;\F)\in H^{2k}(B;\RR)
\]
are defined.
\end{thm}

\begin{rem}[how $\t_k$ is defined]\label{sec2:rem:how higher FG torsion is defined}
These higher torsion invariants are given by pulling back the
universal invariants (\ref{sec1:eq:universal FR torsion
invariants}) along the canonical cone
\[
    C(C_\ast(f)):B\to\Wh\bu^h(M_r(\CC),U(r))
\]
of the cellular chain complex $C_\ast(f)$ associated to a
fiberwise framed function \[f:(E,\d_0E)\to(I,0).\] If $f$ does not
exist we can take the product with a large even dimensional sphere
and divide by $2$:
\[
    \t_k(E,\d_0E;\F)=\frac12\t_k(E\times S^{2N},\d_0E\times
    S^{2N};\pi^\ast\F).
\]
Here $\pi^\ast(\F)$ is the pull-back along $\pi:E\times S^{2N}\to
E$ of the Hermitian coefficient system $\F$. (See
Lemma~\ref{sec2:lem:product formula} below.)
\end{rem}

Another basic fact that we need in this paper is the
\emph{involution property}. This relates the higher FR torsion
invariant given by a fiberwise GMF $f$ with that of its negative
$-f$. We need to assume that $E$ is \emph{fiberwise oriented} in
the sense that $M$ is an oriented manifold and the action of
$\pi_1B$ on $M$ preserves the orientation (otherwise, $-f$ would
not be fiberwise oriented).

\begin{thm}[Involution property]\label{sec2:thm:involution
property} Suppose that $E\to B$ is a fiberwise oriented smooth
bundle and $\F$ is a Hermitian coefficient system on $E$. Let
\[
    f:(E,\d_0E,\d_1E)\to(I,0,1)
\]
be any fiberwise oriented GMF. Then the higher torsion invariants
associated to $f$ and $-f$, if defined, differ by the sign
$(-1)^{n-1}$ where $n=\dim M$:
\[
    \t_{k}(CC_\ast(f);\F)=(-1)^{n-1}\t_{k}(CC_\ast(-f);\F).
\]
\end{thm}

\begin{proof} We outline the proof given in \cite{[I:BookOne]}, Theorem
5.8.1. The proof has two parts
\begin{enumerate}
  \item First consider the case $n=1$. Geometrically, this
  corresponds to the case when $M$ is a circle. Algebraically,
  we can imagine that we have a family of invertible matrices
  parametrized by $B$. The operation of \emph{inverse conjugate
  transpose} is an involution on $GL_n(\CC)$ which reverses the
  sign of the Kamber-Tondeur form. Since inversion acts as $-1$ on additive cohomology classes, \emph{conjugate
  transpose} does not change the sign of the higher torsion so the
  theorem holds for $n=1$.
  \item In general, we use the two-index theorem. This puts the
  cellular chain complex of $f$ into degrees $0,1$ and that of
  $-f$ into degrees $n,n-1$. Then we desuspend the cellular chain complex of $-f$ $n-1$ times to
  reduce to the $n=1$ case. By the suspension theorem this changes
  the sign by $(-1)^{n-1}$.
\end{enumerate}
\end{proof}

\subsection{Suspension Theorem}\label{sec2:subsec:suspension
theorem}

One of the fundamental properties of the higher torsion is that it
anti-commutes with suspension. We will see later that this is
geometrically obvious. (See Corollary~\ref{sec2:cor:easy geometric
suspension}.)

Given a chain complex $(C_\ast,\d)$, the ``alternating''
\emph{suspension} $(\Sig C,\d^\Sig)$ is given by $(\Sig
C)_n=C_{n-1}$ and $\d_n^\Sig=\d_{n-1}$. We do not alternate the
signs as in the usual definition. This version of the suspension
map induces simplicial functors
\[
    \Sig:\Wh\bu(R,G)\to\Wh\bu(R,G),\quad \Sig:\Wh\bu^h(R,G)\to\Wh\bu^h(R,G)
\]
(If $-1\in G$, the usual suspension operation will also induce
such a functor and the two suspension functors will be related by
a natural transformation and thus be homotopic.)

\begin{thm}[Suspension Theorem]\label{sec2:thm:suspension theorem}
The functor
\[
    \Sig+id:\Wh\bu^h(R,G)\to\Wh\bu^h(R,G)
\]
sending $\xi$ to $\Sig\xi\oplus\xi$ is null homotopic.
\end{thm}

Algebraic suspension is given by geometric negative suspension.
This is simply the upside down version of a positive suspension.

Suppose that
\[
    f:(E,\d E)\to (I,0)
\]
is a fiberwise oriented GMF. Let $C_\ast(f;\r)$ be the
$(R,G)$-expansion functor for $f$ with coefficients in a
representation $\r:\pi_1E\to G\subseteq R^\times$.

\begin{defn}[negative suspension]\label{sec2:def:negative suspension}
Let $\pi:D\to E$ be the unit disk bundle of an oriented $m$-plane
bundle on $E$. Then the \emph{negative suspension}
\[
    h=\s_-f:D\to\RR
\]
of $f$ is the map given by
\[
    h(y)=f(x)-\length{y}^2
\]
where $x=\pi(y)$.
\end{defn}

\begin{lem}[negative suspension lemma]\label{sec2:lem:negative
suspension lemma}  $\s_-f$ is an oriented fiberwise GMF whose
cellular chain complex functor is the $m$-fold suspension of
$C_\ast(f;\r)$.
\[
    C_\ast(\s_-f;\r)=\Sig^m C_\ast(f):\simp B\to\Wh\bu(R,G).
\]
\end{lem}

\begin{proof}
The proof is the same as for the positive suspension lemma. The
degree shift comes from the Thom isomorphism theorem.

We need to use the fact that the Thom isomorphism is given
functorially at the chain level by cap product with the Thom
class. The total singular complex $C_\ast$ of $D|\t$ together with
the system of subcomplexes given by all faces of $\t$ and closed
subsets $Q^+$ of $P^+(\t)$ is filtered chain homotopic to the
subcomplex $C_\ast'$ generated by simplices transverse to the zero
section. On this subcomplex we have the Thom isomorphism given as
a degree shifting quasi-isomorphism
\[
    C_\ast'(D|\t,(\d D)|\t;\r)\simeq\Sig^mC_\ast(E|\t,(\d E)|\t;\r)
\]
which preserves the filtration structure and induces a filtered
equivalence.
\end{proof}

If the disk bundle $D$ is trivial ($D=E\times D^m$) then the
negative suspension of any framed function will be framed. Thus we
can combine the above statements with the additivity property
(Theorem~\ref{sec2:thm:additivity wrt direct sum} below) to give
the following.

\begin{thm}[geometric suspension theorem]\label{sec2:thm:geometric
suspension}
\[
    \t_k(E\times D^n,\d(E\times D^n);\pi^\ast\F)=(-1)^m\t_k(E,\d E;\F)
\]
when the terms are defined.
\end{thm}

\subsection{Additivity, Splitting
Lemma}\label{sec2:subsec:splitting lemma}

Higher FR torsion satisfies two additivity properties. One comes
from coefficients and the other from filtrations of chain
complexes.

By additivity of coefficients we mean:
\[
    \t_k(E;\F_1\oplus\F_2)=\t_k(E;\F_1)+\t_k(E;\F_2).
\]
This is a special case of the following theorem.

\begin{thm}[additivity]\label{sec2:thm:additivity wrt direct sum}
Let $\xi,\eta$ be two families of acyclic based chain complexes
\[
    \xi,\eta:B\to\Wh\bu^h(M_n(\CC),U(n))
\]
Then the higher torsion of their direct sum $\xi\oplus\eta$ is
equal to the sum of their torsions:
\[
    \t_k(\xi\oplus\eta)=\t_k(\xi)+\t_k(\eta).
\]
\end{thm}

These two properties follow from the additivity of the
Kamber-Tondeur form under direct sum of matrices: If $f=f_1\oplus
f_2$ then $h=h_1\oplus h_2$ and
\[
    \Tr((h^{-1}dh)^{2k+1})=\Tr((h_1^{-1}dh_1)^{2k+1})+\Tr((h_2^{-1}dh_2)^{2k+1}).
\]
The polynomial correction terms are also additive with respect to
direct sum.

One easy corollary of this is the following.

\begin{cor}[torsion of disjoint union]\label{sec2:cor:torsion of
disjoint union} The higher torsion of a disjoint union is the sum
of the torsions, i.e.,
\[
    \t_k(E_1\smalldisj
    E_2,\d_0E_1\smalldisj\d_0E_2;\F_1\smalldisj\F_2)
    =\t_k(E_1,\d_0E_1;\F_1)+\t_k(E_2,\d_0E_2;\F_2),
\]
assuming the terms are defined.
\end{cor}

The following basic lemma allows us to extend the additivity of
torsion to filtered families of chain complexes. The proof is
trivial but the lemma is difficult to state. Let $(P,\f)\in
\Delta^k(R)$ be a $\Delta^k$-family of chain complexes and suppose
that $Q$ is a closed subset of $P$. Then the endomorphisms
$\f(\s)$ will not send any element of $Q$ to any element of the
complement $P/Q$. Thus it can be written as
\[
    \f(v_\ast)=
    \begin{pmatrix}
        \f_Q(\s) & f(\s)\\
        0 & \c(\s)
    \end{pmatrix}
\]
where $\c=\f/\f_Q$ is the quotient twisted cochain. Note that
$f(v_\ast)$ is a homomorphism
\[
    f(\s):(P/Q)R\to QR.
\]

\begin{lem}[local splitting lemma]
\label{sec2:lem:local splitting lemma} Suppose that every element
of $Q$ is less than every element of $P/Q$. Suppose that either
the subcomplex $(Q,\f_Q)$ or quotient complex $(P/Q,\c)$ is
acyclic. Suppose also that $f(\s)=0$ for all proper faces $\s$ of
$\Delta^k$. Then there is a simplicial homotopy of $(P,\f)$ which
is the identity on $Q$, on $P/Q$ and on all proper faces $\s$ of
$\Delta^k$ and so that, at the end, $f(\Delta^k)$ is also zero:
\[
    (P,\f)\simeq(Q,\f_Q)\oplus(P/Q,\c).
\]
\end{lem}

\begin{rem}[uniqueness of local
splitting]\label{sec2:rem:uniqueness of local splitting} It
follows from the lemma, applied to the union of two homotopies,
that the simplicial homotopy described in the lemma is unique up
to all higher homotopies, i.e., the space of such homotopies forms
a contractible space.
\end{rem}

This lemma is meant to be used as follows. Given any
$(R,G)$-expansion functor $\xi$ and subfunctor $\eta$ so that
either $\eta$ or $\xi/\eta$ is acyclic, we can deform $\xi$ into a
direct sum
\[
    \xi\simeq\eta\oplus\xi/\eta
\]
in two steps:
\begin{enumerate}
  \item Deform the poset functor (with expansions) so that
  every generator of $\eta$ is less than every generator of
  $\xi/\eta$.
  \item Apply the lemma one simplex at a time on the domain
  complex to kill the incidence $f(\s)$ between $\xi/\eta$ and
  $\eta$.
\end{enumerate}

A jazzed-up version of this simple argument gives the following.

\begin{thm}[Splitting Lemma]\label{sec2:thm:splitting lemma}
Suppose that $R$ is a ring $G$ is a subgroup of $R^\times$ and $B$
is a connected simplicial set with base point $v\in B_0$. Let
\[
    \xi:B\to \Wh\bu(R,G)
\]
be a family of based free chain complexes having a family of
subcomplexes given by a functor $\eta$. Suppose that either
\begin{enumerate}
  \item the homology of $\xi(v)$, $\eta(v)$ and $\xi(v)/\eta(v)$ are
projective $R$-modules on which $\pi_1B$ acts trivially or
  \item $R$ is a field and the homology of $\xi(v)$, $\eta(v)$ and $\xi(v)/\eta(v)$
  are upper triangular as $\pi_1B$-modules.
\end{enumerate}
Then the canonical cone of $\xi$ splits up to homotopy:
\[
    C(\xi)\simeq C(\eta)\oplus C(\xi/\eta).
\]
\end{thm}

In the special case when $R=M_n(\CC)$ and $G=U(n)$ this, together
with additivity, implies that the higher torsion invariant is a
sum:
\[
    \t_k(C(\xi))=\t_k(C(\eta))+\t_k(C(\xi/\eta)).
\]

\subsection{Applications of the Splitting
Lemma}\label{sec2:subsec:appl of splitting lemma}

In order to apply the Splitting Lemma we need to construct a
subfunctor of the expansion functor
\[
    \xi:\simp B \to\Wh\bu(R,G)
\]
given by a fiberwise oriented GMF
\[
    f:E\to\RR.
\]
The simplest case is if there is a regular value $c_t$ of
$f_t:M_t\to\RR$ which varies smoothly with $t\in B$. Then the
critical points of $f_t$ with critical value $<c_t$ generate a
family of subcomplexes $D_t$ of $C_\ast(f_t;\F)$ giving a
subfunctor
\[
    \eta:\simp B\to\Wh\bu(R,G)
\]
of $\xi$. Geometrically, $E$ decomposes as a union of two bundles
$E=E_1\cup E_2$ along $\d_1E_1=\d_0E_2$ where
\begin{align*}
    E_1&=\{x\in E\st f_t(x)\leq c_t (t=p(x)\}\\
    E_2&=\{x\in E\st f_t(x)\geq c_t\}\\
    \d_1E_1=d_0E_2&=\{x\in E\st f_t(x)= c_t\}
\end{align*}
We refer to this situation by saying that $E_2$ is \emph{stacked
on top of} $E_1$.

\begin{cor}[stacking lemma]\label{sec2:cor:torsion of
stacks} Suppose that $E_1,E_2$ are smooth bundles over $B$ so that
$\d_1E_1=\d_0E_2$. Take the union of $E_1$ and $E_2$ along
$\d_1E_1=\d_0E_2$. Let $\F_1,\F_2$ be the restrictions to
$E_1,E_2$ of a Hermitian coefficient system $\F$ on $E_1\cup E_2$.
Then
\[
    \t_k(E_1\cup E_2,\d_0E_1;\F)=\t_k(E_1,\d_0E_1;\F_1)+\t_k(E_2,\d_1E_2;\F_2)
\]
assuming the terms are defined.
\end{cor}

\begin{proof} Suppose first that there is a fiberwise framed
function $f:E_1\cup E_2\to [0,2]$ so that $E_i=f^{-1}[i-1,i]$.
Then the Splitting Lemma says that the canonical cone of
$\xi=C_\ast(f;\F)$ splits:
\[
    C(\xi)\simeq C(\eta)\oplus C(\xi/\eta).
\]
But $\eta\cong C_\ast(f|E_1;\F_1)$ since $E_1^Q\simeq E^Q$ for
every closed subset $Q$ of $P_1^+(\s)$ and $\xi/\eta\cong
C_\ast(f|E_2;\F_2)$ by excision of $E_1$. ($E_1|\s$ is a closed
subset of every $E^Q|\s$ when $Q$ contains $P_1^+(\s)$.)

If there is no fiberwise framed function we can first take the
product with $S^{2N}$ and divided all torsions by $2$.
\end{proof}

As a consequence of this corollary we get a \emph{gluing formula}.
Suppose that $E_1, E_2$ are smooth bundles over $B$ so that part
of the boundary of $E_1$ agrees with part of the boundary of
$E_2$:
\[
    \d_0E_1=\d_0E_2.
\]
Let $E_1\cup E_2$ denote the union of $E_1$, $E_2$ along
$\d_0E_1=\d_0E_2$.

\begin{cor}[gluing formula]\label{sec2:cor:gluing formula}
\[
    \t_k(E_1\cup
    E_2;\F)=\t_k(E_1;\F_1)+\t_k(E_2;\F_2)-\t_k(\d_0E_1;\F_0)
\]
where $\F$ is a Hermitian coefficient system on $E_1\cup E_2$
which restricts to $\F_1,\F_2,\F_0$ and where we assume that the
terms are all defined.
\end{cor}

\begin{proof}
Choose a collar neighborhood $D$ of $\d_0E_1$ in $E_1\cup E_2$.
Then $D$ is a disk bundle over $E_1$ so, by the positive
suspension lemma (\ref{sec1:lem:positive suspension lemma},
\ref{sec1:thm:def of cellular chain complex for d0E empty}),
\[
    \t_k(D;\F|D)=\t_k(\d_0E_1;\F_0),
\]
assuming this torsion is defined.

The union $E_1\cup E_2$ is given by stacking $E_1\coprod E_2$ on
top of $D$. Therefore, by the stacking lemma, we have
\begin{align*}
    \t_k(E_1\cup E_2;\F)&=\t_k(D;\F|D)+\t_k(E_1\smalldisj
    E_2;\F)\\
    &=\t_k(\d_0E_1;\F_0)+\t_k(E_1;\F_1)+\t_k(E_2;\F_2),
\end{align*}
where we use the obvious torsion of disjoint union formula
(\ref{sec2:cor:torsion of disjoint union}).
\end{proof}

Another example of the same idea is the decomposition of the
higher torsion of a manifold bundle into the torsion of its
boundary and the relative torsion.

\begin{cor}[relative torsion]\label{sec2:cor:boundary torsion and relative torsion}
Suppose that $E\to B$ is a smooth manifold bundle whose fiber $M$
has boundary $\d M=\d_0 M\coprod\d_1 M$ where each $\d_i M$ is the
fiber of a subbundle $\d_iE$ of $E$. Then
\[
    \t_{k}(E;\F)=\t_{k}(E,\d_0 E;\F)+\t_{k}(\d_0 E;\F)
\]
when the terms are defined.
\end{cor}

\begin{proof}
This follows from the gluing formula since $E$ is equivalent to
the union of $\d_0E\times I$ and $E$ along their common boundary
$d_0E\times 0$.
\end{proof}

One final example of the same idea shows why the Suspension
Theorem is ``geometrically obvious.''

\begin{cor}[easy geometric suspension]\label{sec2:cor:easy
geometric suspension} Suppose that $E\to B$ is a bundle with
closed manifold fibers and $\F$ is a Hermitian coefficient system
on $E$. Then
\[
    \t_k(E\times D^m,E\times S^{m-1};\pi^\ast\F)=(-1)^m\t_k(E;\F)
\]
assuming the terms are defined.
\end{cor}

\begin{proof} The theorem is a tautology for $m=0$. So suppose by
induction that the theorem holds for $m\geq0$. Then, writing
$E\times S^m$ as a union of $E\times D_-^m$ and $E\times D_+^m$
along $E\times S^{m-1}$, the stacking lemma gives:
\begin{align*}
    \t_k(E\times S^m)&=\t_k(E\times D_-^m)+\t_k(E\times
    D_+^m,E\times S^{m-1})\\
    &=\t_k(E)+(-1)^m\t_k(E)
\end{align*}
by induction and Theorem \ref{sec1:thm:def of cellular chain
complex for d0E empty}. So, the relative torsion formula above
gives
\begin{align*}
    \t_k(E\times D^{m+1},E\times S^m)&=\t_k(E\times
    D^{m+1})-\t_k(E\times S^m)\\
    &=[1-(1+(-1)^m)]\t_k(E)=(-1)^{m+1}\t_k(E).
\end{align*}
\end{proof}

Suppose that $\Sig(f)$ is a disjoint union of topologically closed
subsets
\[
    \Sig(f)=X_1\smalldisj X_2\smalldisj\cdots\smalldisj X_m
\]
so that $X[i]=X_1\coprod\cdots\coprod X_i$ is poset closed in
$\Sig(f)$ over every simplex in $B$, i.e., so that there is no
trajectory of $\grad f$ going from any point in $X_i$ to any point
in $X_j$ for $j>i$ . Then there are corresponding subfunctors
$\xi[i]$ of $\xi=C_\ast(f;\F)$. If the subquotients
$\xi_i=\xi[i]/\xi[i-1]$ are acyclic then, by the Splitting Lemma
the cellular chain complex functor of $f$ is homotopic to a
disjoint union of these subquotients:
\[
    C_\ast(f;\F)=\xi\simeq\xi_1\oplus\cdots\oplus\xi_m.
\]

We need to know that the algebraic subquotients, $\xi_i$, are
determined by the geometry of the trajectories of $\grad f$ coming
out of the critical points in $X_i$. We call this the local
equivalence lemma.

\subsection{Local equivalence lemma}\label{sec2:subsec:local
equivalence}

We state the lemma first and give the definitions afterwards. The
\emph{local equivalence lemma} assumes that we have two fiberwise
oriented GMF's $f:E\to\RR$, $f':E'\to\RR$ which are
\emph{geometrically equivalent} on ``convex'' subsets $S\subseteq
E$ and $S'\subseteq E'$ in the sense that there is a
diffeomorphism $\c:S\to S'$ over $B$ which is an isometry on each
fiber covered by an isomorphism of coefficient sheaves
$\widetilde{\c}:\F\to\F'$ so that $f'\c=f|S+c_t$ where $c_t$ is a
function of $t=p(x)\in B$. In other words, the following diagram
commutes. {\dgVERTPAD=1em\relax \dgHORIZPAD=2em\relax
\dgARROWLENGTH=2em\relax
\[
\begin{diagram}
    \node{\F}\arrow{e}\arrow{s,l}{\widetilde{\c}}
    \node{S}\arrow{e,t}{\subseteq}\arrow{s,r}{{\c}}
    \node{E}\arrow{se,t}{(p,f+cp)}\\
    \node{\F'}\arrow{e}\node{S'}\arrow{e,t}{\subseteq}
    \node{E'}\arrow{e,b}{(p',f')}\node{B\times\RR}
\end{diagram}
\]}

\noi Note that, since $\c$ is a fiberwise isometry, its derivative
takes $\grad f$ to $\grad f'$.

\begin{lem}[local equivalence lemma]\label{sec2:lem:local
equivalence} Under these conditions, the local subquotient of
$C_\ast(f;\F)$ corresponding to $S$ is equivalent to the portion
of $C_\ast(f';\F')$ corresponding to $S'$. In particular, the
cellular chain complex functor of $S$ is well defined without
reference to the ambient bundle $E$.
\end{lem}

The idea and basic properties of convex sets comes from
discussions with Gabriel Katz.

\begin{defn}[convex set]\label{sec2:def:convex sets} Let $T$ be a compact
subset of $B$. Then a compact subset $S$ of $E|T$ will be called
\emph{convex} over $T$ (with respect to $f:E\to\RR$ and the
vertical metric) if there are functions $a,b:T\to\RR$, written
$t\mapsto a_t,b_t$ with $a_t<b_t$ for all $t\in T$, satisfying the
following.
\begin{enumerate}
  \item The fiber $S_t$ of $S$ over $t\in T$ is a codimension $0$
submanifold of $M_t$ (with corners) which varies smoothly with
$t\in T$. (See Figure~\ref{sec2:fig:octagon}.)
\begin{figure}
\includegraphics{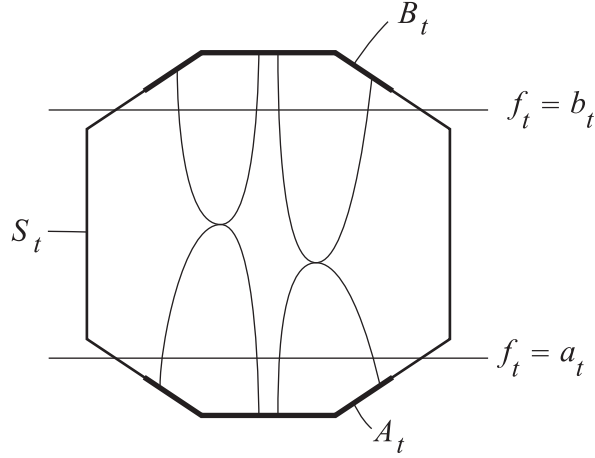}
\caption{$S_t$ (octagon) is called \emph{convex} if
$A_t,B_t,a_t,b_t$ exist.}\label{sec2:fig:octagon}
\end{figure}
  \item All critical points $y$ of $f_t$ in $S_t$ have critical value
  \[
    a_t<f_t(y)<b_t.
  \]
  \item The vertical gradient of $f$ is nowhere zero on the vertical boundary of
  $S$, i.e., $\grad f_t\noteq0$ along $\d(S_t)$ for each $t\in T$.
  \item There are compact subsets $A_t,B_t$ of
$\d(S_t)$ which vary smoothly with $t\in T$ satisfying the
following.
\begin{enumerate}
  \item $f_t(x)\leq a_t$ for all $x\in A_t$.
  \item $f_t(x)\geq b_t$ for all $x\in B_t$.
  \item Every point in $\d(S_t)$ which is incident under a
  critical point of $f_t$ in $S_t$ lies in
  $A_t$ and $\grad f_t$ points inward at that point (so that the
  smooth flow of $\grad f_t$ carries $A_t$ into the interior of
  $S_t$).
  \item Every point in $\d(S_t)$ which is incident over a
  critical point of $f_t$ in $S_t$ lies in
  $B_t$ and $\grad f_t$ points outward at that point (so that the
  smooth flow of $-\grad f_t$ carries $B_t$ into the interior of
  $S_t$).
\end{enumerate}
\end{enumerate}
\end{defn}

One easy example of a convex set is given by
\[
    S_t=f_t^{-1}[a_t,b_t]
\]
where $a_t<b_t$ are regular values of $f_t$ which vary smoothly
with $t$. Suppose also that $\d_0M\subseteq f^{-1}(-\infty,a_t)$
and $\d_1M_t\subseteq f^{-1}(b_t,\infty)$. Then we can take
$A_t=f_t^{-1}(a_t)$ and $B_t=f_t^{-1}(b_t)$.

Another example is given, in two steps, as follows. First, let
$S_t$ be a cobordism from $A_t$ to $B_t$ where these are
codimension zero submanifolds of $f_t^{-1}(a_t), f_t^{-1}(b_t)$
respectively where $a_t<b_t$ are regular values of $f_t$ and
\[
    \d S_t=A_t\cup B_t\cup\d A_t\times I.
\]
Here, the side set $\d A_t\times I$ should be a union of
trajectories of $\grad f_t$ going from $\d A_t$ up to $\d B_t\cong
\d A_t$. There is a problem with this example. Along the corner
set $\d A_t\times\d I$, the gradient of $f_t$ is tangent to part
of the boundary. So, $S_t$ is not convex with this choice of
$A_t,B_t$. But, if we delete a small collar neighborhood of the
boundary of both $A_t$ and $B_t$, the same set $S_t$ will be
convex.

The main nontrivial example is the following.

\begin{eg}[main example]\label{sec2:eg:main example} We start with
a smooth bundle $E\to B$ with fiber $Y$. Suppose that $g:E\to \RR$
is a fiberwise oriented GMF with the property that $\d_0E\subset
g^{-1}(-\infty,0)$, $\d_1E\subseteq g^{-1}(3\delta,\infty)$ and,
for each $t\in B$, the critical values of $g_t$ in the closed
interval $[0,3\delta]$ lie in the open interval $(\delta,2\delta)$
for some $\delta>0$. Let
\[
    E[0,3\delta]=g^{-1}[0,3\delta]
\]
with fibers $Y_t[0,3\delta]=g_t^{-1}[0,3\delta]$ (and, more
generally, $Y_t(J)=g_t^{-1}(J)$ for any $J\subseteq\RR$).

Let $\g$, $\g'$ be vector bundles over $B$ of dimension $m$, $m'$
with structure groups $O(m)$, $O(m')$ and let $\e D$, $\e D'$ be
the $\e$-disk bundles of $\g$, $\g'$ where $\e$ is a positive real
number so that $\e^2\geq2\delta$.

Let $S$ be the fiber product of $E[0,3\delta]$, $\e D$ and $\e D'$
over $B$ with the fiberwise product metric. Thus each fiber is a
product
\[
    S_t=Y_t[0,3\delta]\times \e D_t\times \e D'_t
\]
with the product metric. Let $f:S\to\RR$ be given on each fiber
$S_t$ by
\[
    f_t(y,x,z)= g_t(y)-\length{x}^2+\length{z}^2.
\]
Then
\[
    \grad f_t(y,x,z)=(\grad g_t(y),-2x,2z).
\]
Consequently, the critical points of $f_t$ are the same as those
of $g_t|Y_t[0,3\delta]$.
\[
    \Sig(f_t)=\Sig(g_t|Y_t[0,3\delta])\times0\times0.
\]
Then we claim that $S$ is a convex set with respect to $f$ where
\begin{enumerate}
  \item $T=B$.
  \item $a_t=0$.
  \item $b_t=3\delta$.
  \item $A_t=Y_t(0)\times \e D_t\times0\ \cup\ Y_t[0,2\delta]\times\d(\e
  D_t)\times0$.
  \item $B_t=Y_t(3\delta)\times0\times \e D'_t\ \cup\
  Y_t[\delta,3\delta]\times0\times\d(\e D_t')$.
\end{enumerate}
The definition is satisfied since all trajectories of $\grad f_t$
emanating from critical points of $f_t$ lie in
$Y_t[\delta,3\delta]\times 0\times\e D_t'$ and the intersection of
this set with $\d S_t$ is exactly $B_t$. Similarly, trajectories
of $\grad f_t$ going towards $\Sig(f_t)$ lie in
$Y_t[0,2\delta]\times\e D_t\times0$ which meets $\d S_t$ along
$A_t$. The inequality $2\delta\leq \e^2$ insures that $f_t\leq0$
on $\ Y_t[0,2\delta]\times\d(\e D_t)\times0\subseteq A_t$ and
$f_t\geq3\delta$ on $Y_t[\delta,3\delta] \times0\times\d(\e
D_t')\subseteq B_t$.
\end{eg}

Sometimes we need to deform the function $f_t$ in the main
example. This happens, for example, when we need $f_t$ to be
fiberwise framed.

\begin{eg}[deformation of main example]\label{sec2:eg:1st example
deformed} Take $S_t=Y_t[0,3\delta]\times \e D_t\times \e D'_t$ and
$f_t(y,x,z)= g_t(y)-\length{x}^2+\length{z}^2$ as in the main
example. Note that
\[
    \Sig(f_t)\subseteq Y_t(\delta,2\delta)\times0\times0.
\]
Suppose that $\widetilde{f}_t$ is a deformation of $f_t$ with
support in a small neighborhood of $\Sig(f_t)$, say
\begin{equation}\label{sec2:eq:neighborhood of sig(f) in main eg}
    Y_t[\delta,2\delta]\times\b D_t\times\b D_t'
\end{equation}
for some small $\b>0$. If
\[
    2\b^2<\delta
\]
then we claim that $S$ will be a convex set for the new function
$\widetilde{f}$ where
\begin{enumerate}
  \item $T=B$.
  \item $a_t={\delta}/{2}$.
  \item $b_t=5\delta/2$.
  \item $A_t=Y_t(0)\times \e D_t\times\b D_t'\ \cup\ Y_t[0,2\delta]\times\d(\e
  D_t)\times\b D_t'$.
  \item $B_t=Y_t(3\delta)\times\b D_t\times \e D'_t\ \cup\
  Y_t[\delta,3\delta]\times\b D_t\times\d(\e D_t')$.
\end{enumerate}
Any trajectory of $\grad \widetilde{f}_t$, before reaching $\d
S_t$ must pass through the boundary of the set
(\ref{sec2:eq:neighborhood of sig(f) in main eg}). At that point
it becomes equal to $\grad f_t$ and thus goes up to $B_t$ as does
any point in the set $Y_t[\delta,3\delta]\times\b D_t\times\e
D_t'$. Since the $x$-coordinate is now allowed to be in $\b D_t$
we need to subtract $\b^2$ from $b_t$.
\end{eg}

Another example that we will use is the following.

\begin{eg}[second example]\label{sec2:eg:second example}
Let $g:E\to \RR$, $E[0,3\delta]$ with fibers $Y_t[0,3\delta]$ be
as in the main example. Also, let $\g$, $\g'$ be as before and let
$\e D$, $\e D'$ be the $\e$-disk bundles of $\g$, $\g'$ where $\e$
is a positive real number so that
\[
    3\delta\leq\e^3.
\]

As before, let $S$ be the fiber product of $E[0,3\delta]$, $\e D$
and $\e D'$ over $B$ with the fiberwise product metric. Take the
product with $J^2$ where $J=[-\e,\e]$. Then $W=S\times J^2$ is a
bundle over $B\times J$ whose fiber over $(t,u)\in B\times J$ is
\[
    W_{t,u}=S_t\times J=Y_t[0,3\delta]\times \e D_t\times \e
    D'_t\times J.
\]
We take the fiberwise product metric on $W$. Let $f:W\to\RR$ be
given on each fiber $W_{t,u}$ by
\[
    f_{t,u}(y,x,z,w)= g_t(y)-\length{x}^2+\length{z}^2+w^3-uw.
\]
Then
\[
    \grad f_{t,u}(y,x,z,w)=(\grad g_t(y),-2x,2z,3w^2-u).
\]
So,
\[
    \Sig(f_{t,u})=\Sig(g_t|Y_t[0,3\delta])\times0\times0\times\pm\sqrt{u/3}.
\]

The function $f$ has a problem. It may have a singularity of the
form $y^3+w^3$. These bad singularities occur when $x,z,w,u$ are
all zero and $x$ is a b-d singularity of $g_t$. To fix this we
choose $\b>0$ so that $\b<2\e/3$ and $2\b^2<\delta$ and note that
the bad set lies in the interior of the set
\begin{equation}\label{sec2:eq:bad set for example 2}
    Y[\delta,2\delta]\times\b D\times\b D'\times [-\b,\b]\times
    B\times[-\tfrac{2\e}{3},\tfrac{2\e}{3}].
\end{equation}
Suppose that $h$ is a deformation of $f$ with support in the
interior of this set so that $h$ is a fiberwise oriented GMF. (For
example, if $g$ is fiberwise framed then so is $f$ outside of the
set (\ref{sec2:eq:bad set for example 2}) so we can use the
$C^1$-local framed function theorem (\ref{sec1:rem:C1 local FFT})
to find $h$.)

We claim that $W$ is a convex set with respect to the new function
$h$ where
\begin{enumerate}
  \item $T=B\times J$.
  \item $a_{t,u}=\delta/2$.
  \item $b_{t,u}=5\delta/2$.
  \item $A_{t,u}
  =\d S_t\cap\left(Y_t[0,2\delta]\times \e D_t \times\b D_t'\times[-\e,\frac{2\e}{3}]\right)$.
  \item $B_{t,u}
  =\d S_t\cap\left(Y_t[\delta,3\delta]\times\b D_t\times \e D_t' \times
  [-\frac{2\e}{3},\e]\right)$.
\end{enumerate}
The definition is satisfied since all trajectories of $\grad
h_{t,u}$ emanating from critical points of $h_{t,u}$ lie in
$Y_t[\delta,3\delta]\times \b D_t\times\e D_t'\times
[-\frac{2\e}{3},\e]$ and similarly, those going towards
$\Sig(h_{t,u})$ lie in $Y_t[0,2\delta]\times\e D_t\times\b
D_t'\times [-\e,\frac{2\e}{3}]$.
\end{eg}

Before proving the local equivalence lemma we make some basic
observations about convex sets. A convex set $S$ will contain all
\emph{multiple trajectories} of $\grad f_t$ (i.e., a sequence of
trajectories so that the target of each is the source of the next)
passing between two critical points in $S_t$. The reason is that
if a trajectory of $\grad f_t$ coming from a critical point of
$f_t$ in $S_t$ passes out of $S_t$ then it goes through a point in
$\d S_t$ where $\grad f_t$ points outward (away from $S_t$) and
$f_t\geq b_t$. It has no chance to arrive later at a point in $\d
S_t$ where $\grad f_t$ points inward and $f_t\leq a_t<b_t$.

Suppose that $T=\s$ is a small simplex in $B$. Then this property
of convex sets implies that the subset $Q$ of the poset $P(\s)$
given by the critical points of $f|\s$ in $S|\s$ is \emph{half
closed} (the intersection of an open and closed subset) provided
that $\s$ is sufficiently small.

Combinatorially, a subset $Q$ of a poset $P$ is half closed if and
only if $x<y<z$ in $P$ with $x,z\in Q$ implies that $y\in Q$.
Since we are using the minimal partial ordering, the critical
points of $f_t$ in $S_t$ for $t$ fixed satisfies this property as
we explained above. However, for critical points in
$Q\times\s\subseteq\Sig(f)|\s$ it could happen that
\[
    (x,t_1)<(y,t_1)\quad\text{and}\quad(y,t_2)<(z,t_2)
\]
for two different points $t_1,t_2\in\s$. But this would imply that
$f_{t_1}(y)>b_t$ and $f_{t_2}(y)<a_t$ which will not happen if
$\s$ is chosen sufficiently small.

\begin{defn}[local subquotient]\label{sec2:def:local subquotient functor}
A \emph{local subquotient} of an expansion functor
\[
    \xi:\simp X\to\Wh\bu(R,G)
\]
is defined to be the expansion functor generated by a \emph{family
of half closed subsets}. These are given by $Q(x)\subseteq P(x)$
for each $(x,[k])$ in $\simp X$ so that the structure morphisms
\[
    (P(y),\f(y))\to (P(x),\f(x))|y
\]
send $Q(y)$ onto $Q(x)$ and so that, if one member of a collapsing
pair lies in $Q(y)$, so does the other. (Recall that $P(y)$ is a
disjoint union of collapsing pairs $x_i^+,x_i^-$ and another set
$A$ which maps bijectively onto $P(x)$.)
\end{defn}

The purpose of local subquotients is to decompose the expansion
functor $\xi$ as illustrated by the following proposition.

\begin{prop}[acyclic local subquotients]\label{sec2:prop:acyclic
local subquotients} Suppose that $\xi$ is an $(R,G)$-expansion
functor as above and $Q_1,\cdots Q_m$ are families of half closed
subsets of $P$ so that
\[
    P(x)=\coprod Q_i(x)
\]
for every $(x,[k])\in\simp X$. Suppose further that the local
subquotient $\eta_i$ generated by $Q_i$ is acyclic for every $i$.
Then $\xi$ is homotopic to a direct sum
\[
    \xi\simeq\bigoplus\eta_i.
\]
\end{prop}

\begin{proof}
Use the local splitting lemma (\ref{sec2:lem:local splitting
lemma}). Over each simplex $(x,[k])$ in $X$ the $\Delta^k$-family
of chain complexes $\xi(x)$ has a filtration with acyclic
subquotients $\eta_i(x)$. Therefore, by induction on $k$, we can
eliminate (make $0$) all cross terms in the twisted cochain
$\f(x)$. This defines a simplicial homotopy $\f\simeq
\bigoplus\f_i$ where $\f_i$ is the twisted cochain of $Q_i$. But
we also have a monomial equivalence
\[
    \bigoplus\eta_i(x)=(\vee Q_i(x),\bigoplus \f_i(x))\to (P(x),\bigoplus \f_i(x)).
\]
\end{proof}

\begin{proof}[Proof of local equivalence lemma]
By making $a_t$ slightly larger and $b_t$ slightly smaller we may
assume that $f_t<a_t$ on $A_t$ and $f_t>b_t$ on $B_t$. Then the
smooth flow generated by $\grad f_t$ translates $A_t$ and $B_t$
into compact subsets $A_t'$ and $B_t'$ of the level surfaces
$X_t=f_t^{-1}(a_t)$ and $Y_t=f_t^{-1}(b_t)$. Since $\grad f_t$
points into the interior of $S_t$ along $A_t$, the set $A_t'$ will
lie in the interior of $S_t$ as will the set $B_t'$.

Let $V$ be the union over all $t\in B$ of the set $V_t$ of all
points $x\in S_t$ so that
\[
    a_t\leq f_t(x)\leq b_t.
\]
Let $Z=\cup_{t\in B}Z_t$ where $Z_t$ is the set of points in $V_t$
which are incident over $\Sig(f_t|S_t)$. Then $V,Z$ are compact
and\[Z_t\cap Y_t\subseteq B_t'\subseteq Y_t\cap\interior S_t.
\]Note that $V,Z$ are intrinsic to $S$, i.e., independent of the
embedding $S\subseteq E$. So we can use them to construct the
intrinsic subquotient functor given by $S$.

For every sufficiently small transverse simplex $\s$ in $B$ let
$Q^+(\s)$ be the half-closed subset of $P^+(\s)$ corresponding to
the singularities of $f|\s$ in $S|\s$. ($Q^+(\s)$ is the set of
components of $\Sig(f|\s)\cap S$.) For every closed subset $L$ of
$Q^+(\s)$ let $V^{L}|\s$ be the subset of $V|\s$ consisting of all
points which are not incident over $Q^+(\s)/L$. In particular,
$V^\emptyset(\s)=(S-Z)|\s$. Take the total singular complex of the
pair $(V|\s,V^\emptyset|\s)$ together with the filtration given by
the subsets $V^{L}|\t$ for all faces $\t$ of $\s$. We claim that
this gives a filtered chain complex which is filtered
quasi-isomorphic to the total complex of the subquotient functor.

We need to compare this intrinsically defined filtered complex
with the actual subquotient functor. Let $K^+(\s)$ be the closure
in $P^+(\s)$ of the set of all components $x$ containing elements
with $f_t(x)< a_t$. We also add to $K^+(\s)$ members of collapsing
and expanded pairs whose mate is already in $K^+(\s)$. Then
$K^+(\s)$ will be disjoint from $Q^+(\s)$ and $K^+(\s)\coprod
Q^+(\s)$ will be closed. We also know, by
Proposition~\ref{sec1:prop:minimal ordering}, that $f_t(x)< b_t$
for every $t\in \s$ and every $x\in\Sig(f_t)$ in any component in
$K^+(\s)$. By definition of a convex set we also know that the set
of all points in $S|\s$ which are incident over a point in the
complement of $Q^+(\s)\cup K^+(\s)$ is a closed subset disjoint
from $Z_t$.

Given any closed subset $L$ of $Q^+(\s)$ and any face $\t$ of
$\s$, the inclusion maps
\[
    (V^L|\t,V^\emptyset|\t)\to(E^{L\cup
    K^+(\s)}\cap W|\t,E^{K^+(\s)}\cap W|\t)\xrarrow{\simeq}(E^{L\cup
    K^+(\s)}|\t,E^{K^+(\s)}|\t)
\]
induce isomorphisms in homology. Here $W$ is given by
\[
    W_t=f_t^{-1}(-\infty,b_t]
\]
and the deformation retraction onto the intersection with $W$ is
given by the smooth flow of $-\grad f$. The first mapping is
excisive.

By Theorem~\ref{sec1:thm:existence of twisted cochains}, the
twisted cochain generated by $Q^+(\s)$ is intrinsically defined.
So the subquotient functor is isomorphic to the one for $S'$.
\end{proof}

If we take the construction of the intrinsically defined
subquotient functor for $S$ applied to the main example
(\ref{sec2:eg:main example}), we get an $m$-fold suspension of the
expansion functor of $f$ by the Thom isomorphism theorem (as in
the case of the negative suspension lemma \ref{sec2:lem:negative
suspension lemma}).

\begin{prop}[subquotient functor of main example]
\label{sec2:prop:local suspension lemma} The intrinsic cellular
chain complex functor for the convex set in
Example~\ref{sec2:eg:main example} is the $m$-fold suspension of
the expansion functor for $f:E\to I$:
\[
    C_\ast(h;\r)=\Sig^m C_\ast(f;\r).
\]
\end{prop}

Now we discuss the second example (\ref{sec2:eg:second example}).
We will not obtain a complete description of the subquotient
functor because there is a problem with the b-d points. The
problem is that
\[
    h(x,y)=x^3+y^3
\]
is not an allowable singularity. Suspensions of these
singularities occur in the second example at the points where both
$f$ and $g$ have b-d points.


\subsection{Product formula}\label{sec2:subsec:product formula}

We are now ready to prove one of the main formulas.

\begin{lem}[product formula]\label{sec2:lem:product formula}
Suppose $Y$ is a compact smooth manifold so that $\d_0E\times\d
Y=\emptyset$, i.e., either $Y$ is closed or $\d_0E$ is empty.
Then,
\[
    \t_k(E\times Y,\d_0E\times Y;\F')=\chi(Y)\t_k(E,\d_0E;\F)
\]
assuming that $\t_k(E,\d_0E;\F)$ is defined. Here $\F'$ is the
pull-back of $\F$ along the projection $E\times Y\to E$.
\end{lem}

\begin{rem}[definition of higher torsion]\label{sec2:rem:the product formula is a key lemma}
This is the key lemma that we are using in this paper to
circumvent one of the difficulties in the definition of higher FR
torsion in the relative case. The proof is the same as in
\cite{[I:BookOne]} except that we must be careful to avoid
circular reasoning.
\end{rem}

\begin{proof}The assumption $\d_0E\times\d Y=\emptyset$ implies
that the boundary of $E\times Y$ is a disjoint union
\[
    \d(E\times Y)=\d_0E\times Y\smalldisj (\d_1E\times Y\cup E\times\d
    Y).
\]

Choose a Morse function $g:Y\to[0,n+1]$ so that $g^{-1}(n+1)=\d Y$
with critical points $y_0,y_1,\cdots,y_n$ having distinct integer
critical values $g(y_i)=i$. Next, suppose there exists a fiberwise
framed function
\[
    f:(E,\d_0E)\to(I,0)
\]
so that $f^{-1}(1)=\d_1E$. Consider the smooth function
\[
    h:E\times Y\to[0,n+2]
\]
given by $h(x,y)=f(x)+g(y)$.

For each $t\in B$ we get a smooth function $h_t:M_t\times Y\to\RR$
with critical set $\Sig(h_t)=\Sig(f_t)\times\Sig(g)$. At each
critical point $(x,y_i)$ we can rearrange the local coordinates so
that the negative quadratic terms for $y_i$ go to the beginning.
This gives $h$ the structure of a fiberwise framed function.

Since the gradient of $h_t$ points inward along $\d_0M_t\times Y$
and outward along
\[
    \d_1M_t\times Y\cup M_t\times\d Y
\]
for all $t\in B$ we can modify $h$ in a small neighborhood of $\d
(E\times Y)$ without introducing new critical points so that
\[
    h^{-1}(0)=\d_0E\times Y,\quad\d_1E\times Y\cup E\times\d Y\subseteq h^{-1}[n+1,n+2].
\]

Then, $h^{-1}[0,i+1]$, $i=0,\cdots,n$ gives a filtration of
$E\times Y$ resulting in a filtration
\[
    \xi_0\subseteq\xi_1\subseteq\cdots\subseteq\xi_n=\xi
\]
of the family of chain complexes $\xi=C_\ast(h;\F')$. Since $y_0$
has index $0$, $\xi_0\cong C_\ast(f)$ and, for each $i$, the
quotient $\xi_i/\xi_{i-1}$ is isomorphic to a suspension of
$C_\ast(f;\F)$
\[
    \xi_i/\xi_{i-1}\cong\Sig^{\ind(y_i)}C_\ast(f;\F).
\]
By the Suspension Theorem this has torsion
\[
    \t_k(C(\xi_i/\xi_{i-1}))=(-1)^{\ind(y_i)}\t_k(CC_\ast(f;\F)).
\]
The lemma now follows from the Splitting Lemma
\begin{align*}
    \t_k(E\times Y,\d_0E\times
    Y;\F')&={\t_k(C(\xi))=\sum\t_k(C(\xi_i/\xi_{i-1}))}\\
    &=\sum(-1)^{\ind(y_i)}\t_k(CC_\ast(f;\F))=\chi(Y)\t_k(E,\d_0E;\F).
\end{align*}

Now suppose that $\dim M<\dim B$ and there is no fiberwise framed
function for $(E,\d_0E)$. To avoid circular reasoning we cannot
use the usual trick of taking the product with an even dimensional
sphere. Instead we take the product with two spheres
$S^{2N},S^{2N'}$ with $N,N'$ very large.

Since $E\times S^{2N}$, $E\times S^{2N'}$ have large dimensional
fibers, the lemma holds for them. Thus:
\begin{align*}
    \t_k(E\times S^{2N},\d_0E\times
    S^{2N};\F')&=\frac12\t_k(E\times S^{2N}\times S^{2N'},\d_0E\times S^{2N}\times
    S^{2N'};\F''')\\
    &=\t_k(E\times S^{2N'},\d_0E\times S^{2N'};\F'').
\end{align*}
Therefore, $\t_k(E\times S^{2N},\d_0E\times S^{2N};\F')$ is
\emph{independent of the choice of $N$} as long as $N$ is
sufficiently large. So, we can \emph{define} $\t_k(E,\d_0E;\F)$ to
be half of this:
\begin{equation}\label{sec2:eq:define t to be half of t(ExS2N)}
  \t_k(E,\d_0E;\F):=\frac12\t_k(E\times S^{2N},\d_0E\times S^{2N};\F').
\end{equation}
The above argument shows that this equation holds whenever
$(E,\d_0E)$ admits a fiberwise framed function.

Now we can complete the proof of the lemma (suppressing the
coefficients and the $\d_0E$ terms from the notation).
\[
    \t_k(E\times Y)=\frac12\t_k(E\times Y\times
    S^{2N})=\frac12\chi(Y)\t_k(E\times S^{2N})=\chi(Y)\t_k(E).
\]
\end{proof}

\subsection{Transfer for coverings}\label{sec2:subsec:transfer}

A Morse function on a manifold $M$ gives a Morse function on any
covering of $M$ and the associated chain complexes are related by
a simple transfer formula. We use this elementary observation to
explain why the Riemann zeta function appears in the formula for
the higher FR torsion invariant.

This argument is a simplification of an argument which appears in
\cite{[I:BookOne]}. I should thank William Hoffman and Neil
Stoltzfus for helping me to clarify this idea.

Suppose that $G$ is a finite group which acts freely on the total
space $E$ so that the action is smooth and commutes with the
projection $p:E\to B$. Then we get an induced bundle
\[
    E_G=E/G\to B
\]
with fiber $M_G=M/G$. Note that $E$ is a finite covering space of
$E_G$. If $H$ is any subgroup of $G$, then $E_H=E/H$ is also a
finite covering space of $E_G$ and we have covering maps
\[
    E\to E_H\to E_G.
\]

Conversely, any finite connected covering $\widetilde{X}\to X$ has
this form where $G$ is the quotient of $\pi_1X$ by the largest
normal subgroup contained in the image of $\pi_1\widetilde{X}$.

\begin{thm}[transfer for coverings of $E$]\label{sec2:thm:transfer formula for groups}
Suppose that $V$ is a unitary representation of $H$ so that
$\pi_1B$ acts trivially on $H_\ast(M_H;V)\cong
H_\ast(M_G;\Ind_H^GV)$. Then
\[
    \t_k(E_G;\Ind_H^GV)=\t_k(E_H;V).
\]
\end{thm}

\begin{proof}
Suppose first that there is a fiberwise framed function
$f_G:E_G\to\RR$. Composing with the covering maps $E\to E_H\to
E_G$ we get induced fiberwise framed functions $f_H:E_H\to \RR$
and $f:E\to \RR$. The cellular chain complex $C_\ast(f)$ is a
family of free $\CC G$-complexes over $B$ and the higher torsion
invariants $\t_k(E_H;V)$ are induced by
\begin{align*}
    C_\ast(f_H;V)&= C_\ast(f)\otimes_{\CC H} V\\
    &\cong C_\ast(f)\otimes_{\CC G}\CC G\otimes_{\CC H} V\\
    &=C_\ast(f)\otimes_{\CC G}\Ind_H^GV=C_\ast(f_G;\Ind_H^GV)
\end{align*}
considered as a mapping from $B$ into the space of acyclic chain
complexes (using the canonical cone construction). Consequently,
they have the same torsion:
\[
    \t_k(E_G;\Ind_H^GV)=\t_k(C_\ast(f_G;\Ind_H^GV))=\t_k(C_\ast(f_H;V))=\t_k(E_H;V).
\]

If the dimension of $M$ is not large enough to admit a fiberwise
framed function we first take the product of all bundles with a
large even dimensional sphere $S^{2N}$. Then there is a fiberwise
framed function on $E_G\times S^{2N}$ and the above argument tells
us that
\[
    \t_k(E_G\times S^{2N};\Ind_H^GV)=\t_k(E_H\times S^{2N};V).
\]
Divide by $2$ to get the theorem.
\end{proof}

Now we specialize to cyclic groups. Let $Z_n=\ZZ/n$ and for $z$
any $n$-th root of unity let
\[
    \r_z:Z_n\to U(1)
\]
be the $1$-dimensional representation of $Z_n$ which sends the
generator of $Z_n$ to $z$. Then, by Frobenius reciprocity, we have
\[
    \Ind_{Z_n}^{Z_{nm}}(\r_z)=\bigoplus_{\r|H=\r_z}\r=\bigoplus_{\z^m=z}\r_{\z}
\]
Consequently we have:

\begin{cor}[transfer for cyclic coverings]\label{sec2:cor:transfer formula for cyclic groups}
\[
    \t_k(E/Z_n;\r_z)=\sum_{\z^m=z}\t_k(E/Z_{nm};\r_\z).
\]
\end{cor}

Now consider the special case where $M=S^1=U(1)$ and $E\to B$ is
the universal $S^1$-bundle
\[
    EU(1)\to BU(1)=\CC P^\infty.
\]
Then $H^\ast(BU(1))=\ZZ[c]$ and the circle bundle $EU(1)/Z_m\to
BU(1)$ is classified by a mapping
\[
    \f_m:BU(1)\to BU(1)
\]
satisfying $\f_m^\ast(c)=mc$ ($\Rightarrow\f_m^\ast(c^k)=m^kc^k$).

More generally, there is a morphism of $S^1$-bundles:

\[
\begin{CD}
    {EU(1)/Z_{nm}} @>>>{EU(1)/Z_n}\\
    @VVV @VVV\\
    {BU(1)} @>>{\f_m}> {BU(1)}
\end{CD}
\]
Consequently, by naturality of the higher FR torsion invariant, we
have
\[
    \t_k(EU(1)/Z_{nm};\r_z)=m^k\t_k(EU(1)/Z_n;\r_z)
\]
for any $n$-th root of unity $z$. Since this holds in the
universal case, it holds for all oriented $S^1$-bundles. Combining
this with the second transfer formula (\ref{sec2:cor:transfer
formula for cyclic groups}) we get:
\[
    \t_k(E/Z_s;\r_z)=m^k\sum_{\z^m=z}\t_k(E/Z_s;\r_\z).
\]
Here $s=nm$. Now we use the elementary fact that the polylogarithm
function
\begin{equation}\label{sec2:eq:polylogarith function}
    L_{k+1}(z)=\Re\left(\frac1{i^k}\sum_{n=1}^\infty\frac{z^n}{n^{k+1}}\right)
\end{equation}
is the only smooth function on $z\in U(1)-\{1\}$ satisfying the
equation
\[
    L_{k+1}(z)=m^k\sum_{\z^m=z}L_{k+1}(\z)
\]
for all integers $m$. The conclusion is that
\begin{equation}\label{sec2:eq:tau k is a k L(k+1)}
  \t_k(EU(1)/Z_n;\r_z)=n^k a_kL_{k+1}(z)c^k
\end{equation}
where $a_k$ is a real number depends only on $k$. See
\cite{[I:BookOne]} for the missing steps in this argument and for
the proof that:
\begin{equation}\label{sec2:eq:value of a k}
  a_k=-\frac{1}{k!}
\end{equation}
for all $k$.

\begin{thm}[torsion of circle bundles]\label{sec2:thm:torsion of
circle bundles} Let $S^1(\ll)$ be the circle bundle associated to
the complex line bundle $\ll$ over $B$. Then
\[
    \t_k(S^1(\ll)/Z_n;\r_z)=-\frac{n^k}{k!}L_{k+1}(z)c_1(\ll)^k.
\]
\end{thm}

Taking the special case $n=1=z$ we get the higher torsion of the
universal oriented $S^1$-bundle:
\begin{equation}\label{sec2:eq:real torsion of universal S1 bundle}
  \t_{2k}(EU(1))=-\frac1{(2k)!}L_{2k+1}(1)c^{2k}=(-1)^{k+1}\z(2k+1)ch_{2k}(\g)
\end{equation}
where $\g$ is the canonical complex line bundle over $BU(1)$.

\begin{cor}[torsion of lens space bundles]\label{sec2:cor:torsion
of lens space bundles} Let $S^{2n-1}(\xi)\to B$ be the sphere
bundle associated to a complex $n$-plane bundle $\xi$. Then
\[
    \t_k(S^{2n-1}(\xi)/Z_m;\r_z)=-m^k L_{k+1}(z)ch_k(\xi)
\]
where $L_{k+1}(z)$ is given in (\ref{sec2:eq:polylogarith
function}) above.
\end{cor}

\begin{proof}
This is Theorem 5.7.12 in \cite{[I:BookOne]}. An outline of the
proof goes as follows.

Using the splitting principle we may assume that $\xi$ is a sum of
line bundles $\xi=\oplus\ll_i$. This implies that the sphere
bundle $S^{2n-1}(\xi)$ is a fiberwise join of the circle bundles
associated to the $\ll_i$ and
\[
    S^{2n-1}(\xi)/Z_m\cong\ast_i S^1(\ll_i)/Z_m.
\]
Consequently, by the Splitting Lemma, we have
\[
    C_\ast(S^{2n-1}(\xi)/Z_m)\cong\bigoplus_i
    C_\ast(S^1(\ll_i)/Z_m).
\]
So,
\begin{align*}
    \t_k(S^{2n-1}(\xi)/Z_m;\r_z)&=\sum_{i=1}^n\t_k(S^1(\ll_i)/Z_m;\r_z)\\
    &=-\frac{m^k}{k!}L_{k+1}(z)=-m^kL_{k+1}(z)ch_k(\xi).
\end{align*}
\end{proof}

\subsection{More transfer formulas}\label{sec2:subsec:more
transfer}

Another way to do the transfer is on the base. If
\[
    \pi:\widetilde{B}\to B
\]
is a finite covering map then we recall that the \emph{transfer}
or \emph{push-down}
\[
    tr=tr^{\widetilde{B}}_B=\pi_\ast: H^q(\widetilde{B};A)\to
    H^q(B;A)
\]
with any additive group coefficients $A$ is defined at the chain
level by
\begin{equation}\label{sec2:eq:transfer for covering maps}
    \<tr(c_q),\s\>=\sum_{\widetilde{\s}}\<c_q,\widetilde{\s}\>
\end{equation}
for any $q$-cocycle $c_q$ on $\widetilde{B}$ (with coefficients in
$A$) and any singular $q$-simplex $\s:\Delta^q\to B$. The sum is
over all lifting $\widetilde{\s}$ of $\s$ to $\widetilde{B}$.

Suppose that $p:E\to \widetilde{B}$ is a smooth bundle with $\d
E=\d_0E\coprod\d_1E$ and $\F$ is a Hermitian coefficient system on
$E$ so that $H_\ast(M,\d_0M;\F)=0$. Then the higher relative FR
torsion of $(E,\d_0E)$ as a bundle over both $\widetilde{B}$ and
$B$ are defined and related by transfer.

\begin{thm}[transfer for coverings of $B$]\label{sec2:thm:transfer
for coverings of B}
\[
    \t_k(E,\d_0E;\F)_B=tr^{\widetilde{B}}_B(\t_k(E,\d_0E;\F)_{\widetilde{B}}).
\]
\end{thm}

\begin{rem}[torsion commutes with transfer]\label{sec2:rem:torsion commutes with transfer} At the cellular
chain complex level this theorem says that the torsion of the
transfer is the transfer of the torsion, i.e.,
\begin{equation}\label{sec2:eq:transfer commutes with torsion}
    \t_k(tr^{\widetilde{B}}_B(\xi))=tr^{\widetilde{B}}_B(\t_k(\xi))
\end{equation}
where the \emph{transfer}, $tr^{\widetilde{B}}_B(\xi)$, of an
expansion functor $\xi$ on ${\widetilde{B}}$ is the expansion
functor on $B$ given by
\begin{equation}\label{sec2:eq:transfer of xi}
    tr^{\widetilde{B}}_B(\xi)(\s)=\bigoplus_{\widetilde{\s}}\xi(\widetilde{\s}).
\end{equation}
where the direct sum is over all lifting $\widetilde{\s}$ of $\s$
to $\widetilde{B}$.
\end{rem}

\begin{proof}
The proof is the same as that of the additivity theorem. We choose
a fiberwise framed function
\[
    f:(E,\d_0E)\to (I,0),
\]
taking the product of $E$ with $S^{2N}$ if necessary. The cellular
chain complex $C_\ast(f;\F)$ gives a $(M_r(\CC),U(r))$-expansion
functor
\[
    \xi:\simp \widetilde{B}\to\Wh\bu^h(M_r(\CC),U(r)).
\]

For any small simplex $\s$ in $B$, the inverse image of $\s$ in
$E$ is a disjoint union,
\[
    E|\s=\coprod_{\widetilde{\s}}E|\widetilde{\s},
\]
over all liftings $\widetilde{\s}$ of $\s$ to $\widetilde{B}$.
Consequently, using the same function $f$, the resulting expansion
functor on $\widetilde{B}$ is the transfer of $\xi$ as given in
(\ref{sec2:eq:transfer of xi}). Thus it suffices to show that
transfer commutes with torsion (\ref{sec2:eq:transfer commutes
with torsion}).

Using the $2$-index theorem (Theorem~\ref{sec1:thm:2-index thm}),
we can deform $\xi$ into a family of complex matrices $\eta$
parametrized by $\widetilde{B}$. The homotopy $\xi\simeq\eta$ give
a homotopy between the push-down expansion functors
\[
    \overline{\xi}=tr^{\widetilde{B}}_B(\xi):\simp B\to \Wh\bu^h(M_r(\CC),U(r))
\]
\[
    \overline{\eta}=tr^{\widetilde{B}}_B(\eta):\simp B\to
    \Wh\bu^{h[0,1]}(M_r(\CC),U(r))
\]
and it remains only to show that the higher torsion of
$\overline{\eta}$ is the transfer of the higher torsion of $\eta$.

The higher FR-torsion of $\overline{\xi}$ is given by integrating
the Kamber-Tondeur form for $\overline{\eta}$ plus polynomial
correction terms. Both of these obviously satisfy the property
that they are additive with respect to direct sum. (They are both
integrals of traces of products of matrix $1$-forms which preserve
the direct sum structure of the input matrices $\overline{\eta}$.)
In other words,
\[
    \<KT_k(\overline{\eta}),\s\>=\sum_{\widetilde{\s}}\<KT_k(\eta),\widetilde{\s}\>
\]
where $KT_k(\overline{\eta})$ is the $2k$-cocycle on $\simp B$
which gives the higher FR-torsion $\overline{\eta}^\ast(\t_k)$ of
$(E,\d_0E;\F)$ over $B$. Thus,
$\t_k(\overline{\eta})=tr^{\widetilde{B}}_B(\t_k(\eta))$ at the
chain level.
\end{proof}

Now suppose that
\[
    f:E\to\RR
\]
is a \emph{fiberwise framed Morse function}, i.e., a fiberwise
framed function without b-d points. Then the singular set
$\Sig=\Sig(f)$ is a disjoint union
\[
    \Sig(f)=\coprod_{i=0}^n\Sig^i(f)
\]
where each $\Sig^i(f)$ is a finite covering space of $B$.

Suppose we have a smooth bundle
\[
    Y\to D\to E
\]
and $\F$ is a Hermitian coefficient system over $D$ with respect
to which $Y$ is acyclic, i.e., $H_\ast(Y;\F)=0$. We consider $D$
as a bundle over $B$. This bundle has corners which we can ignore
since we will be enclosing the singular set of a fiberwise framed
function in a disjoint union of convex sets disjoint from the
corner set.

The transfer theorem for coverings of $B$
(Theorem~\ref{sec2:thm:transfer for coverings of B}) applies to
the restriction of $D$ to $\Sig^i=\Sig^i(f)$ to give
\begin{equation}\label{sec2:eq:transfer for D res Sig i}
  \t_k(D|\Sig^i;\F)_B=tr^{\Sig^i}_B(\t_k(D|\Sig^i;\F)_{\Sig^i}).
\end{equation}
We claim that the alternating sum of these terms gives the higher
torsion of $D$ as a bundle over $B$. This is a special case of the
main theorem of this paper.

\begin{thm}[transfer without b-d points]\label{sec2:thm:transfer
for twisted Morse bundles} Given a fiberwise framed Morse function
$f:E\to\RR$, a smooth bundle $D\to E$ and Hermitian coefficient
system $\F$ over $D$ as above, we get the following.
\[
    \t_k(D;\F)_B=\sum_{i=0}^n(-1)^i\t_k(D|\Sig^i;\F)_B.
\]
\end{thm}

\begin{rem}[push-down without b-d points]\label{sec2:rem:pushdown
without bd points} By (\ref{sec2:eq:transfer for D res Sig i})
this can be expressed as
\[
    \t_k(D;\F)=p_\ast^{\Sig}(\t_k(D|\Sig;\F)_\Sig)
\]
where the \emph{push-down operator},
\[
    p_\ast^\Sig:H^{2k}(\Sig(f);\RR)\to H^{2k}(B;\RR),
\]
is given by the alternating sum of transfer maps:
\[
    p_\ast^\Sig:=\sum_{i=0}^n(-1)^itr^{\Sig^i}_B.
\]
\end{rem}

\begin{proof} We will enclose the singular set of $f$ in a subset
$T$ of $E$ and find a fiberwise framed function $h$ on $D$ which
is locally equivalent over $T$ to the convex set in the main
example (\ref{sec2:eg:main example}). Then
Proposition~\ref{sec2:prop:local suspension lemma} will give us
the expansion functor of $h$.

Let $T^i\subseteq E$ be a tubular neighborhood of
$\Sig^i=\Sig^i(f)$ so that the fiber of the projection map
\[
    \pi:T^i\to\Sig^i
\]
is a product of two disks, an $i$-disk and a twisted $(n-i)$-disk.
\[
    T^i \cong\e D^i\times\e D^{n-i}(\g_i)
\]
where $\e D^i$ is the disk of radius $\e$ in $\RR^i$:
\[
    \e D^i =\{x\in\RR^i\st \length{x}\leq\e\}
\]
and $\e D^{n-i}(\g_i)$ is the $\e$-disk bundle of an $(n-i)$-plane
bundle $\g_i$ over $\Sig^i$. ($\g_i$ is the positive eigenspace
bundle of $D^2f$.)

By the parametrized Morse lemma we may assume that $f|T^i$ is
given by
\begin{equation}\label{sec2:eq:formula for f on Ti}
  f(x,z)=f(s)-\length{x}^2+\length{z}^2
\end{equation}
where $s=\pi(z)\in\Sig^i\subseteq E$ (assuming $\e>0$ is
sufficiently small).
\begin{figure}
\includegraphics{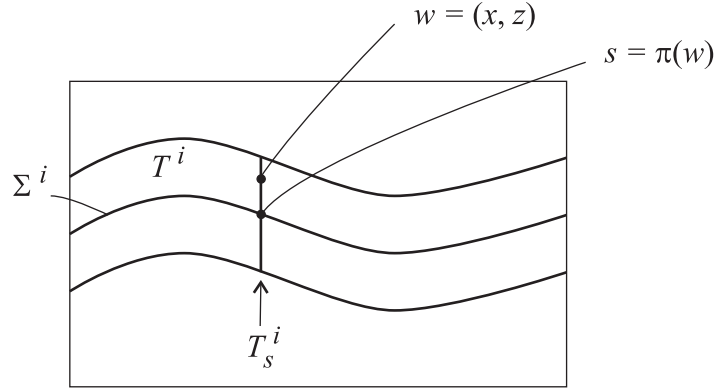}
\caption{$w=(x,z)\in T_s^i=\pi^{-1}(s)$ in $E$
(box).}\label{sec2:fig:morse transfer}
\end{figure}

Since $\pi:T^i\to \Sig^i$ is a bundle with contractible fibers,
the smooth bundle
\[
    Y\to D|T^i\to T^i
\]
is fiber diffeomorphic to $\pi^\ast(D|\Sig^i)$, i.e., we have a
smooth family of diffeomorphisms
\[
    \f_w:Y_w\xrarrow{\approx}Y_s\quad (s=\pi(w))
\]
for all $w\in T^i$ so that $\f_w=\f_s=id$ if $w=s\in\Sig^i$. These
diffeomorphisms give
\begin{equation}\label{sec2:eq:product formula for D over Tsi}
  D|T_s^i\cong \e D^i\times \e D_s^{n-i}\times Y_s
\end{equation}
where $T_s^i=\pi^{-1}(s)$ and $Y_s=D|s$ is the fiber of $D$ over
$s\in\Sig^i$.

Let $g:D|\Sig\to [1,3]$ be a fiberwise framed function, i.e.,
giving a family of framed functions
\[
    g_s:(Y_s,\d Y_s)\to ([1,3],3),\quad s\in\Sig^i
\]
with critical values in $[1,2]$. Then $g$ extends to a fiberwise
framed function on $D|T$ ($T=\coprod T^i$) by
\begin{equation}\label{sec2:eq:extending g to D over T}
  g_w=g_{s}\circ\f_w:Y_w\to Y_{s}\to I.
\end{equation}
Let $\overline{g}:D\to I$ be any smooth extension of $g$
satisfying (\ref{sec2:eq:extending g to D over T}) for all $w\in
T$.

Let $h:D\to\RR$ be given by
\[
    h(y)=f(w)+\delta\overline{g}(y)
\]
if $y\in Y_w$ where $0<\delta<\e^2$ is fixed. If $y\in D|T^i$
then, combining (\ref{sec2:eq:formula for f on Ti}) and
(\ref{sec2:eq:extending g to D over T}) we have, in terms of the
local coordinates (\ref{sec2:eq:product formula for D over Tsi}),
\[
    h(x,z,y)=f(s)-\length{x}^2+\length{z}^2+\delta g_s(y).
\]
The fiberwise gradient of $h$ (over $B$) is equal to zero (on
$D|T^i$) iff $x=0$, $z=0$ and $y\in\Sig(g)$.

Adding the tangent vectors $\grad x_1,\cdots,\grad x_i$ in front
of the framing vectors for $g$, we obtain a framed structure on
$h$ restricted to $D|T$. If $\delta>0$ is chosen sufficiently
small then $h$ will have no other critical points. (The gradient
of $f$ outside of $T$ is bounded below and the gradient of
$\overline{g}$ is bounded above.) Consequently, $h$ is fiberwise
framed with the same critical set as $g$:
\[
    \Sig(h)=\Sig(g).
\]

For each $i$, the fiberwise framed function $g$ on $D|\Sig^i$
gives an expansion functor
\[
    \xi_i:\simp \Sig^i\to \Wh\bu^h(M_r(\CC),U(r)),
\]
where $r$ is the dimension of the fibers of $\F$. Since $D|T$ is a
convex set equivalent to the main example (\ref{sec2:eg:main
example}), the local suspension lemma (\ref{sec2:prop:local
suspension lemma}) tells us that $\Sig^i(\xi_i)$ is a local
subquotient functor of the cellular chain complex functor of $h$.
Since $\Sig^i(\xi_i)$ is acyclic, we can use the local splitting
lemma as in Proposition~\ref{sec2:prop:acyclic local subquotients}
(acyclic local subquotients) to simplicially homotope
$C_\ast(h;\F)$ into a local direct sum:
\begin{equation}\label{sec2:eq:C(h) is the alternating sum of pushdowns of xi}
  C_\ast(h;\F)\simeq \bigoplus_{i=0}^n
  tr^{\Sig^i}_B(\Sig^i(\xi_i)).
\end{equation}

Equation (\ref{sec2:eq:C(h) is the alternating sum of pushdowns of
xi}) proves the theorem since, by additivity, suspension theorem
and the ``transfer commutes with torsion'' formula
(\ref{sec2:eq:transfer commutes with torsion}), we have:
\begin{align*}
\t_k\left( \bigoplus tr^{\Sig^i}_B(\Sig^i(\xi_i))
    \right)&=\sum tr^{\Sig^i}_B(\t_k(\Sig^i(\xi_i)))\\
    &=\sum(-1)^itr^{\Sig^i}_B(\t_k(\xi_i))\\
    &=\sum(-1)^itr^{\Sig^i}_B(\t_k(D|\Sig_i;\F)_{\Sig^i}).
\end{align*}
\end{proof}

%% file: FramPr.tex
\vfill\eject\section{The Framing Principle}

\begin{enumerate}
  \item Statement for Morse bundles
  \item General statement
  \item Push-down/transfer
  \item The Framing Principle
\end{enumerate}

\subsection{Statement for Morse bundles}\label{FP:subsec:statement
for Morse bundles}

We first state the Framing Principle in the case when there is a
fiberwise oriented Morse function $f:E\to\RR$. In that case the
singular set $\Sig(f)$ is a finite covering space of $B$. It is a
disjoint union of covering spaces
\[
    \Sig(f)=\coprod\Sig^i(f).
\]
For each $i\geq0$, we have an $i$-dimensional vector bundle $\g^i$
over $\Sig^i(f)$ given by the negative eigenspace of the second
derivative of $f$ (equal to the tangent plane to the descending
manifold).

If each of these vector bundles $\g^i$ is trivial then $f$ admits
a framed structure so it can be used to define the higher torsion
of the smooth bundle $E$.
\[
    \t_{2k}(E):=\t_{2k}(C_\ast(f))\in H^{4k}(B;\RR)
\]
where $\t_{2k}(C_\ast(f)):=\t_{2k}(C(C_\ast(f)))$ is the torsion
of the canonical cone (Lemma~\ref{sec2:lem:canonical cone}) of the
family of chain complexes defined by $f$. If the bundles $\g^i$
are nontrivial then this formula must be modified.
\begin{equation}\label{FP:eq:Framing Principle for Morse bundles}
    \t_{2k}(E)=\t_{2k}(C_\ast(f))+
    \sum_{i\geq0}(-1)^{i+k}\z(2k+1)\left(p|\Sig^i(f)\right)_\ast\left(
    \tfrac12 ch_{2k}(\g^i\otimes\CC)\right)
\end{equation}
where $ch_{2k}$ is the degree $4k$ component of the \emph{Chern
character} and
\[
    \left(p|\Sig^i(f)\right)_\ast:H^\ast(\Sig^i(f);\RR)\to H^\ast(B;\RR)
\]
is the \emph{push-down} operator given by summing over the sheets
of the covering:
\[
    \<\left(p|\Sig^i(f)\right)_\ast(c),\s\>
    =\sum_{\widetilde{\s}}\<c,\widetilde{\s}\>.
\]
(A cocycle $c$ is evaluated on a singular simplex $\s$ in $B$ by
taking the sum of the values of $c$ on the liftings
$\widetilde{\s}$ of $\s$ to $\Sig^i(f)$.)

Equation (\ref{FP:eq:Framing Principle for Morse bundles}) is
called the \emph{Framing Principle}. Thus equation
(\ref{FP:eq:Framing Principle for Morse bundles}) says that the
Framing Principle holds for fiberwise oriented Morse functions.

The basic example of the use of the Framing Principle is the
following (from \cite{[I:BookOne]}). U. Bunke
\cite{[Bunke:spheres]} had already proven the analogous statement
for higher analytic torsion.

\begin{thm}[torsion of sphere bundles]\label{FP:thm:torsion of oriented sphere bundles}
Let $S^{n-1}\to E(\xi)\to B$ be the sphere bundle of an oriented
$n$-plane bundle $\xi$ over $B$. Then
\[
    \t_{2k}(E(\xi))=(-1)^{k+n-1}
    \z(2k+1)\left(
    \tfrac12 ch_{2k}(\xi\otimes\CC)\right).
\]
\end{thm}

\begin{rem}[relation to Pontrjagin classes]\label{FP:rem:Newton
polynomial}To express this formula in terms of the Pontrjagin
classes $p_i(\xi)$ of $\xi$ (as in the introduction), we need to
substitute
\[
    \tfrac12
    ch_{2k}(\xi\otimes\CC)=\frac1{(2k)!}N_k(p_1(\xi),p_2(\xi),\cdots)
\]
where $N_k$ is the $k$-th \emph{Newton polynomial}, i.e., the
integer polynomial which expresses the power sum
\[
    \sum x_i^k
\]
in terms of the elementary symmetric functions of $x_i$.
\end{rem}

\begin{rem}[logic is reversed]\label{FP:rem:logic is reversed}
This theorem is actually a lemma in the proof of the Framing
Principle. The higher torsion of oriented linear sphere bundles is
computed directly in \cite{[I:BookOne]} and this calculation was
used there to prove a restricted version of the Framing Principle.
Thus, the following argument actually shows that
Theorem~\ref{FP:thm:torsion of oriented sphere bundles} is
equivalent to the Framing Principle in the special case of
fiberwise oriented Morse functions.
\end{rem}

\begin{proof}
Let $D(\xi)$ be the disk bundle of $\xi$. Then $E(\xi)=\d D(\xi)$.
So, Corollary~\ref{sec2:cor:boundary torsion and relative torsion}
gives:
\[
    \t_{2k}(E(\xi))=\t_{2k}(D(\xi))-\t_{2k}(D(\xi),E(\xi)).
\]
To compute $\t_{2k}(D(\xi))$ we take the fiberwise Morse function
$f:D(\xi)\to\RR$ given by
\[
    f(x)=\length{x}^2.
\]
This is fiberwise framed. It has a constant family of chain
complexes ($0\to\ZZ\to 0$ on every fiber). Thus
\[
    \t_{2k}(D(\xi))=0.
\]
To compute the relative torsion $\t_{2k}(D(\xi),E(\xi))$ we take
the fiberwise Morse function
\[
    1-f:(D(\xi),E(\xi))\to(I,0).
\]
The chain complex of the fiber is again constant. However, $1-f$
is not framed since it has a critical point of index $n$ with
descending disk bundle $\g^n=\xi$. Thus the Framing Principle
(\ref{FP:eq:Framing Principle for Morse bundles}) gives:
\[
    \t_{2k}(D(\xi),E(\xi))=(-1)^{n+k}
    \z(2k+1)\left(
    \tfrac12 ch_{2k}(\xi\otimes\CC)\right).
\]
The theorem follows.
\end{proof}

\subsection{General statement}\label{FP:subsec:general statement}

Suppose now that $f:(E,\d_0E)\to(I,0)$ is a fiberwise oriented
GMF. Then the singular set $\Sig(f)$ is a smooth submanifold of
$E$. It is a union of closed subsets
\[
    \Sig(f)=\bigcup\overline{\Sig}^i(f)
\]
which intersect along the birth-death sets
\[
    \Sig^i_1(f)=\overline{\Sig}^i(f)\cap\overline{\Sig}^{i+1}(f).
\]

The negative eigenspace bundle $\g^i$ is defined over all of
$\overline{\Sig}^i(f)$ and is a subbundle of $\g^{i+1}$ over
$\Sig^i_1(f)$. The complementary line bundle is the kernel of
$D^2(f)$. This has an intrinsic orientation given by the intrinsic
third derivative of $f$. By definition of an oriented GMF, the
orientation of $\g^{i+1}$ is given by the orientation of $\g^i$
together with the intrinsic orientation of the complementary line
bundle.

\begin{defn}[stable bundle]\label{FP:def:gamma on Sig(f)}
For any oriented fiberwise GMF $f:(E,\d_0E)\to(I,0)$ we define the
\emph{stable bundle} $\g_f$ be the $n$-dimensional vector bundle
over $\Sig(f)$ (where $n=\dim M$) given by the following clutching
construction. The restriction of $\g_f$ to $\overline{\Sig}^i(f)$
is equal to $\g^i$ plus a trivial $n-i$ plane bundle spanned by
vector fields $v_{i+1},\cdots,v_n$. Along the birth-death set
$\Sig^i_1(f)$, let $v_{i+1}$ be identified with the positive cubic
direction of the kernel of $D^2(f)$ and let $v_{i+2},\cdots,v_n$
be identified so that they extend continuously to
$\overline{\Sig}^i(f)\cup\overline{\Sig}^{i+1}(f)$. Thus each
vector field $v_j$ is continuous over its domain.
\end{defn}

\begin{prop}[stable $+$ unstable]\label{FP:prop:gamma f + gamma -f give nu Sigma}
If $f:(E,\d_0E,\d_1E)\to(I,0,1)$ is an oriented fiberwise GMF so
that $-f$ is also oriented (i.e., the fibers $M_t$ are oriented in
a coherent way), then $\g_f\oplus\g_{-f}$ is stably isomorphic to
the normal bundle of $\Sig(f)$ in $E$.
\end{prop}

\begin{rem}[normal=vertical tangent bundle]\label{FP:rem:normal bundle of Sig in E is the restriction
of the vertical tangent bundle} Since $\Sig(f)$ is the transverse
inverse image of the zero section of the vertical tangent bundle
of $E$, the normal bundle of $\Sig(f)$ in $E$ is isomorphic to the
restriction to $\Sig(f)$ of the vertical tangent bundle of $E$.
\end{rem}

\begin{proof}
Along the index $i$ Morse point set $\Sig_0^i(f)$, the normal
bundle is
\[
    \nu(\Sig_0^i(f))=\g_f^i\oplus\g_{-f}^{n-i}.
\]
If we add a trivial $n$-plane bundle $\e^n$ to $\nu(\Sig_0^i(f))$
we get an isomorphism
\[
    \f:(\g_f\oplus\g_{-f})|_{\Sig_0^i(f)}\xrarrow{\approx}\nu(\Sig_0^i(f))\oplus\e^n
\]
by sending the auxiliary vectors
\[
    v_{i+1},v_{i+2},\cdots,v_n,v^-_{n-i+1},\cdots,v^-_n
\]
to the basis vectors $e_{i+1},\cdots,e_n,e_i,\cdots,e_0$,
respectively. If we use the notation
\[
    u_j=v^-_{n-j+1}
\]
then the isomorphism $\f$ sends $u_j$ to $e_j$ and $v_k$ to $e_k$.

Near the birth-death set $\Sig_1^{i-1}(f)$ the auxiliary vectors
$v_i$ and $u_i$ are identified with $\frac{\d}{\d x_i}\in\g_f^i$
over $\Sig^i(f)$ and $-\frac{\d}{\d x_i}\in\g_{-f}^{n-i+1}$ over
$\Sig^{n-i+1}(-f)=\Sig^{i-1}(f)$, respectively. (See
Figure~\ref{FG:fig:rotate}.)
\begin{figure}
\includegraphics{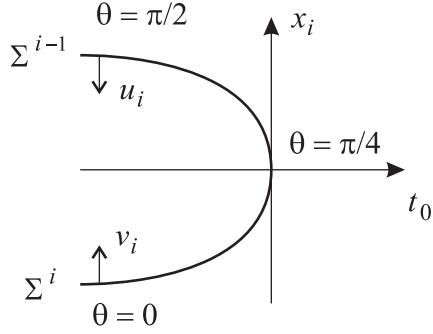}
\caption{$v_i$ rotates to $u_i$.}\label{FG:fig:rotate}
\end{figure}

The isomorphism $\f:\g_f\oplus\g_{-f}\to\nu(\Sig(f))\oplus\e^n$
can thus be defined near $\Sig_1^{i-1}$ by
\[
    \f(x,y,\sum a_ju_j,\sum b_kv_k)=x+y+(-a_i\sin\th+b_i\cos\th)\frac{\d}{\d x_i}
\]
\[
    -\sin\th\cos\th\frac{\d}{\d
    t_0}+\sum_{j=1}^{i-1}a_je_j+(a_i\cos\th+b_i\sin\th)e_i+\sum_{k=i+1}^nb_ke_k.
\]
\end{proof}

If $f$ is a fiberwise framed function then $\g_f$ is a trivial
bundle over $\Sig(f)$. If $f$ is not framed then the Framing
Principle says:
\begin{equation}\label{FP:eq:the general Framing Principle}
    \t_{2k}(E,\d_0E)=\t_{2k}(C_\ast(f))+(-1)^{k}\z(2k+1)p_\ast\left(
    \tfrac12 ch_{2k}(\g_f\otimes\CC)\right)
\end{equation}
where the push-down operator $p_\ast$ is defined below.

\subsection{Push-down/transfer}\label{FP:subsec:push-down
transfer}

In deRham cohomology, the push-down operator
\[
    p_\ast^E:H^{q+n}(E,\d E;\RR)\to H^q(B;\RR)
\]
is given by integration along the fibers which must be oriented.
(See, e.g., \cite{[BottTu]}.) Thus, if $\s$ is a smooth
$q$-simplex in $B$ and $\w$ is a closed $q+n$ form on $E$ which is
zero on $\d E$ then the value of $\w$ on $\s$ is just the integral
of its restriction to $p^{-1}\s\cong \s^q\times M^n$. It is
well-known \cite{[MoritaBook]} that this can also be given
integrally as the \emph{Gysin homomorphism}:
\begin{equation}\label{FP:eq:Gysin homomorphism}
  H^{q+n}(E,\d E;\ZZ)\twoheadrightarrow
    E_\infty^{q,n}\hookrightarrow E_2^{q,n}\cong H^q(B;H^n(M;\ZZ))
    \cong H^q(B;\ZZ)
\end{equation}
where we use the orientation of $M$ to identify $H^n(M;\ZZ)$ with
$\ZZ$.

\begin{defn}[relative Euler class]\label{FP:def:relative Euler class}
The \emph{relative Euler class}
\[
    e_{E,\d_0}=e(E,\d_0E)\in H^n(E,\d E;\ZZ)
\]
is defined to be the pull-back of the Thom class of the verical
tangent bundle $T^vE$ along any section, i.e., vertical vector
field, $v:E\to T^vE$ which points inward along $\d_0E$ and outward
along $\d_1E$. If $\d_0E$ is empty we write $e_E=e(E,\emptyset)$.
\end{defn}

The relative Euler class depends on the choice of $\d_0E$. For
example, $e(\d_0E\times I,\d_0E)=0$ but $e_{\d_0E\times I}=\Sig
e_{\d_0E}$ where $\Sig$ is the isomorphism
\[
    \Sig:H^{n-1}(\d_0E;\ZZ)\xrarrow{\approx}H^n(\d_0E\times
    I,\d_0E\times\d I;\ZZ).
\]

If we compose the Gysin homomorphism with the cup product with
$e_{E,\d_0}\in H^n(E,\d E;\ZZ)$ we get the \emph{relative transfer
homomorphism}
\begin{equation}\label{FP:eq:relative transfer for oriented fibers}
  tr^{E,\d_0}_B:H^q(E;\ZZ)\xrarrow{\cup e_{E,\d_0}}H^{q+n}(E,\d
  E;\ZZ)\xrarrow{p_\ast^E}H^q(B;\ZZ)
\end{equation}
which can be expressed in terms of the usual transfer by
\[
    tr^{E,\d_0}_B=tr^E_B-j^\ast tr^{\d_0E}_B
\]
where $j:\d_0E\to E$ is the inclusion map. The transfer does not
require oriented fibers. It is defined when the fiber is just a
finite cell complex \cite{[Becker-Gottlieb:Adams-conj]}.

In the Framing Principle (\ref{FP:eq:the general Framing
Principle}) we use the push-down operator
\begin{equation}\label{FP:eq:pushdown operator pSigma}
  p_\ast^\Sig:H^{q}(\Sig(f);\RR)\to H^{q}(B;\RR)
\end{equation}
for $q=4k$ which can be described (for all $q\geq0$) in two ways.

The first definition is in terms of differential forms. It does
not require $f$ to be a fiberwise oriented GMF. We only need $f$
to be a generic smooth function on $E$ whose fiberwise gradient
points inward along $\d_0E$ and outward along $\d_1 E$. Since
$\Sig(f)$ is normally oriented with codimension $n$, it has a dual
form $\eta_\Sig$ which is a closed $n$-form with support in a
small tubular neighborhood $T$ of $\Sig(f)$ in $E$ representing
the Thom class of the normal bundle of $\Sig(f)$ in $E$. We define
the push-down operator (\ref{FP:eq:pushdown operator pSigma}) to
be given on a closed $q$-form $\w$ on $\Sig(f)$ by:
\begin{equation}\label{FP:eq:def of push-down using forms}
  p_\ast^\Sig(\w)=p_\ast^E( \pi^\ast(\w)\wedge \eta_\Sig)
\end{equation}
where $\pi^\ast(\w)$ is the pull-back of $\w$ along the projection
$\pi:T\to\Sig(f)$.

Although (\ref{FP:eq:def of push-down using forms}) is in terms of
differential forms it actually gives an integral operation
\[
  p_\ast^\Sig:H^{q}(\Sig(f);\ZZ)\to H^{q}(B;\ZZ).
\]
To see this note that the singular set $\Sig(f)$ is the inverse
image of the zero section of the vertical tangent bundle of $E$
under a transverse mapping $\grad f$ which points inward along
$\d_0E$ and outward along $\d_1 E$. Thus the Thom class of the
normal bundle of $\Sig(f)$ in $E$, which is an integral class
\[
    [\eta_\Sig]\in H^n(T,\d T;\ZZ)
\]
maps to the relative Euler class $e_{E,\d_0}\in H^n(E,\d E;\ZZ)$
under the mapping
\[
    k^\ast:H^n(T,\d T;\ZZ)\cong H^n(E,E-\interior T;\ZZ)\to H^n(E,\d E;\ZZ)
\]
given by excision and restriction. Thus (\ref{FP:eq:def of
push-down using forms}) is the composition of the maps
\begin{equation}\label{FP:eq:integral form of 1st def of p ast}
  H^q(\Sig(f))\xrarrow{\pi^\ast} H^q(T)\xrarrow{\cup[\eta_\Sig]}
  H^{q+n}(T,\d T)\xrarrow{k^\ast}H^{q+n}(E,\d
  E)\xrarrow{p_\ast^E}H^q(B).
\end{equation}

\begin{prop}[push-down/transfer diagram]\label{FP:prop:push-down transfer diagram}
The following diagram commutes where $i^\ast,j^\ast$ are
restriction maps.
\[\begin{diagram}
\node[2]{H^q(T;\ZZ)}\arrow{e,t}{\cup[\eta_\Sig]}\node{H^{q+n}(T,\d T;\ZZ)}
    \arrow{s,r}{k^\ast}\\
    \node{H^q(\Sig(f);\ZZ)}\arrow{ne,t}{\pi^\ast}\arrow{se,b}{p_\ast^\Sig}
    \node{H^q(E;\ZZ)}\arrow{s,r}{tr^{E,\d_0}_B}\arrow{n,r}{j^\ast}
    \arrow{e,t}{\cup e_{E,\d_0}}\arrow{w,t}{i^\ast}
    \node{H^{q+n}(E,\d E;\ZZ)}\arrow{sw,b}{p_\ast^E}\\
\node[2]{H^q(B;\ZZ)}
\end{diagram}\]
\end{prop}

We now give another description of the integral push-down operator
$p_\ast^\Sig$ in the case when $f$ is a fiberwise oriented GMF.
Basically, it is the alternating sum of the push-downs on
$\Sig^i(f)$, the set of index $i$ critical points of $f$.

The projection map $p:E\to B$ gives a codimension one immersion on
the birth-death set:
\[
    p:\Sig^i_1(f)\to B.
\]
This immersion is \emph{two-sided} in the sense that its normal
bundle is oriented. For $t$ on the \emph{inside} of
$p(\Sig^i_1(f))$ the function $f_t$ has two cancelling critical
points of index $i$ and $i+1$. The \emph{outward} direction is the
normal direction in which these critical points vanish. However,
there is still a ``ghost'' of a critical point which we define to
be a lifting $Gh^i$ of an exterior immersed open collar of
$p(\Sig^i_1(f))$ in $B$. This lifting should start at the set
$\Sig^i_1(f)$ and extend outward. (See Figure~\ref{FG:fig:ghost
01}.)
\begin{figure}
\includegraphics{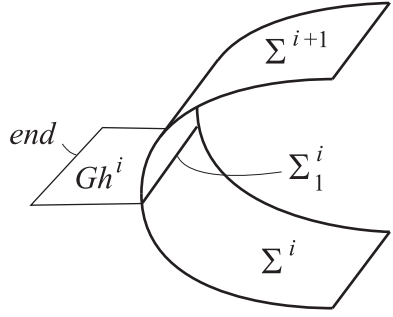}
\caption{The ghost set $Gh^i$ is attached to $\Sig^i_1(f)$.}
\label{FG:fig:ghost 01}
\end{figure}

In local coordinates we have:
\[
    f_t(x)=-x_1^2-\cdots-x_i^2+x_{i+1}^3+t_0x_{i+1}+x_{i+2}^2+\cdots+x_n^2
\]
where $dt_0$ is the outward normal direction to $p(\Sig^i_1(f))$.
For $t_0<0$ we have two critical points which converge to the
birth-death point $x=0$ at $t_0=0$. For $0<t_0<\e$ we take the
point $x=0$ to be the ghost point. The limit points (where
$t_0=\e$) form the \emph{end} of the ghost set.

Let $\widehat{\Sig^i}(f)$ be the union of $\overline{\Sig}^i(f)$
with the two exterior open collars $Gh^i$ and $Gh^{i-1}$:
\[
    \widehat{\Sig^i}(f)=\overline{\Sig}^i(f)\cup Gh^i\cup
    Gh^{i-1}.
\]
Then $\widehat{\Sig^i}(f)$ is a topological manifold with the same
dimension as $B$ and
\[
    p:\widehat{\Sig^i}(f)\to B
\]
is locally a homeomorphism. (In Figure \ref{FG:fig:ghost 01},
$Gh^i$ and $\Sig^i$ form the displayed part of
$\widehat{\Sig^i}(f)$ and $Gh^i\cup\Sig^{i+1}$ is the displayed
part of $\widehat{\Sig^{i+1}}(f)$. The ghost set $Gh^i$ is
transverse to the singular set so
$\widehat{\Sig^i}(f),\widehat{\Sig^{i+1}}(f)$ are not smooth
manifolds.)

For any $q$-cocycle $c^q$ defined in a neighborhood of $\Sig(f)$
in $E$ let $p_\ast(c)$ be the $q$-cocycle on $B$ whose value on
any sufficiently small (diameter $<\e$) singular simplex
\[
    \s:\Delta^q\to B
\]
is given by the alternating sum of the values of $c$ on the
liftings $\widetilde{\s}$ of $\s$ to $\widehat{\Sig^i}(f)$:
\begin{equation}\label{FP:eq:integral pushdown formula}
  \<p_\ast(c),\s\>=\sum_{i=0}^n(-1)^i\sum_{\widetilde{\s}}\<c,\widetilde{\s}\>.
\end{equation}
If $\widetilde{\s}:\Delta^q\to\widehat{\Sig^i}(f)$ has image in a
ghost set $Gh^i$ or $Gh^{i-1}$ then it is cancelled by another
lifting to $\widehat{\Sig^{i+1}}(f)$ or $\widehat{\Sig^{i-1}}(f)$,
respectively. Therefore, we may take the sum in
(\ref{FP:eq:integral pushdown formula}) over only those lifting
$\widetilde{\s}$ of $\s$ to $\widehat{\Sig^i}(f)$ whose image does
not lie entirely in the ghost set.

\begin{prop}[equivalence of push-downs]
\label{FP:prop:equivalence of two pushdown formulas}
These two definitions (\ref{FP:eq:def of push-down using forms}),
( \ref{FP:eq:integral pushdown formula}) of the push-down operator
(\ref{FP:eq:pushdown operator pSigma}) agree.
\end{prop}

\begin{rem}[orientation not needed]\label{FP:rem:orientation is
not needed} The \emph{integral push-down formula}
(\ref{FP:eq:integral pushdown formula}) does not require an
orientation of the fibers of $E$.
\end{rem}

\begin{proof}
This is an elementary exercise. However, since it plays a crucial
role in the proof of the Framing Principle, we will go through the
proof quickly.

First, we smooth out the sets $\widehat{\Sig^i}(f)$ as indicated
in Figure \ref{FG:fig:smooth out ghosts} so that they become
smooth manifolds locally diffeomorphic to $B$. We call these
$\widetilde{\Sig^i}(f)$.
\begin{figure}
\includegraphics{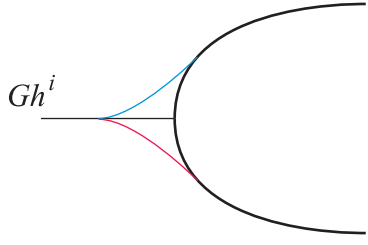}
\caption{Smoothings of {\colorr$\widehat{\Sig^i}(f)$} and
{\colorb$\widehat{\Sig^{i+1}}(f)$}.}\label{FG:fig:smooth out
ghosts}
\end{figure}
Since the smooth versions $\widetilde{\Sig^i}(f)$ of the
topological manifolds $\widehat{\Sig^i}(f)$ are vertical
deformations of these topological local sections with differences
only in an arbitrarily small neighborhood of the birth-death set
of $f$, we may use lifting of $\s$ to $\widetilde{\Sig^i}(f)$
instead of $\widehat{\Sig^i}(f)$ in the formula
(\ref{FP:eq:integral pushdown formula}) without changing the
cohomology class of $p_\ast(c)$. As before, the terms in which
$\widetilde{\s}$ lie entirely in the ghost set cancel.

Now suppose that the cocycle $c^q$ in (\ref{FP:eq:integral
pushdown formula}) is given by integration of the $q$-form
$\pi^\ast\w$ in (\ref{FP:eq:def of push-down using forms}). Since
each $\widetilde{\Sig^{i}}(f)$ is a smooth submanifold of $E$, it
has a dual $n$-form $\eta_i$ which is defined in the complement of
the ends of the ghost sets $Gh^i$ and $Gh^{i-1}$ and which has
support in $T$ (the tubular neighborhood of $\Sig(f)$). This form
has the property that
\begin{equation}\label{FP:eq:c(s) in terms of omega eta-i}
  \sum_{\widetilde{\s}}\<c,\widetilde{\s}\>=\int_{{\s}^\ast
  E}\pi^\ast\w\wedge\eta_i+\<b_i,\d\widetilde{\s}\>
\end{equation}
where $b_i$ is a $q-1$ form which is equal to $b_{i+1}$ on
$\widetilde{\Sig^{i}}(f)\cap\widetilde{\Sig^{i+1}}(f)$.

Since the manifolds $\widetilde{\Sig^{i}}(f)$ and
$\widetilde{\Sig^{i+1}}(f)$ agree near the end of $Gh^i$, the
forms $\eta_i$ and $\eta_{i+1}$ can be chosen to agree near the
end of $Gh^i$. Consequently, the alternating sum
\[
    \eta=\sum_{i=0}^n(-1)^i\eta_i
\]
is a well-defined $n$-form on all of $E$ and it has support in the
tubular neighborhood $T$ of $\Sig(f)$. Combining (\ref{FP:eq:c(s)
in terms of omega eta-i}) and (\ref{FP:eq:integral pushdown
formula}) we get
\begin{equation}\label{FP:eq:p(c),s in terms of omega}
    \<p_\ast(c),\s\>=\int_{\s^\ast
    E}\pi^\ast\w\wedge\eta+\<b,\d\s\>
\end{equation}
where the value of $b^{q-1}$ on a small singular $q-1$ simplex
$\t$ in $B$ is given by the alternating sum of the values of $b_i$
on the liftings $\widetilde{\t}$ of $\t$ to
$\widetilde{\Sig^i}(F)$:
\[
    \<b,\t\>=\sum_{i=0}^n(-1)^i\sum_{\widetilde{\t}}\<b_i,\widetilde{\t}\>.
\]

The bundle map $\pi:T\to\Sig(f)$ can be chosen to be vertical
(parallel to the fibers of $E$) except near the birth-death set
$\Sig_1(f)$. Consequently, we may choose $\eta_\Sig$ to be equal
to $\eta$ outside a tubular neighborhood of $\Sig_1(f)$. Since
$\Sig_1(f)$ has codimension $n+1$, the closed $n$-forms
$\eta_\Sig$ and $\eta$ must differ by an exact form $d\a$. So
\[
    \int_{\s^\ast
    E}\pi^\ast\w\wedge\eta=\int_{\Delta^q}\s^\ast
    p_\ast^E(\pi^\ast\w\wedge\eta)
    =\int_{\Delta^q}\s^\ast
    p_\ast^E(\pi^\ast\w\wedge\eta_\Sig)+\<a,\d\s\>
\]
where $a$ is given by integrating $p_\ast^E(\pi^\ast\w\wedge\a)$.

Combining this with (\ref{FP:eq:p(c),s in terms of omega}) we see
that the $q$-cocycle $p_\ast(c)-\delta(a+b)$ is given by
integrating the $q$-form $p_\ast^\Sig(\w)$ of (\ref{FP:eq:def of
push-down using forms}).
\end{proof}

\subsection{The Framing Principle}\label{FP:subsec:framing
principle}

We can now state the general Framing Principle.

\begin{thm}[Framing Principle]\label{FP:thm:Framing Principle}
Suppose that $p:E\to B$ is a smooth bundle with compact fiber $M$
and $\F$ is an $r$-dimensional Hermitian coefficient system over
$E$ so that either
\begin{enumerate}
  \item[a)]$\pi_1B$ acts trivially on $H_\ast(M,\d_0M;\F)$ or
  \item[b)]$r=1$ and $H_\ast(M,\d_0M;\F)$ is $\pi_1B$-upper triangular.
\end{enumerate}
Suppose also that $f:(E,\d_0)\to (I,0)$ is an oriented fiberwise
GMF. Then the higher relative FR-torsion invariants
$\t_\ast(E,\d_0E;\F)$ of $(E,\d_0E)$ with coefficients in $\F$
differs from the higher torsion given by $f$ in the following way.
\[
    \t_{2k}(E,\d_0E;\F)=\t_{2k}(C_\ast(f;\F))+r(-1)^k\z(2k+1)p_\ast^\Sig
    \left(\tfrac12 ch_{2k}(\g_f\otimes\CC)\right)
\]
where $\g_f$ is the negative eigenspace bundle over $\Sig(f)$
(Definition~\ref{FP:def:gamma on Sig(f)}) and $p_\ast^\Sig$ is the
push-down operator defined in the previous section. The correction
term is zero in the other degrees:
\[
    \t_{2k+1}(E,\d_0E;\F)=\t_{2k+1}(C_\ast(f;\F)).
\]
\end{thm}

%% file: ProofFP.tex
\vfill\eject\section{Proof of the Framing
Principle}\label{Proof:section}

\begin{enumerate}
  \item Transfer theorem.
  \item Stratified deformation lemma.
  \item Proof of transfer theorem.
  \item Proof of Framing Principle.
\end{enumerate}

This section contains the proof of both the general Framing
Principle and the transfer theorem for higher torsion. These
theorems are closely related and have the same proof. We will
first state and prove the transfer theorem. Then we modify the
prove of the transfer theorem to give a proof of the Framing
Principle.

\subsection{Transfer theorem}\label{Proof:subsec:transfer}

Assume as before that $p:E\to B$ is a smooth bundle with compact
fiber $M$ where $B$ is a closed manifold. Then we define the
\emph{transfer}
\[
    tr^E_B:H^q(E;\ZZ)\to H^q(B;\ZZ)
\]
{from $E$ to $B$} as follows. If there exists a fiberwise framed
function $f:E\to \RR$ we let $tr^E_B$ be the composition
\[
    H^q(E;\ZZ)\xrarrow{i^\ast}H^q(\Sig(f);\ZZ)\xrarrow{p_\ast^\Sig}H^q(B;\ZZ)
\]
of restriction to $\Sig(f)$ and the integral push-down operator
$p_\ast^\Sig$ given in (\ref{FP:eq:integral pushdown formula}). If
the fibers of $E$ are oriented, this is the same as the previously
defined transfer map by the push-down/transfer diagram
(Proposition~\ref{FP:prop:push-down transfer diagram}).

If $E$ is not fiberwise oriented and there is no fiberwise framed
function on $E$ then take a linear disk bundle $D^N\to D\to E$
where $N$ is large enough to admit a fiberwise framed function
$g:D\to\RR$ and pull back the transfer from $D$:
\[
    H^q(E;\ZZ)\cong H^q(D;\ZZ)\xrarrow{tr^D_B}H^q(B;\ZZ).
\]

To see that this is independent of the choice of $D$, suppose that
$D'$ is another linear disk bundle over $E$. Then $g$ gives a
fiberwise framed function on the fiber product of $D$ and $D'$
over $E$:
\[
    h:D\times_E D'\to\RR,\quad h(x,y)=g(x)+\length{y}^2.
\]
But $\Sig(h)=\Sig(g)\times0$ so the transfer from $D$ to $B$ is
the pull-back of the transfer from $D\times_ED'$ to $B$. By
symmetry the same holds for $D'$ so $tr^D_B$ and $tr^{D'}_B$ pull
back to the same mapping on $H^q(E;\ZZ)$.

Suppose that $D$ is a smooth bundle over $E$ with fiber $Y$.
Suppose that $\F$ is a Hermitian coefficient system on $D$ so that
$H_\ast(Y;\F)=0$. Let $F$ be the fiber of $D$ over $B$. Then $F$
is also a bundle over $M$ with fiber $Y$.
\[
\begin{diagram}
\node{Y}\arrow{e}\arrow{s,l}{=}\node{F}\arrow{e}\arrow{s}\node{M}\arrow{s}\\
\node{Y}\arrow{e}\node{D}\arrow{e}\arrow{s}\node{E}\arrow{s}\\
\node[2]{B}\arrow{e,t}{=}\node{B}
\end{diagram}
\]
Then $H_\ast(F;\F)=0$ so both terms in the following formula are
define.
\begin{thm}[transfer theorem]\label{Proof:thm:transfer theorem}
Let $\t_k(D;\F)_E\in H^{2k}(E;\F)$ denote the higher torsion of
$D$ as a bundle over $E$ and similarly for $E$ replaced with $B$.
Then
\begin{equation}\label{Proof:eq:transfer formula}
    \t_k(D;\F)_B=tr^E_B(\t_k(D;\F)_E).
\end{equation}
\end{thm}

In the special case when $D$ is a fiber product $D=E\times_BE'$ we
get the following.

\begin{cor}[fiber product formula]\label{Proof:cor:fiber product
formula} Let $E'\to B$ be a smooth bundle with fiber $Y$. Let $\F$
be a Hermitian coefficient system on $E'$ so that $H_\ast(Y;\F)=0$
and let $\F'$ be the pull-back of $\F$ to the fiber product
$E\times_BE'$. Then
\[
    \t_k(E\times_BE';\F')=\chi(M)\t_k(E';\F).
\]
\end{cor}

\begin{proof}
By the transfer theorem we have
\[
    \t_k(E\times_BE';\F')_B=tr^E_B(\t_k(E\times_BE';\F')_E)=tr^E_Bp^\ast(\t_k(E';\F)).
\]
Then, use the following well-known fact which follows trivially
from the integral push-down formula.
\end{proof}

\begin{lem}[transfer of pull-back]\label{Proof:lem:mult by chi(M)}The composition
\[
    H^q(B;\ZZ)\xrarrow{p^\ast} H^q(E;\ZZ)\xrarrow{tr^E_B} H^q(B;\ZZ)
\]
is multiplication by $\chi(M)$.
\end{lem}

\subsection{Stratified deformation lemma}\label{Proof:subsec:stratified deformation lemma}

The first step in the proof of the transfer theorem is to
stabilize so that we can use the Framed Function Theorem. To make
$\dim M>\dim B$ we replace $E$ with a disk bundle $E'$ with the
corners rounded. Then we have an isomorphism \[
    \f^\ast:H^{2k}(E)\xrarrow{\approx}H^{2k}(E')
\]with any coefficients. Let $D'$ be the pull-back of $D$ to $E'$. Then
\[
    tr^{E'}_B\t_k(D';\F')_{E'}=tr^{E'}_B\f^\ast(\t_k(D;\F)_E)=tr^{E}_B(\t_k(D;\F)_E)
\]
and
\[
    \t_k(D';\F')_B=\t_k(D;\F)_B
\]
by the product formula (Lemma~\ref{sec2:lem:product formula}).

To make $\dim Y>\dim B$ we take the product of $Y$ with an even
dimensional sphere. The torsion of $D\times S^{2N}$ is twice the
torsion of $D$ over both $E$ and $B$ so the transfer equation for
$D\times S^{2N}$ is equivalent to the one for $D$.

Now that $\dim M>q=\dim B$, there is a fiberwise framed function
\[
    f:E\to\RR.
\]
The fiberwise singular set $\Sig=\Sig(f)$ is a smooth closed
$q$-submanifold of $E$ which is \emph{stratified} by index and
\emph{codimension}. (Birth-death points have codimension $1$.
Morse points have codimension $0$.) The restriction of $p:E\to B$
to $\Sig$ gives a smooth mapping $\Sig\to B$ which is a local
diffeomorphism on the codimension $0$ set and has \emph{fold
singularities} along the codimension $1$ set $\Sig_1=\Sig_1(f)$.
In other words, it has the local form
\[
    (x_1,x_2,\cdots,x_q)\mapsto(x_1^2,x_2,\cdots,x_q)
\]
where $\Sig_1$ is given by $x_1=0$, $\Sig^i$ is given by
$x_1\leq0$ and $\Sig^{i+1}$ is given by $x_1>0$.

By a \emph{stratified deformation} of $\Sig$ over $B$ we mean a
compact $q+1$ dimensional smooth manifold $W$ over $B\times I$
stratified by index and codimension as described above so that
$\d_0W=W|B\times0$ and $\d_1W=W|B\times1$ are stratified closed
manifolds over $B$ and $\d W=\d_0W\coprod\d_1W$.

The restriction of $Y\to D\to E$ to $\Sig$ is a smooth bundle
\[
    Y\to D|\Sig\to\Sig
\]
which is classified by a continuous map
\[
    \c_\Sig:\Sig\to BDif\!f(Y;\F)
\]
where $Dif\!f(Y;\F)$ is the group of diffeomorphisms $\f$ of $Y$
together with an automorphism $\widetilde{\f}$ of the Hermitian
coefficient system $\F$ covering $\f$ (i.e.,
$\widetilde{\f}:\F\cong\f^\ast\F$).

Note that $\c_\Sig$ is the restriction to $\Sig$ of the mapping
\[
    \c_E:E\to BDif\!f(Y;\F)
\]
which classifies the bundle $Y\to D\to E$. The higher torsion
class $\t_k(D;\F)_E$ is the pull-back along $\c_E$ of the
\emph{universal higher torsion class}
\[
    \t_k\in H^{2k}(BDif\!f(Y;\F);\RR).
\]

By a \emph{stratified deformation} of the \emph{singular pair}
$(\Sig,\c_\Sig)$ we mean a stratified deformation $W$ of $\Sig$ as
above together with a continuous mapping
\[
    {\c_W}:W\to BDif\!f(Y;\F)
\]
extending $\c_\Sig$.

\begin{lem}[stratified deformation lemma]
\label{Proof:lem:stratified deformation lemma} Both sides of the
transfer formula (\ref{Proof:eq:transfer formula}) are functions
of the stratified deformation class of the singular pair
$(\Sig,\c_\Sig)$. Furthermore, they are both additive with respect
to disjoint union of such pairs.
\end{lem}

\begin{rem}[Becker-Gottlieb
transfer]\label{Proof:rem:Becker-Gottlieb} This lemma can be
interpreted as saying that both sides of the transfer formula
depend only on the homotopy class of the mapping:
\[
    B\to Q(E_+)\to Q(BDif\!f(Y;\F)_+)
\]
where the first map $B\to Q(E_+)$ is the \emph{Becker-Gottlieb
transfer} \cite{[Becker-Gottlieb:Adams-conj]}.
\end{rem}

\begin{proof}
\underline{Case 1} (algebraic torsion). First consider the right
hand side (RHS) of the transfer formula (\ref{Proof:eq:transfer
formula}). By definition of the transfer $tr^E_B$, the RHS is
given by the push-down of the pull-back of the universal torsion
class $\t_k$:
\[
    tr^E_B(\t_k(D;\F)_E)=p_\ast^E\c_E^\ast(\t_k)=p_\ast^\Sig\c_\Sig^\ast(\t_k).
\]
We call this the \emph{algebraic torsion} of $(\Sig,\c_\Sig)$. To
easily see that this is a stratified deformation invariant, let
$(W,\c_W)$ be any stratified deformation of $(\Sig,\c_\Sig)$. Then
\[
    p_\ast^\Sig\c_\Sig^\ast(\t_k)=p_\ast^W\c_W^\ast(\t_k)|B\times0
    =p_\ast^W\c_W^\ast(\t_k)|B\times1
\]
where $p_\ast^W:H^{2k}(W;\RR)\to H^{2k}(B\times I;\RR)$ is given
by the integral push-down formula (after attaching an ghost set
along $W_1$). Additivity with respect to disjoint union is
obvious.

\underline{Case 2} (topological torsion). Now take the left hand
side (LHS) of the transfer formula. This side is defined
topologically, using the Framed Function Theorem.

\underline{Case 2a} ($f$ is Morse). Suppose first that $\Sig_1$ is
empty, i.e., $f$ has no birth-death singularities. Then the
statement follows from Theorem~\ref{sec2:thm:transfer for twisted
Morse bundles}. We review the argument briefly since we need to
generalize it.

We took tubular neighborhoods $T^i\subseteq E$ of the singular
sets $\Sig^i(f)$ and constructed a fiberwise framed function
$h:D\to \RR$ so that $D|T^i$ is a convex set equivalent to the
main example (\ref{sec2:eg:main example}). By the local
equivalence lemma (\ref{sec2:lem:local equivalence}), this
identifies a subquotient functor of the cellular chain complex
functor for $h$. Since $Y$ is acyclic, the local splitting lemma
(\ref{sec2:lem:local splitting lemma}) applies. So,
$C_\ast(h;\F)$, whose torsion is the LHS of
(\ref{Proof:eq:transfer formula}), locally decomposes as a direct
sum of subquotient functors so it is isomorphic to the covering
space push-down of the suspensions of the cellular chain complex
of the fiberwise framed function $g:D|\Sig\to\RR$.
\[
    C_\ast(h;\F)\simeq\bigoplus
    p_\ast^\Sig(\Sig^i(C_\ast(g|\Sig^i))).
\]
Thus, by the suspension theorem (\ref{sec2:thm:suspension
theorem}), we have
\[
    \t_k(D;\F)_B=\t_k(C_\ast(h;\F))=p_\ast^\Sig\left(\sum(-1)^i
    \c_\Sig^\ast(\t_k)\right).
\]
This is a stratified deformation invariant, proving the lemma (and
the theorem) in this case.


\underline{Case 2b} ($f$ has b-d points.) In this case we use two
kinds of convex sets. We envelop the b-d sets $\Sig_1^i(f)$ in
subsets $L^i$ and the Morse points $\Sig_0^i(f)$ in subsets $T^i$
as in the proof of Case 2a (Theorem~\ref{sec2:thm:transfer for
twisted Morse bundles}). (See Figure~\ref{Proof:fig:bd transfer}.)
\begin{figure}
\includegraphics{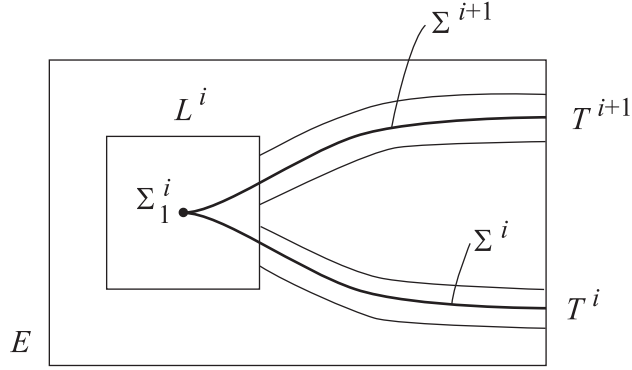}
\caption{$L^i$ (square) envelops the b-d set
$\Sig_1^i$.}\label{Proof:fig:bd transfer}
\end{figure}

Let $L^i$ be a small tubular neighborhood in $E$ of the index $i$
birth-death set $\Sig_1^{i}=\Sig_1^{i}(f)$. Let
$\pi:L^i\to\Sig_1^i$ be the projection map. Since $f$ is fiberwise
framed with index $i$ along $\Sig_1^{i}$, we can choose $L^i$ to
have a twisted product structure
\[
    L^i=\Sig_1^{i}\times J\times \e D^i\times_{\Sig_1^i}\e
    D(\g')\times J
\]
where $\e D^i$ is the $i$-disk of radius $\e$ and $\e D(\g')$ is
the radius $\e$ disk bundle of the unstable bundle of $f$ over
$\Sig_1^{i}$ and $J=\e D^1=[-\e,\e]$.

Then $L^i$ is a bundle over the codimension $0$ immersed
submanifold $\Sig_1^i\times J$ of $B$ where the fiber over
$(t,u)\in\Sig_1^i\times J$ is
\[
    L^i_{t,u}=\e D^i\times\e D(\g')\times J.
\]
Using the parametrized Morse lemma we may assume that $f|T$ is
given by
\[
    f_{t,u}(x,z,w)=c_t-\length{x}^2+\length{z}^2+w^3-uw.
\]

Since $\dim Y>q=\dim\Sig$, there is a fiberwise framed function
$g:D|\Sig_1\to [0,3]$ which we write as a family of framed
functions
\[
    g_t:(F_t,\d F_t)\to ([0,3],3)
\]
on the fibers $F_t$ of $D$ over $t\in\Sig_1$. We can assume that
the critical values of $g_t$ lie in the open interval $(1,2)$ for
all $t\in\Sig_1$. We can extend $g_t$ to $L=\coprod L^i$ by
composing with the projection $\pi:L^i\to\Sig_1^i$. Then we can
extend to a fiberwise framed function over the entire singular set
$\Sig=\Sig(f)$ and then smoothly to all $t\in E$ but $g_t$ will
only be framed in a small neighborhood of $\Sig$.

The fiber of $D$ over $(t,u)\in \Sig_1^i\times J$ is
\[
    D_{t,u}=Y_t\times L_{t,u}^i=Y_t\times \e D^i\times\e D(\g')\times J.
\]
Thus $h:D\to\RR$ can be given by $h=\delta g+f$, i.e.,
\[
    h_{t,u}(y,x,z,w)=\delta g_t(y)+c_t-\length{x}^2+\length{z}^2+w^3-uw.
\]
If $\delta$ is sufficiently small, $h_{t,u}$ will not have
singularities outside of the $D|\Sig$.

As in Example~\ref{sec2:eg:second example}, we need to modify the
function $h_{t,u}$ in a small neighborhood of the bad set where
$u,x,z,w$ are all zero. Since we have stabilized our bundles, we
can keep $x,z$ fixed and deform the function
\[
    h_{t,u}(y,0,0,w)=\delta g_t(y)+c_t+w^3-uw
\]
so that it is fiberwise framed. Adding back the terms
$-\length{x}^2+\length{z}^2$, we get a fiberwise framed function
which is locally equivalent to Example~\ref{sec2:eg:second
example}.

On the rest of the singular set $\Sig=\Sig(f)\subseteq E$ we
construct tubular neighborhoods $T^i$ as in the proof of Case 2a
(Theorem~\ref{sec2:thm:transfer for twisted Morse bundles}).

Since $Y$ is acyclic, the expansion functor given by $h$ locally
decomposes as a direct sum of the local subquotients. These in
turn depend only on the stratified set $\Sig$ together with the
bundle $D|\Sig$ and the coefficient system $\F$. In other words,
it depends only on the singular pair $(\Sig,\psi_\Sig)$. A
stratified deformation of this pair carries enough information to
deform the convex sets and therefore determines the local
subquotients. The lemma therefore holds.
\end{proof}

\subsection{Proof of transfer theorem}\label{Proof:subsec:proof of
transfer formula}

The proof of the transfer theorem now proceeds in three steps
\begin{enumerate}
  \item[a)] Both sides of the transfer formula (\ref{Proof:eq:transfer
  formula}) are zero if $\c_\Sig:\Sig\to BDif\!f(Y;\F)$ is
  trivial.
  \item[b)] The transfer formula holds if the $\c_\Sig$ is trivial
  on the b-d set $\Sig_1$.
  \item[c)] The singular pair $(\Sig,\c_\Sig)$ can be deformed by a
  stratified deformation to a pair which lies in two indices
$i,i+1$ and satisfies (b).
\end{enumerate}

To clarify what we mean by ``trivial,'' we choose a base point
(each one called $\ast$) for every component of $BDif\!f(Y;\F)$.
By \emph{trivial} we mean ``Each component (of $\Sig$ or $\Sig_1$)
maps to the base point of the corresponding component of
$BDif\!f(Y;\F)$.''

The logical order of the proof of these statements is:
(c),(a),(b), but we discuss them in the order listed. So, suppose
for a moment that we have already proved (c). Then (a) follows
from Lemma~\ref{Proof:lem:stratified deformation lemma}. The
reason is that, with $\c_\Sig$ trivial, both sides of the transfer
formula are functions of the stratified deformation class of
$\Sig$ which is equivalent to a mapping of $B$ into $Q(S^0)$.
Since $Q(S^0)$ is rationally trivial (except for
$\pi_0(Q(S^0))=\ZZ$ which gives the Euler characteristic of $Y$
which must be zero since it is acyclic), some positive integer
multiple of $(\Sig,\ast)$ must be stratified null-deformable.

To prove the statement in the last paragraph, we use (c) which
says that we may assume that $\Sig$ is in the two indices $i,i+1$.
We view $\Sig^i,\Sig^{i+1}$ as positive and negative particles,
depending on the signs $(-1)^i$, $(-1)^{i+1}$. Then, $\Sig$ gives
a mapping from $B$ into the configuration space of positive and
negative particles which are allowed to cancel two at a time. This
configuration space is known to be homotopy equivalent to $Q(S^0)$
\cite{[Caruso81]}. Caruso also has a short proof of the unstable
version of this statement.

Now we prove (b) assuming (a). The proof is in two steps. First,
we show that the transfer formula holds for what we call
``immersed lenses.'' Then, we show that the singular pair
$(\Sig,\c_\Sig)$ has a stratified deformation to a linear
combination of such ``immersed lenses'' plus a singular pair with
trivial $\c_\Sig$.

Suppose that $V$ is a compact connected $q$-manifold with boundary
which is immersed in the closed $q$-manifold $B$. Let
$\c_1,\c_2:V\to BDif\!f(Y;\F)$ be two mappings which are trivial
(map to the base point $\ast$ of the same component) on $\d V$.
Then the \emph{immersed lens}
\[
    L_i(V,\c_1,\c_2)=(L,\c_L)
\]
is given by letting $L$ be the double of $V$ (two copies glued
along $\d V$) with indices $i,i+1$ where $\c_L$ is equal to $\c_1$
on $L^i$ and $\c_2$ on $L^{i+1}$.

\begin{lem}[immersed lenses]\label{Proof:lem:immersed lenses}
The transfer formula holds for immersed lenses.
\end{lem}

\begin{proof}
First, consider the case when $V$ is embedded in $B$. Then we can,
by a stratified deformation, eliminate all the b-d points of $L$.
(Just double the set $V\times I\cup B\times[\frac12,1]$ with
corners rounded.) By Theorem~\ref{sec2:thm:transfer for twisted
Morse bundles}, the transfer formula holds.

For the general case, we use the push-down homomorphism
\[
    H^{2k}(V,\d V;\RR)\to H^{2k}(B;\RR)
\]
given by the same formula as the transfer for covering maps
(\ref{sec2:eq:transfer for covering maps}). The algebraic and
topological torsion of the singular pair $(L,\c_L)$ are push-downs
of the corresponding torsion classes in the relative cohomology of
$(V,\d V)$ given by the identity embedding $V\to V$. Therefore, it
suffices to show that the transfer formula holds for an embedded
lens in the relative case.

Now, we view the singular pair $(L,\c_L)$ over $(V,\d V)$ as a
mapping
\[
    V/\d V\to Q(BDif\!f(Y;\F)_+).
\]
By changing this mapping on the base point we can eliminate the
birth-death points. This reduces us to Case 2a where the transfer
formula already holds.
\end{proof}

To complete the proof of (b), we need the following lemma. The
wording of the lemma assumes that we are working in the abelian
group generated by singular pair with addition given by disjoint
union. The additivity of both sides of the transfer formula under
disjoint union allows us to do this.

\begin{lem}[decomposition into
lenses]\label{Proof:lem:decomposition into lenses} Any singular
pair $(\Sig,\c_\Sig)$ over $B$ so that $\c_\Sig$ is constant on
the b-d set $\Sig_1$ is equivalent to a linear combination of
immersed lenses and a trivial pair (where $\c_\Sig$ is trivial).
\end{lem}

\begin{proof}
The singular pair $(\Sig,\c_\Sig)$ will have a finite number of
Morse components on which $\c_\Sig$ is nontrivial. If this number
is zero then the pair is trivial. So, it suffices to reduce this
number by one by adding and deleting immersed lenses.

Let $V$ be a component of $\Sig^i$ on which $\c_\Sig$ is
nontrivial. So
\[
    \c_V=\c_\Sig|V:V/\d V\to BDif\!f(Y;\F)
\]
is essential. Then we can add a trivial immersed lens
$L_{i+1}(V,\ast,\ast)$ (being stratified null deformable) and
cancel along a parallel copy of $\d V$ in the interior of $V$ as
indicated in Figure~\ref{Proof:fig:lenses}. This isolates $\c_V$
in an immersed lens isomorphic to $L_i(V,\c_V,\ast)$. Remove this
lens and replace it with a trivial lens $L_i(V,\ast,\ast)$. Then
go backwards in the deformation of Figure~\ref{Proof:fig:lenses}.
\begin{figure}
\includegraphics{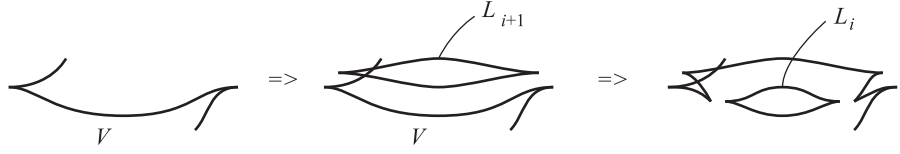}
\caption{Isolate $\c_V$ in a lens.} \label{Proof:fig:lenses}
\end{figure}

The end result is that $\Sig$ is unchanged and $\c_\Sig$ has been
trivialized on $V$.
\end{proof}

The last step (c) comes logically before (a) and (b). The proof of
(c) will be in two steps. First, we deform any singular pair
$(\Sig,\c_\Sig)$ into two indices. Next, we deform the singular
pair so that each component of $\Sig_1$ lies in a contractible
subset of $\Sig$. It will then follow that the map $\c_\Sig$ can
be changed by a homotopy so that it is trivial on $\Sig_1$.

To deform the singular pair $(\Sig,\c_\Sig)$ into two indices we
use ``{twisted lenses}'' which are defined as follows. Let $V$ be
a compact $q$-manifold with boundary which is immersed in $B$ and
let $\c_V:V\to BDif\!f(Y;\F)$ be any continuous mapping. Then the
\emph{twisted lens} is given by
\[
    L_i(V,\c_V,\c_V)=(L,\c_L)
\]
where $L=L^i\cup L^{i+1}$ is the double of $V$ and $\c_L:L\to
BDif\!f(Y;\F)$ is given by $\c_V$ on both $L^i$ and $L^{i+1}$. In
other words it is the composition of the folding map $L\to V$ with
$\c_V$.

\begin{figure}
\includegraphics{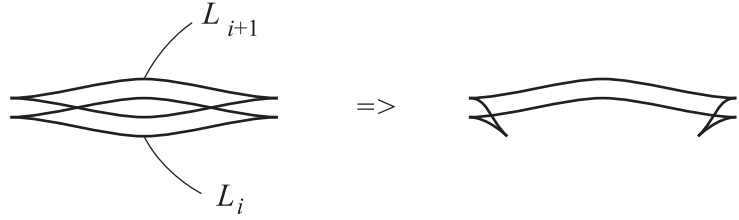}
\caption{Cancel the bottoms of two twisted lenses.}
\label{Proof:fig:lenses2}
\end{figure}

It is obvious that every twisted lens is stratified null
deformable. So, we can add and delete them at will. In particular,
we can add two twisted lenses $L_i(V,\c_V,\c_V)$ and
$L_{i+1}(V,\c_V,\c_V)$. The bottom portion of each of these
twisted lenses can be cancelled by a stratified deformation to
form a singular pair as shown in Figure~\ref{Proof:fig:lenses2}.

\begin{figure}
\includegraphics{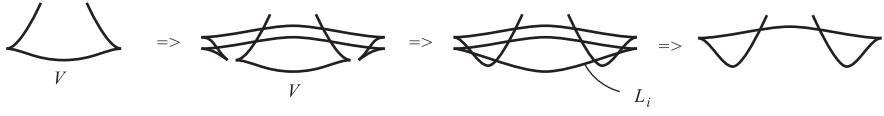}
\caption{Place $V$ in twisted lens $L_i$ and eliminate.}
\label{Proof:fig:lenses3}
\end{figure}

Now let $i$ be minimal so that $\Sig^i$ is nonempty. Let $V$ be a
component of $\Sig^i$. Let $\c_V=\c_\Sig|V$. Use this to construct
the singular pair in Figure~\ref{Proof:fig:lenses2}. Then, along
$\d V$, the mappings to $BDif\!f(Y;\F)$ agree so we can perform
the stratified deformation shown in
Figure~\ref{Proof:fig:lenses3}. This places $V$ inside a new
twisted lens $L_i(V,\c_V,\c_V)$ which we can eliminate. The result
is that the number of components of $\Sig^i$ decreases by one.
(And a new component of $\Sig^{i+1}$ is introduced.) By induction,
we can deform the singular pair into two indices.

Assume now that $\Sig=\Sig^i\cup \Sig^{i+1}$. The final step in
the proof of the transfer theorem is to show that $\c_\Sig$ can be
trivialized on $\Sig_1$. In fact, we will arrange for every
component of $\Sig^i$ to be contained in a contractible subset of
$\Sig$. Since $\Sig_1=\d\Sig^i$, this will prove the theorem.

Choose a triangulation of $\Sig^i$ so that each simplex maps
monomorphically into $B$. Then, we will ``cut apart'' the set
$\Sig^i$ by deleting a tubular neighborhood of each simplex
starting at the lowest dimension. If we let $S$ denote what is
left of the original set $\Sig^i$, then, at the end, the set $S$
will be a disjoint union of $q$-disk ($q=\dim B$).

At each step of the deformation, a new component of $\Sig^i$ will
be introduced which will be contained in a $q$-disk disjoint from
the set $S$. Every step of the deformation will alter the set
$\Sig$ only in an arbitrarily small neighborhood of the set $S$.
Consequently, the $q$-disks containing the new components of
$\Sig^i$ will not be affected.

Now we describe the deformation. The deformation is almost the
same as the one given in \cite{[I:FF]} and illustrated in Figures
E,F,G,H,I,J in \cite{[I:FF]}. (But here we have $i+1$ instead of
$i-1$.) Suppose first that $v$ is vertex in the interior of
$S\subseteq\Sig^i$. Then we introduce a trivial lens $L$ over the
image of $v$ in $B$ mapping to the point $\c_\Sig(v)$. Then $L^i$
is a new component of $\Sig^i$ which is contractible. We will not
touch this set again. Instead we cancel part of $L^{i+1}$ with a
$q$-disk neighborhood of $v$ in $S$.

Suppose by induction that all internal simplices of $S$ of
dimension $<m$ have been removed. Let $D^m$ be an $m$-disk
embedded in $S$ with boundary $S^{m-1}\subseteq\d S$. Suppose for
a moment that the sphere $S^{m-1}$ also bounds an $m$-disk
$\Delta$ in $\Sig^{i+1}$ which lies over the image of $D^m$ in
$B$. Suppose also that $D^m\cup\Delta$ forms a twisted lens, i.e.,
the restrictions of $\c_\Sig$ to $\Delta$ and $D^m$ agree if we
compose with a diffeomorphism $\Delta\cong D^m$. Then we can
cancel $\Delta$ with $D^m$ eliminating this $m$-disk from $S$.

\begin{figure}
\includegraphics{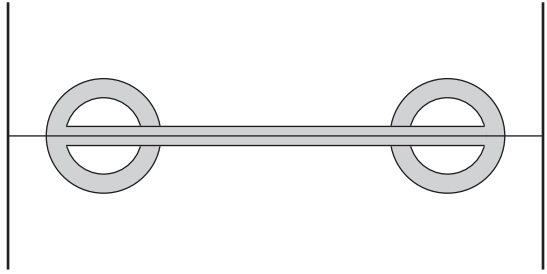}
\caption{New component of $\Sig^i$ is shaded.}
\label{Proof:fig:surgery}
\end{figure}

In general, we do not have such a disk $\Delta$. We need to
construct it. We take a sphere $S^{'m-1}$ in the interior of
$\Sig^{i+1}$ which is parallel to $S^{m-1}$ and lies over the
image of $D^m$ in $B$. Over that sphere, we create a ``{tube}''
$T$, a product of an $m-1$ sphere with a trivial lens of dimension
$q-m+1$. Then $T^i\cong T^{i+1}\cong D^{q-m+1}\times S^{m-1}$. The
mapping to $BDif\!f(Y;\F)$ should agree with $\c_\Sig|S^{'m-1}$.
Then we can cancel a tubular neighborhood of $S^{'m-1}$ with a
tubular neighborhood of the core $\ast\times S^{m-1}$ of $T^i$.
Then the new component of $\Sig^i$ will be equivalent to
\[
    S^{m-1}\times(D^{q-m+1}-\tfrac12 D^{q-m+1}).
\]
This is shown as two shaded annuli in
Figure~\ref{Proof:fig:surgery}.

Along the inner edge of $\Sig_1$ given by $S^{m-1}\times\frac12
e_1$ where $e_1$ is the unit vector in the first coordinate axis,
we do surgery on this $m-1$ sphere producing an $m-1$ disk of
points in both $\Sig^i$ and $\Sig^{i+1}$. This is shown as a
horizontal shaded bar in Figure~\ref{Proof:fig:surgery}. We extend
the mapping $\c_\Sig$ along this disk so that it matches
$\c_\Sig|D^m$. This gives a parallel disk $\Delta$ as described in
the previous paragraph. ($\Delta$ is shown as a line bisecting the
shaded region in Figure~\ref{Proof:fig:surgery}. However, it lies
in $\Sig^{i+1}$, not in $\Sig^i$.)

By construction, we can cancel $\Delta$ with $D^m$. The
deformation alters the original set $\Sig$ only in a small
neighborhood of $D^m$ since half of it occurs in the set $\Delta$
which is newly added.

The new component of $\Sig^i$ lies in a set diffeomorphic to
$D^{q-m+1}\times S^{m-1}$ with an $m$-handle added, making it
contractible. (In Figure~\ref{Proof:fig:surgery}, this
contractible set is given by adding to the shaded region the two
disks of $\Sig^{i+1}$ components with boundary the two outer
circles. Spin this about the central vertical axis to get the case
$m=2$.) Thus, the deformation works in the required way, proving
(c) and thus the transfer theorem.

\subsection{Proof of Framing
Principle}\label{Proof:subsec:proof of FP}

Now we are ready to prove the Framing Principle. Suppose that
$p:E\to B$ is a smooth bundle with $\d E=\d_0E\coprod\d_1E$ and
$\F$ is a r-dimensional Hermitian coefficient system over $E$ so
that the relative torsion $\t_k(E,\d_0E;\F)$ is defined, i.e.,
either
\begin{enumerate}
  \item $\pi_1B$ acts trivially on $H_\ast(M;\F)$ or
  \item $r=1$ and $H_\ast(M;\F)$ is $\pi_1B$-upper triangular.
\end{enumerate}
Let
\[
    f:(E,\d_0E)\to(I,0)
\]
be a fiberwise oriented GMF and let $\g_f$ be the negative
eigenspace bundle of $f$. Then $\g_f$ is an oriented vector bundle
over $\Sig(f)$.

If $\g_f$ is a stably trivial vector bundle then the cellular
chain complex functor $C_\ast(f;\F)$ gives the higher FR torsion
of $(E,\d_0E)$. When $\g_f$ is nontrivial the Framing Principle
(\ref{FP:thm:Framing Principle}) tells us the difference between
$C_\ast(f;\F)$ and the true higher torsion of $(E,\d_0E)$.

The first step in proving the Framing Principle is to stabilize
the bundle $\g_f$ and the complementary bundle $\g_{-f}$.

\begin{lem}[reduction to stable case]\label{Proof:lem:reduction to
stable case} It suffices to prove the Framing Principle in the
case when the fibers $M_t$ have dimension much larger than $\dim
B$ and where the singular set $\Sig(f)$ is a disjoint union of
\emph{framed} and \emph{unframed} components where the framed
components admit a framed structure and the unframed components
have critical points of very large index and coindex.
\end{lem}

\begin{proof}
To do this without introducing technical complications such as
corners, we take the product with $Y=S^{2N}\times S^{2N}$ for $N$
large. Then the product formula (\ref{sec2:lem:product formula})
tells us that
\[
    \t_k(E\times Y,\d_0\times Y;\F')
    =\chi(Y)\t_k(E,\d_0E;\F)=4\t_k(E,\d_0E;\F)
\]
where $\F'$ is the pull back of $\F$ to $E\times Y$. Thus it
suffices to compute the higher torsion of $(E\times Y,\d_0\times
Y)$.

Imitating the proof of the product formula, we choose a fixed
Morse function\[g:Y\to[0,9]\]having exactly four critical points
$y_0,y_1,y_2,y_3$ with critical values $g(y_i)=3i$. Thus $y_0$ is
a minimum, $y_3$ a maximum and $y_1,y_2$ are critical points of
index $2N$.

Let
\[
    h:E\times Y\to[0,10]
\]
be given by $h(x,y)=f(x)+g(y)$. Then the critical points of $h$
will lie over the points $y_i$ since
\[
    \Sig(h)=\Sig(f)\times \Sig(g)\subseteq\coprod E\times y_i.
\]

The critical values of $h_t$ lie in the intervals $(3i,3i+1)$ so
$3i+2=2,5,8$ are regular values of $h_t$ for all $t\in B$. The
sets $S(i)=h^{-1}[3i-1,3i+2]$ for $i=0,1,2,3$ are convex. As in
the proof of the product formula, the corresponding subquotient
functors are
\[
    \xi_i/\xi_{i-1}\cong\Sig^{\ind(y_i)}C_\ast(f;\F).
\]
Since the indices of $y_i$ are all even, the torsion of each of
these subquotients is equal to the torsion of $f$.

By the $C^1$-local framed function theorem
(Remark~\ref{sec1:rem:C1 local FFT}), we can approximate $h|S(i)$
by a fiberwise framed function for each $i$. Let $\widetilde{h}$
be the fiberwise oriented GMF which is equal to the original
function $h$ on $S(1)\cup S(2)$ but equal to the fiberwise framed
$C^1$-approximation on $S(0)$ and $S(3)$. Then $\widetilde{h}$
satisfies the condition of the lemma so we may assume that the
Framing Principle holds for this function.

The unframed components of the singular set $\Sig(\widetilde{h})$
lie in $S(1)$ and $S(2)$. Since the indices of $y_1,y_2$ are even,
the correction terms in the Framing Principle are the same for
both $S(1)$ and $S(2)$ as they are for $f$. Call this term
$C_{FP}$. Then the Framing Principle for $\widetilde{h}$ gives
\begin{equation}\label{Proof:eq:FP for widetilde h}
    4\t_k(E,\d_0E;\F)=\t_k(\widetilde{h};\F')+2C_{FP}.
\end{equation}

By the Splitting Lemma (\ref{sec2:thm:splitting lemma}) we have
\[
    \t_k(\widetilde{h};\F')=\sum_{i=0}^3\t_k(\widetilde{h}|S(i);\F').
\]
We claim that the first and last summands give twice the torsion
of $(E,\d_0E)$. To see this, we can repeat the process above with
$Y=S^{4N}$ and $g:Y\to[0,9]$ having only two critical points with
critical values $0,9$. Then $S(0)=h^{-1}[0,2]$ and
$S(3)=h^{-1}[8,10]$ remain unchanged. So,
\[
    \t_k(\widetilde{h}|S(0);\F')+\t_k(\widetilde{h}|S(3);\F')
    =\t_k(E\times S^{4N},\d_0E\times S^{4N};\F')
    \]
    \[
    =\chi(S^{4N})\t_k(E,\d_0E;\F)=2\t_k(E,\d_0E;\F).
\]
Therefore, subtracting this from both sides of (\ref{Proof:eq:FP
for widetilde h}) we get
\begin{equation*}
    2\t_k(E,\d_0E;\F)=
    \t_k(\widetilde{h}|S(1);\F')+\t_k(\widetilde{h}|S(2);\F')+2C_{FP}.
\end{equation*}

However, on $S(1)$ and $S(2)$, $\widetilde{h}$ is equal to $h$
which is an even suspension of the original function $f$.
Consequently, it has the same torsion and we get
\[
    2\t_k(E,\d_0E;\F)=2\t_k(f;\F)+2C_{FP}.
\]
Divide by $2$ to get the Framing Principle.
\end{proof}

We now assume the conditions of the reduction lemma
(\ref{Proof:lem:reduction to stable case}). The next step is to
prove the stratified deformation lemma.

\begin{lem}[stratified deformation lemma for
FP]\label{Proof:lem:stratified deformation lemma for FP} If we
rewrite the Framing Principle as
\begin{equation}\label{Proof:eq:FP rewritten}
  \t_k(E,\d_0E;\F)-\t_k(f;\F)=C_{FP}
\end{equation}
then both sides of the equation are additive stratified
deformation invariants of the singular pair $(\Sig(f),\c_f)$ where
\[
    \c_f=(\c_f^-,\c_f^+):\Sig(f)\to BSO\times BO
\]
is the classifying map of the stable negative and positive
eigenspace bundles $(\g_f,\g_{-f})$ of $f$. (By \emph{additive} we
mean additive with respect to disjoint union.)
\end{lem}

\begin{rem}[comparison class]\label{Proof:rem:comparison class}
The left hand side of (\ref{Proof:eq:FP rewritten}) is always well
defined even though the individual terms may not be defined. This
difference class is called the \emph{comparison class}. More
details can be found in \cite{[I:BookOne]}.
\end{rem}

Suppose for a moment that the lemma holds. Then, as in the proof
of the transfer theorem (\ref{Proof:thm:transfer theorem}), we
need to verify the following. (Step (a) is unnecessary.)
\begin{enumerate}
  \item[b)] The comparison formula holds if $\c_f$ is trivial
  on the b-d set $\Sig_1$.
  \item[c)] The singular pair $(\Sig,\c_f)$ can be deformed by a
  stratified deformation to a pair which satisfies (b).
\end{enumerate}
However, both of these conditions have already been verified. The
proof of the transfer theorem shows that (c) holds for any
additive stratified deformation invariant. The first statement (b)
is the version of the Framing Principle proved in
\cite{[I:BookOne]} where we assumed that $f$ was framed along
$\Sig_1(f)$. (This is equivalent to $\c_f^-$ being trivial on
$\Sig_1$ under the assumptions of the reduction lemma
\ref{Proof:lem:reduction to stable case}.)
Proposition~\ref{sec2:prop:locally equivalent coefficients} was
used in that proof.

Therefore, the general Framing Principle follows from the
stratified deformation lemma above.

\begin{proof}[Proof of Lemma~\ref{Proof:lem:stratified deformation lemma for
FP}] Since the correction term $C_{FP}$ is a characteristic class
of $\g_f$ on $\Sig$, it is an invariant of the stratified
deformation class of $(\Sig,\c_f)$. So, we consider only the left
hand side of the comparison formula (\ref{Proof:eq:FP rewritten}).

The argument will be the same as the proof of the product formula
(\ref{sec2:lem:product formula}) and the reduction lemma. Assume
that we are given a fiberwise oriented GMF
\[
    f:(E,\d_0E,\d_1E)\to (I,0,1)
\]
satisfying the conditions of the reduction lemma.

Take the product of $E$ with a circle $Y=S^1$. Let
\[
    g:S^1\to[0,6]
\]
be a Morse function with four critical points $y_0,y_1,y_2,y_3$
having critical values $0,2,2+\delta,6$, resp., where $\delta$
will be chosen later. Then these points will be arranged in the
order $y_0,y_2,y_1,y_3$ around the circle as shown in
Figure~\ref{Proof:fig:circle}. Let
\[
    h:E\times S^1\to [0,7]
\]
be given by $h(x,y)=f(x)+g(y)$. Then $h$ will be a fiberwise
oriented GMF. (Take any orientation of $y_2,y_3$.)
\begin{figure}
\includegraphics{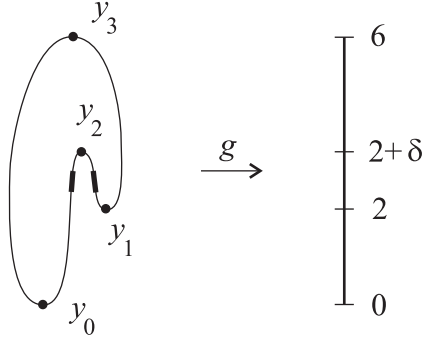}
\caption{Two components of
$g^{-1}[2+\tfrac{\delta}{3},2+\tfrac{2\delta}{3}]\subseteq S^1$
are darkened.} \label{Proof:fig:circle}
\end{figure}

The vertical gradient of $h$ with respect to the product metric is
given by $\grad h_t=(\grad f_t,\grad g_t)$. Consequently,
\[
    \Sig(h)=\Sig(f)\times\{y_0,y_1,y_2,y_3\}.
\]
Since $f$ is a fiberwise oriented GMF, so is $h$. Also, $\frac32$
and $\frac{11}{2}$ are regular values of $h_t$ for all $t\in B$ so
$h^{-1}[\frac32,\frac{11}{2}]$ is a convex set.

At the good points, where $f$ is framed, so is $h$. Of the points
where $h$ is not framed we will modify the ones in $E\times y_1$.
Then we claim that the resulting family of fiberwise oriented
GMF's, $\widetilde{h}$, has the property that
$\widetilde{h}^{-1}[\frac32,\frac{11}{2}]$ is convex with acyclic
subquotient functor whose torsion is the left hand side of
(\ref{Proof:eq:FP rewritten}). We prove this by taking $\delta>1$.

We also claim that the higher torsion of $\widetilde{h}$ is a
stratified deformation invariant of the singular pair
$(\Sig(f),\c_f)$. We prove this by taking $\delta$ very small.
Finally, we show that the subquotient functor of
$\widetilde{h}^{-1}[\frac32,\frac{11}{2}]$ is independent of the
choice of $\delta$. This will prove the lemma.

To prove the first assertion, take $\delta=2$. Then $\frac72$ will
be a regular value of $\widetilde{h}_t$ for each $t\in B$. The
convex set $\widetilde{h}_t^{-1}[\frac72,\frac{11}{2}]$ will have
subquotient functor equal to the suspension of $C_\ast(f)$ since
$y_2$ has index 1. Also the convex set
$\widetilde{h}_t^{-1}[\frac{3}{2},\frac72]$ will have the same
torsion as the original bundle $(E,\d_0E)$ since it is a positive
suspension of $(E,\d_0E)$ with a fiberwise framed function by
construction of $\widetilde{h}$. Since we are assuming that this
torsion is defined, we know, by the Splitting Lemma, that
\begin{align*}
    \t_k(\widetilde{h}|\widetilde{h}^{-1}[\tfrac32,\tfrac{11}{2}];\F)
    &=\t_k(\widetilde{h}|\widetilde{h}^{-1}[\tfrac32,\tfrac{7}{2}];\F)
    +\t_k(\widetilde{h}|\widetilde{h}^{-1}[\tfrac72,\tfrac{11}{2}];\F)\\
    &=\t_k(E,\d_0E;\F)-\t_k(f;\F)
\end{align*}
which is the left hand side of the comparison formula.

Next, we switch the roles of $f$ and $g$. Choose tubular
neighborhoods $L^i$ for the b-d sets $\Sig_1^i(f)$ and other
neighborhoods $T^i$ for the remaining points of $\Sig^i(f)$ as in
Figure~\ref{Proof:fig:bd transfer} in the proof of the first
stratified deformation lemma (\ref{Proof:lem:stratified
deformation lemma}). Take $\delta<\e^3/3$ where $\e$ is the size
of $L^i$. Then $L^i\times g^{-1}[2-\delta,2+2\delta]$ is a special
case of the second example of a convex set (\ref{sec2:eg:second
example}). This includes the deformation of $h$ to
$\widetilde{h}$. We also take tubular neighborhoods $T^i$ of what
is left of $\Sig^i(f)$. Then $T^i(f)\times
g^{-1}[2-\delta,2+2\delta]$ is a special case of the main example
(\ref{sec2:eg:main example}) of a convex set for $h$. The same set
is convex for $\widetilde{h}$, giving is a special case of
Example~\ref{sec2:eg:1st example deformed}.

Let $Z$ be the component of $g^{-1}[2-\delta,2+2\delta]$
containing $y_1$ and $y_2$. Then $Z$ is contractible, so the
subquotient functors given by $L^i\times Z$ and $T^i\times Z$ are
acyclic and thus are locally a direct summand of the expansion
functor for $\widetilde{h}$. However, $\widetilde{h}$ was a
fiberwise framed function given by a small deformation of $h$ in a
neighborhood of $\Sig(f)\times y_1$. Consequently, the subquotient
functor given by the convex sets $L^i\times Z$ and $T^i\times Z$
are stratified deformation invariants. (Since the bundles are
stable, they completely determine the fiber diffeomorphism type of
the convex sets $L^i$ and $T^i$.)

Finally, we will prove that the cellular chain complex of
$\widetilde{h}$ is independent of the choice of $\delta$. In fact
we will show that the filtration of the bundle $E\times S^1$ given
by $\widetilde{h}$ can be made to be independent of the choice of
$\delta$.

First we take $\delta$ small as given in the second case above.
Use the standard metric on $S^1$ for the second case. Next we
increase $\delta$ to $2$ and change $g$ so that the derivative of
$g$ is nonzero only in the set $D$ consisting of the two darkened
region in Figure~\ref{Proof:fig:circle}. Thus $D$ is the union of
two of the components of $g^{-1}[2+\delta/3,2+2\delta/3]$. Call
the new function $\widetilde{g}$. Since the derivative of
$\widetilde{g}$ is nowhere zero on $D$, we can alter the metric on
$S^1$ so that the gradient of $\widetilde{g}$ with respect to the
new metric is equal to the gradient of $g$ with respect to the
standard metric on $S^1$. This has no affect on our assumption
that the metric be standard at all critical points since these are
regular points of $g$ and thus $E\times D$ will have no critical
points of $\widetilde{h}$.

We take the product metric on $E\times S^1$ for the function $h$
and later deform it so that $\widetilde{h}$ will have critical
points with standard metric. Then we deform the metric on $E\times
D$ as indicated above. This creates two functions, call them
$\widetilde{h}$ and $\widetilde{h}'$ which have identical critical
set and identical gradients with respect to their corresponding
metrics. The minimal partial ordering of the critical set and the
twisted cochain which give the expansion functors of
$\widetilde{h}$ and $\widetilde{h}'$ and the subquotient functor
supported on the union of the sets $L^i$ and $T^i$ depend only on
the vertical gradient vector fields $\grad\widetilde{h}$ and
$\grad\widetilde{h}'$. Consequently, the subquotient functors are
the same. This proves the claim, the lemma and the main theorem of
this paper.
\end{proof}

%% file: FPapp.tex
\vfill\eject\section{Applications of the Framing Principle}

\begin{enumerate}
  \item Torelli group
  \item Even dimensional fibers
  \item Unoriented fibers
  \item Vertical normal disk bundle
\end{enumerate}

With the new Framing Principle we can compute the higher FR
torsion when the fiber is closed and even dimensional (and the
higher torsion is defined). This extends the earlier calculation
of the higher torsion for the Torelli group. We also give a
formula for the higher torsion of a bundle with even dimensional
fibers with boundary.

\subsection{Torelli group}\label{FP:subsec:Torelli}

One of the main motivations for the Framing Principle was to study
the higher torsion of the \emph{Torelli group} $T_g$. This is the
group of isotopy classes of self-diffeomorphisms $\f$ of a closed
oriented surface $\Sig_g$ of genus $g$ so that $\f$ induces the
identity on the integral homology of $\Sig_g$.

It was John Klein \cite{[K:Torelli]} who first realized that the
higher torsion invariants $\t_{2k}(T_g)$ for the Torelli group
could be defined and he conjectured that they were proportional to
the {Miller-Morita-Mumford class} $\k_{2k}$.

Recall that the \emph{Miller-Morita-Mumford classes}
(\cite{[Mumford83:MMM_class]},\cite{[Morita87]},\cite{[Miller86:MMM]})
for closed oriented surface bundle $E\to B$ are the integral
cohomology classes in $B$ given by pushing down the powers of the
Euler class of the vertical tangent bundle of $E$.
\begin{equation}\label{FP:eq:MMM classes}
    \k_k(E):=p^E_\ast(e_E^{k+1})\in H^{2k}(B;\ZZ).
\end{equation}
For surfaces with boundary we use the Euler class $e_E\in H^2(E,\d
E)$ to get the \emph{punctured Miller-Morita-Mumford classes}
\begin{equation}\label{FP:eq:punctured MMM classes}
  \k_k^\d(E)=p_\ast^E(e_E^{k+1})=p_\ast^E(c_1^k\cup
  e_E)=tr^E_B(c_1^k)
\end{equation}
where $c_1\in H^2(E)$ is the restriction of $e_E$ to $E$.

Klein's conjecture was first proved by R. Hain, R. Penner and the
author (unpublished) using B\"{o}kstedt's theorem
\cite{[Bokstedt84]}. In \cite{[I:BookOne]}, \emph{fat graphs}
(finite graphs with cyclic orderings of the half-edges incident to
each vertex) and the Framing Principle were used to find the
proportionality constant. In order to use fat graphs, we needed to
restrict to surfaces with at least one marked point.

\begin{thm}[torsion of $T_g^s$]\label{FP:thm:torsion of Tgs}
Let $T_g^s$ be the Torelli group of genus $g$ surfaces with
$s\geq1$ marked points. Then
\[
    \t_{2k}(T_g^s)=\frac12(-1)^k\z(2k+1)\frac{\k_{2k}}{(2k)!}.
\]
\end{thm}

The version of the Framing Principle proved in \cite{[I:BookOne]}
assumes that the fiberwise GMF $f:E\to\RR$ is framed along the
birth-death set. So, in the proof of \ref{FP:thm:torsion of Tgs},
we needed to choose $f$ to have maxima at the marked points and
minima and saddle points on the fat graph which forms the
canonical core of the complement of the set of marked points.

With the new Framing Principle we can prove this more directly
without using fat graphs.

\begin{thm}[torsion of oriented surface bundles]Let $\Sig\to E\to B$ be a closed oriented surface bundle so
that the action of $\pi_1B$ on $H_\ast(\Sig,\QQ)$ is upper
triangular. Then
\[
    \t_{2k}(E)=\frac12(-1)^k\z(2k+1)\frac{\k_{2k}(E)}{(2k)!}.
\]
\end{thm}

\begin{proof}
This is a special case of Theorem~\ref{FP:thm:even dim closed
fibers} below.
\end{proof}

\subsection{Even dimensional fibers}\label{FP:subsec:even fibers}

The calculation of the higher torsion for surface bundles extends
to all closed oriented even dimensional fibers. However, in order
to state the result we need a generalization of the
Miller-Morita-Mumford classes.

\begin{defn}[the integral class $T_{2k}(E)$]\label{FP:def:general real MMM class}
For any smooth bundle $E\to B$ with compact manifold fibers, let
\[
    T_{2k}(E)=tr^E_B\left(\frac{(2k)!}{2}ch_{2k}(T^vE\otimes\CC)\right)\in
H^{4k}(B;\ZZ)
\]
This is an integral class by Remark~\ref{FP:rem:Newton
polynomial}.
\end{defn}

If $T^vE$ already has a complex structure (as in the case of
oriented surface bundles), we have
\begin{equation}\label{FP:eq:T2k(E)=k2k(E)}
  \frac{(2k)!}{2}ch_{2k}(T^vE\otimes\CC)=(2k)!ch_{2k}(T^vE)=\sum
    c_1(\ll_i)^{2k}.
\end{equation}

\begin{prop}[$T_{2k}$ for surface bundles]\label{FP:prop:T=k}
Let $\Sig\to E\to B$ be an oriented surface bundle. Then
\begin{enumerate}
  \item[a)] $T_{2k}(E)=\k_{2k}(E)$ if $\Sig$ is closed.
  \item[b)] $T_{2k}(E)=\k_{2k}^\d(E)$ if $\Sig$ has a boundary.
  \item[c)]$T_{2k+1}(E)=0$ in both cases.
\end{enumerate}
\end{prop}

\begin{prop}[$T_{2k}$ for $M^{2n+1}$ oriented]\label{FP:prop:tangential invariant for odd dim M}
If the fiber $M$ is an oriented odd dimensional manifold then
\[
    2T_{2k}(E)=T_{2k}(\d E).
\]
\end{prop}

\begin{proof}
First suppose that $M$ is a closed manifold of dimension $2n+1$.
Let $f:E\to\RR$ be any generic smooth function. Then the fiberwise
singular set $\Sig(f)$ is a normally oriented codimension $2n+1$
submanifold of $E$ dual to the Euler class of $T^vE$. But,
$\Sig(-f)=\Sig(f)$ with the normal orientation reversed. So
\[
    e_E=-e_E=0.
\]

Now suppose that $M$ has a boundary. Then $\grad(-f)$ points the
wrong way along $\d E$ so we have to ``fold up'' a collar
neighborhood of $\d E$ by adding the function $g(x)+t^2$ where
$g:\d E\to\RR$ is a smooth function and $t$ is negative the
distance to $\d E$. (See the proof of
Corollary~\ref{sec2:cor:boundary torsion and relative torsion} for
more details.) This adds the fiberwise singular set $\Sig(g)$ to
$\Sig(-f)$. So
\[
    e_E=-e_E+e_{\d E}.
\]
The sign in front of $e_{\d E}$ is the sign of the second
derivative of $t^2$.\
\end{proof}

Let $E^\ast=E\times D^N$ with the corners rounded. Then by the
definition of transfer we have
\begin{equation}\label{FP:eq:T2k is stable}
    T_{2k}(E^\ast)=T_{2k}(E).
\end{equation}

\begin{thm}[torsion for $M^{2n}$ closed, oriented]\label{FP:thm:even dim closed fibers}
Suppose that $M^{2n}$ is a closed oriented even dimensional
manifold and $M\to E\to B$ is a smooth bundle where $\pi_1B$ acts
upper triangularly on the rational homology of $M$. (In
particular, it preserves the orientation of $M$.) Then the higher
FR-torsion invariants $\t_{2k}(E)\in H^{4k}(B;\RR)$ are defined
and given by
\[
    \t_{2k}(E)=\frac12
    (-1)^k\z(2k+1)\frac1{(2k)!}T_{2k}(E).
\]
\end{thm}

\begin{proof}
Suppose first that there is a fiberwise oriented GMF $f:E\to \RR$.
Then $-f$ is also a fiberwise oriented GMF on $E$ so the Framing
Principle gives us:
\[
    \t_{2k}(E)=\t_{2k}(C_\ast(f))+
    (-1)^k\z(2k+1)p_\ast
    \left(\tfrac12 ch_{2k}(\g_f\otimes\CC)\right)
\]
\[
    \t_{2k}(E)=\t_{2k}(C_\ast(-f))+(-1)^k\z(2k+1)p_\ast
    \left(\tfrac12 ch_{2k}(\g_{-f}\otimes\CC)\right)
\]
By the involution property,
$\t_{2k}(C_\ast(-f))=-\t_{2k}(C_\ast(f))$. Consequently, the sum
of the above equations is:
\[
    2\t_{2k}(E)=(-1)^k\z(2k+1)p_\ast
    \left(\tfrac12 ch_{2k}((\g_f\oplus\g_{-f})\otimes\CC)\right).
\]
The theorem then follows from Proposition~\ref{FP:prop:gamma f +
gamma -f give nu Sigma} which implies that
\[
    \g_f\oplus\g_{-f}\cong T^vE|\Sig(f)\oplus\e^{2n}
\]
and the push-down/transfer diagram (
Proposition~\ref{FP:prop:push-down transfer diagram}) which says
that
\[
    p_\ast(c|\Sig(f))=tr^E_B( c).
\]

Now suppose that there is no fiberwise oriented GMF on $E$. Then
we need to stabilize by taking the fiberwise product with a large
even dimensional sphere
\[
    S^{2N}=D_-^{2N}\cup D_+^{2N}.
\]
By the framed function theorem there is a fiberwise framed
function
\[
    f:E\times D_-^{2N}\to (-1,0]
\]
taking the boundary to $0$. By definition, the higher FR torsion
of $C_\ast(f)$ is equal to the higher torsion of $E$:
\[
    \t_{2k}(E)=\t_{2k}(C_\ast(f)).
\]
Taking $-f$ on the upper hemisphere we get a fiberwise oriented
GMF
\[
    f\cup-f:E\times S^{2N}\to\RR.
\]
By the Splitting Lemma and the Suspension Theorem, the higher
torsion of $E\times S^{2N}$ is equal to twice the higher torsion
of $E$:
\[
    \t_{2k}(E\times
    S^{2N})=\t_{2k}(E\times D_-^{2N})+\t_{2k}(E\times D_+^{2N},\d)=(1+(-1)^{2N})\t_{2k}(E).
\]
But the Framing Principle tells us that:
\[
    \t_{2k}(E\times S^{2N})=\t_{2k}(C_\ast(f\cup-f))+
        (-1)^k\z(2k+1)p_\ast
    \left(\tfrac12 ch_{2k}(\g_{-f}\otimes\CC)\right).
\]
The first term is zero by the Splitting Lemma and involution
property:
\[
    \t_{2k}(C_\ast(f\cup-f))=\t_{2k}(C_\ast(f))+\t_{2k}(C_\ast(-f))=0.
\]
In the second term, $\g_{-f}$ is isomorphic to the restriction to
$\Sig(-f)$ of the vertical tangent bundle of $E^\ast=E\times
D_+^{2N}$. Consequently,
\[
    p_\ast
    \left(\tfrac12 ch_{2k}(\g_{-f}\otimes\CC)\right)=
    tr^{E^\ast}_B\left(\tfrac12 ch_{2k}(T^vE^\ast\otimes\CC)
    \right)=\frac1{(2k)!}T_{2k}(E)
\]
by (\ref{FP:eq:T2k is stable}). Now divide both sides by $2$.
\[
    \t_{2k}(E)=\frac12\t_{2k}(E\times S^{2N})
    =\frac12
    (-1)^k\z(2k+1)\frac1{(2k)!}T_{2k}(E).
\]
\end{proof}

Now suppose that the fiber $M$ is even dimensional and oriented
with boundary.

\begin{thm}[torsion for $M^{2n}$ oriented with boundary]\label{FP:thm:even dim oriented fiber with boundary}
Suppose that $M^{2n}$ is an oriented even dimensional manifold and
$M\to E\to B$ is a smooth bundle where the action of $\pi_1B$ on
the rational homology of both $M$ and $\d M$ is upper triangular.
Suppose also that $\pi_1B$ preserves the orientation of $M$. Then
the higher FR-torsion invariants $\t_{2k}(E)\in H^{4k}(B;\RR)$ are
defined and given by
\[
    \t_{2k}(E)=\underbrace{\frac12\t_{2k}(\d E)}_{exotic\ component}+
    \underbrace{\frac12
    (-1)^k\z(2k+1)\frac1{(2k)!}T_{2k}(E)}_{tangential\ component}.
\]
\end{thm}

\begin{rem}[exotic and tangential torsion]\label{FP:rem:exotic and tangential components for even
fibers} The two terms in the above formula can be viewed as the
\emph{exotic} and \emph{tangential} components of the higher FR
torsion. The tangential component is a tangential homotopy
invariant. The exotic component always comes from a bundle with
odd dimensional closed fiber.
\end{rem}

\begin{proof}
Suppose first that there is a fiberwise oriented GMF $f:E\to\RR$.
Then the Framing Principle for $f$ and $-f$ gives us:
\[
    \t_{2k}(E)=\t_{2k}(C_\ast(f))+(-1)^k\z(2k+1)p_\ast
    \left(\tfrac12 ch_{2k}(\g_f\otimes\CC)\right)
\]
\[
    \t_{2k}(E,\d E)=\t_{2k}(C_\ast(-f))+(-1)^k\z(2k+1)p_\ast
    \left(\tfrac12 ch_{2k}(\g_{-f}\otimes\CC)\right).
\]
By the Splitting Lemma we always have:
\[
    \t_{2k}(E)=\t_{2k}(E,\d E)+\t_{2k}(\d E).
\]
Adding these three equations we get:
\[
    2\t_{2k}(E)=\t_{2k}(\d E)+(-1)^k\z(2k+1)p_\ast
    \left(\tfrac12 ch_{2k}(T^vE\otimes\CC)\right)
\]
and we are done. (See the proof of Theorem~\ref{FP:thm:even dim
closed fibers}.)

If the function $f$ does not exist we need to take the product
with a high even dimensional sphere $S^{2N}$. This double all
three terms in our theorem.
\end{proof}

\subsection{Unoriented fibers}\label{FP:subsec:unoriented fibers}

If the fiber $M$ is unoriented we can pass to an oriented disk
bundle over $E$. Here we take an oriented $1$-disk bundle.

Every unoriented manifold $M$ has a $2$-fold oriented covering
$\widetilde{M}$ whose rational homology is equal to the sum
\[
    H_\ast(\widetilde{M};\QQ)\cong H_\ast(M;\QQ)\oplus
    H_\ast(M;\QQ^-)
\]
where $\QQ^-$ is a suitable twisted coefficient system on $M$.
This formula also holds with $\QQ$ replaced by $\RR$ or $\CC$. If
$M$ is closed, the covering $\widetilde{M}$ bounds a
$1$-dimensional disk bundle over $M$ which we call $JM$. Applying
this to each fiber of a smooth bundle $E\to B$ we get a $2$-fold
covering $\widetilde{E}$ of $E$ which forms the boundary of a disk
bundle $JE$.

\begin{lem}[comparing $E$ and $JE$]\label{FP:lem:comparing E and JE}
{\rm a)} $T_{2k}(E)=T_{2k}(JE)$.

{\rm b)} $\t_{2k}(E)=\t_{2k}(JE)$ if it is defined.

{\rm c)} If $\t_{2k}(E)$ and $\t_{2k}(E;\CC^-)$ are both defined
then
\[
    \t_{2k}(\widetilde{E})=\t_{2k}(E)+\t_{2k}(E;\CC^-).
\]
\end{lem}

\begin{proof} (a) follows from the fact that the Pontrjagin
classes of a real line bundle are all zero. (b) is obvious. (c) is
a special case of the transfer formula
Theorem~\ref{sec2:thm:transfer formula for groups} since
$\Ind_1^{Z_2}\CC=\CC\oplus\CC^-$.
\end{proof}

Suppose that $M$ is odd dimensional and closed. Then $JE$ is a
smooth bundle with oriented even dimensional fibers. Consequently,
its higher torsion is given by Theorem~\ref{FP:thm:even dim
oriented fiber with boundary}. Using Lemma~\ref{FP:lem:comparing E
and JE} above, this calculation can be expressed as follows.

\begin{thm}[torsion for $M^{2n+1}$ closed,
unoriented]\label{FP:thm:closed unoriented odd fibers} If
$\t_{2k}(E;\CC)$ and $\t_{2k}(E;\CC^-)$ are both defined then
\[
    \t_{2k}(E;\CC^\pm)=\frac12\t_{2k}(\widetilde{E})\pm
    \frac12(-1)^k\z(2k+1)\frac1{(2k)!}T_{2k}(E).
\]
\end{thm}

If $M$ is even dimensional, closed and unoriented, the best we can
do is the following.

\begin{thm}[torsion for $M^{2n}$ closed,
unoriented]\label{FP:thm:closed unoriented even M}
\[
    \t_{2k}(E)+\t_{2k}(E;\CC^-)=(-1)^k\z(2k+1)\frac1{(2k)!}T_{2k}(E).
\]
\end{thm}

\begin{proof}
By Lemma~\ref{FP:lem:comparing E and JE}(c) and the torsion for
oriented even dimensional fibers (Theorem~\ref{FP:thm:even dim
closed fibers}) we have:
\[
    \t_{2k}(E)+\t_{2k}(E;\CC^-)=\t_{2k}(\widetilde{E})
    =\frac12(-1)^k\z(2k+1)\frac1{(2k)!}T_{2k}(\widetilde{E}),
\]
where
\[
    T_{2k}(\widetilde{E})=2T_{2k}(JE)=2T_{2k}(E)
\]
by Proposition~\ref{FP:prop:tangential invariant for odd dim M}
applied to $JE$ and Lemma~\ref{FP:lem:comparing E and JE}(a).
\end{proof}

If $M$ has a boundary then $\widetilde{E}$ is still a $2$-fold
covering of $E$ but it forms only part of the boundary of $JE$.
Let $J^\ast E$ be the disk bundle $JE$ with the corners rounded.
Then the boundary of $J^\ast E$ is an oriented bundle over $B$
\[
    \d J^\ast E=\widetilde{E}\cup J\d E
\]
with oriented closed manifold fibers $\d J^\ast
M=\widetilde{M}\cup J\d M$. As before, we have
\[
    \t_{2k}(E)=\t_{2k}(J^\ast E)
\]
if the terms are defined. If $M$ is odd dimensional then $J^\ast
E$ is oriented with even dimensional fibers.

\begin{thm}[torsion for $M^{2n+1}$ unoriented with boundary]\label{FP:thm:unoriented fibers with boundary}
If the fiber $M$ of $p:E\to B$ is odd dimensional then
\[
    \t_{2k}(E)=\frac12\t_{2k}(\widetilde{E}\cup J\d
    E)+\frac12(-1)^k\z(2k+1)\frac1{(2k)!}T_{2k}(J^\ast E)
\]
when all terms are defined.
\end{thm}

For the sake of completeness, we make the following trivial
observation.

\begin{prop}[torsion for $M^{2n}$ unoriented]\label{FP:prop:unoriented even dim M} If $M$ is even
dimensional and unoriented then
\[
    \t_{2k}(E)=\t_{2k}(JE)
\]
assuming the terms are defined.
\end{prop}

\begin{rem}\label{FP:rem:about unoriented fibers}
In all cases, the calculation of the higher torsion is reduced to
the case of an oriented odd dimensional fiber (plus the tangential
homotopy invariant $T_{2k}$).
\end{rem}

\subsection{Vertical normal disk bundle}\label{FP:subsec:normal disk
bundle}

We need to generalize this to a higher dimensional disk bundles.
Since $E$ is compact, we can always find a fiberwise embedding of
$E$ into $B\times\RR^{2N}$ for $N$ sufficiently large. Choose such
an embedding and let $D(E)$ be the unit normal disk bundle of $E$
in $B\times\RR^{2N}$. Let $S(E)$ be the unit normal sphere bundle.
We call $D(E), S(E)$ the \emph{vertical normal disk bundle} and
the \emph{vertical normal sphere bundle}, respectively.

The vertical tangent bundle of $D(E)$ is trivial by construction.
Consequently,
\[
    T_{2k}(D(E))=0.
\]
Since $D(E)$ has even dimensional fibers, Theorem~\ref{FP:thm:even
dim oriented fiber with boundary} gives us the following.

\begin{thm}[torsion of $D(E)$ for $M$ closed, oriented]\label{FP:thm:normal disk bundle-closed case}
Suppose that the fiber $M$ is a closed oriented manifold and the
action of $\pi_1B$ on the rational homology of $M$ is upper
triangular. Then
\[
    \t_{2k}(E)=\t_{2k}(D(E))=\frac12\t_{2k}(S(E)).
\]
\end{thm}

If $M$ has a boundary then
\[
    \d D(E)=D(\d E)\cup S(E)
\]
and these two subsets meet at an angle along the set
\[
    D(\d E)\cap S(E)=S(\d E)
\]
whose vertical tangent bundle is stably parallelizable. Since the
fibers of $S(\d E)$ are closed and even dimensional its higher
torsion is trivial if it is defined.

\begin{thm}[torsion of $D(E)$ for $M$ oriented with boundary]\label{FP:thm:normal disk bundle with boundary}
Suppose that $M$ is oriented with boundary and the action of
$\pi_1B$ on the rational homology of both $M$ and $\d M$ is upper
triangular. Then
\[
    \t_{2k}(E)=\frac12\t_{2k}(\d E)+\frac12\t_{2k}(S(E)).
\]
\end{thm}

\begin{proof}Let $D^\ast(E)$ be $D(E)$ with the corners rounded.
Then
\begin{enumerate}
  \item $\t_\ast(E)=\t_\ast(D^\ast(E))$ since $D(E)$ is a linear
  disk bundle over $E$ (\ref{sec1:lem:positive
suspension lemma}).
  \item $\t_\ast(D^\ast(E))=\frac12\t_\ast(\d D^\ast(E))$ since
  $D^\ast(E)$ has even dimensional parallelizable fibers (\ref{FP:thm:even dim oriented fiber with
  boundary}).
  \item Since $\d D^\ast(E)=D(E)|\d E\cup S(E)$, the gluing formula (\ref{sec2:cor:gluing
formula}) gives
\[
    \t_\ast(\d D^\ast(E))=\t_\ast(D(E)|\d
    E)+\t_\ast(S(E))-\t_\ast(S(E)|\d E).
\]
    \item $\t_\ast(D(E)|\d E)=\t_\ast(\d E)$, again by (\ref{sec1:lem:positive
suspension lemma}).
    \item $\t_\ast(S(E)|\d E)=0$ since $S(E)|\d E$ is a bundle
    with stably parallelizable closed even dimensional fibers
    (\ref{FP:thm:even dim closed fibers}).
\end{enumerate}

\end{proof}

If we compare this with Theorem~\ref{FP:thm:even dim oriented
fiber with boundary} we get:

\begin{cor}[torsion of $S(E)$ for $M^{2n}$ oriented]\label{FP:cor:torsion of S(E)}
Suppose that $M$ is oriented and even dimensional and that the
action of $\pi_1B$ on the rational homology of both $M$ and $\d M$
is upper triangular. Then
\[
    \t_{2k}(S(E))=
    (-1)^k\z(2k+1)\frac1{(2k)!}T_{2k}(E).
\]
\end{cor}

\begin{rem}[generalizing the transfer theorem]
Examination of this corollary reveals that it is a transfer
formula. It says that the higher torsion of $S(E)$ over $B$ is the
transfer of the higher torsion of $S(E)$ over $E$. However, it is
not an example of the transfer theorem since the fibers of $S(E)$
over $E$ (odd dimensional spheres) are not acyclic. This suggests
that the transfer theorem may hold whenever the euler
characteristic of the fiber is zero. (It does not hold otherwise.)
\end{rem}